\documentclass{article}

\usepackage{amssymb,amsmath}
\usepackage{amsthm}
\usepackage{layout}

\usepackage{amscd}
\usepackage[all]{xy}
\usepackage[percent]{overpic}

\usepackage{hyperref}

\def\draft{1}

\theoremstyle{plain}
\newtheorem{theorem}{Theorem}[subsection]
\newtheorem*{theorem*}{Theorem}
\newtheorem{proposition}[theorem]{Proposition}
\newtheorem{lemma}[theorem]{Lemma}
\newtheorem{corollary}[theorem]{Corollary}

\newtheorem*{definition}{Definition}
\newtheorem{example}{Example}
\newtheorem*{example*}{Example}
\newtheorem{nonexample}{Non-Example}

\newtheorem*{lemma*}{Lemma}
\newtheorem*{question*}{Question}

\newcommand{\rmk}[1]{\noindent \it Remark:\rm {#1}}

\numberwithin{equation}{section}
\allowdisplaybreaks[1]

\title{Collapsing Manifolds with Boundary}
\author{Jeremy Wong}
\date{}

\begin{document}

\maketitle

\ifnum \draft=0
   \frontmatter
\fi

\begin{abstract}

This manuscript
studies manifolds-with-boundary collapsing in the Gromov-Hausdorff topology. The main aim
is an understanding of the relationship of the topology and geometry of a limiting sequence of manifolds-with-boundary to that of a limit space, which is presumed to be without geodesic terminals.

The main result establishes a disc bundle structure for any manifold-with-boundary having two-sided bounds on sectional curvature and second fundamental form, and a lower bound on intrinsic injectivity radius, which is sufficiently close in the Gromov-Hausdorff topology to a closed manifold.

The second main result identifies Gromov-Hausdorff limits of certain sequences of manifolds-with-boundary as Alexandrov spaces of curvature bounded below.

\end{abstract}

\normalsize

\tableofcontents

\newpage

\section{Introduction} \label{s:0}

\noindent

This paper studies manifolds-with-boundary collapsing in the Gromov-Hausdorff topology.
The only previously published work dealing specifically with Gromov-Hausdorff limits of manifolds-with-boundary seems to be~\cite{Kod:90}, which considered only the non-collapsing regime.

For a manifold-with-boundary there are qualitatively several types of collapses, which may be made precise by stipulating that certain injectivity radii tend to zero:

(1) interior collapse (e.g., rounded cone $\longrightarrow$ cone)

(2) boundary collapse (e.g., $S^{n}(1)\setminus B^{n}(\frac{1}{i}) \longrightarrow S^{n}(1)$)

(3) boundary contact (e.g., $D^{n}(0,2)\setminus B^{n}(1-\frac{1}{i},1) \longrightarrow D^{n}(0,2)\setminus B^{n}(1,1)$, where $D^{n}(x,r)$ (resp. $B^{n}(x,r)$) denotes the closed (resp. open) disc in $\mathbb{R}^{n}$ with center $x$ and radius $r$.)\\

These collapsings can be partial (not occurring globally), may occur on different scales, and may happen simultaneously. (1) may be viewed as a local version of the situation in the closed case, and is somewhat, though far from being completely, well understood. (2), it seems, bears consequences for the interior in direct proportion to the degree of control of the embedding of the boundary (e.g. second fundamental form). With enough of the latter, then, the collapse of the boundary itself may also be viewed and understood from the perspective of the closed case. So it remains to focus on collapses of the form (3).

Manifolds-with-boundary are distinguished from closed manifolds in that there is at almost every point a preferred direction, namely, the shortest route to the boundary. Under the right curvature conditions (see, e.g., Theorem~\ref{thm:Yam}), this feature allows a factorization of a collapse, analogously as in the case of a collapsing sequence of closed manifolds in situations when one is able to scale the fiber in certain directions.
Thereby for a given collapsing sequence of manifolds-with-boundary, (3) may be studied largely independently of (1) and (2).\\

Many (though not all) results concerning existence of limits, regularity of limit metrics, finiteness theorems, etc. for the closed manifold case have analogues in the bordered case, provided supplementary hypotheses are made on the boundary, such as bounds on the second fundamental form, and certain injectivity radii.\\

Several results substantiate the statement that collapses of type (3) do not generate additional topology: (i) the
extension procedure introduced in \cite{Wo:07} (see $\S$\ref{sss:1.2.1})
together with the attendant
homeomorphism finiteness theorem and homotopy characterization of limits,
and (ii)
Theorems~\ref{thm:4.1}-\ref{thm:4.2} and \ref{thm:Top}
of the present paper,
which yield disc bundle structures for a sequence of manifolds collapsing in the sense of (3) under a lower sectional curvature bound and having not-too-concave boundary.\\

Besides their intrinsic interest,
another motivation for studying collapses of manifolds-with-boundary involves understanding geodesic terminals.
A geodesic terminal is a type of singularity.
For a sequence of closed manifolds, a limit fails to be geodesically extendible only if the sectional curvatures are not bounded above
or there is an injectivity radius collapse.
Conversely, starting from a sequence of spaces with singularities, one would like to determine how these singularities may disappear in a limit.\\

For Riemannian manifolds-with-boundary, the boundary itself comprises the set of geodesic terminals.
The boundary is also a so-called extremal subset. In the setting of Alexandrov spaces of curvature bounded below, it was conjectured in~\cite{P:97} that if a limit of Alexandrov spaces had no proper extremal subsets then the collapsing spaces would be fiber bundles over the limit.
Theorem~\ref{thm:4.2} verifies that the conclusion of this conjecture holds, for a special class of manifolds-with-boundary which are themselves not necessarily Alexandrov space of curvature bounded below, but in a sense (made precise in $\S$\ref{ss:1.2}) topologically close to them.\\

\newpage
\noindent
\bf Methods \rm

\vspace{0.2cm}

The chief approach taken in this paper is to study manifolds-with-boundary extrinsically, by extending their boundaries via a codimension-zero extension. This extension is obtained by a non-smooth gluing procedure, though in principle one could just as well try to obtain a $C^{\infty}$ extension by solving a Cauchy initial value problem with the terminal boundary constraint that the new extended boundary be totally geodesic, or at least, locally convex.

The advantage of the exterior approach (compared with 'interior' approaches which would analyze the distance function to the boundary via Morse theory) is that now one can allow the interiors of a sequence of manifolds-with-boundary to collapse in the sense of (3) above.\\

\vspace{0.1cm}

\noindent
\bf Assumptions \rm

\vspace{0.2cm}

Now we shall indicate which curvature bounds will be typically assumed in the sequel, and why.\\

A given smooth Riemannian manifold-with-boundary is an Alexandrov space of curvature bounded above~\cite{ABB:93}, but will not be an Alexandrov space of curvature bounded below unless the boundary is locally convex. For locally convex boundary and non-negative Ricci curvature in the interior (or $curv\geq 0$ for Alexandrov spaces, and then automatically locally convex boundary) there are topological recognition theorems and splitting theorems which identify the manifold isometrically as a cylinder or warped product. Local convexity of the boundary is a strong condition which, if assumed, renders (3) less interesting. In examining (3) some lower bounds on the sectional curvature and second fundamental form are natural. However, to say more about the structure of a limiting sequence in relation to the geometry and topology of its limit, additional hypotheses should be made, about the limit, for instance. The main hypothesis which will be invoked is that the limit is geodesically extendible.\\

\vspace{0.1cm}

\noindent
\bf Notations and Conventions \rm 

\vspace{0.2cm}

Manifolds are assumed to be
metrically complete,
unless specified otherwise.

For an immersion $N\hookrightarrow M$, an inequality of the form $II\geq \lambda$ signifies that all eigenvalues of the associated quadratic form $S: TN\longrightarrow TN$ are $\geq \lambda$. Here $II(X,Y)\overset{def}{=}g({\nabla}_{X}\nu,Y)$, where $\nu$ is the outer normal. By convention, the standard flat disc $D^{2}(r)$ of radius $r$ has $II= \frac{1}{r}\geq 0$ and convex boundary. Immersions with $II\geq 0$ will be called convex, those with $II\leq 0$, concave.

If $N$ is a disconnected Riemannian manifold, an inequality of the form $d(N)\leq d$ will usually be interpreted to mean that every path component of $N$ has an upper intrinsic diameter bound $d$.\\

$\simeq$ $\ \ $ denotes homotopy equivalence, and $\approx$ denotes either homeomorphism or diffeomorphism.

$[xy]_{X}$ $\ \ $ for a length space $X$, denotes a minimizing geodesic segment from the point $x$ to the point $y$.
A geodesic segment is a geodesic which globally realizes the distance between any two of its points. The terms 'minimal geodesic' and 'minimizing geodesic' will also be used synonymously for this.

$B(x,r; X)$ $\ \ $ denotes an open metric ball in $X$ of radius $r$ centered at $x$.

$D^k$ $\ \ $ denotes a closed ball (metric or otherwise) of dimension $k$.

$d_{X}, \, d$ or $|\cdot|_{X}$ interchangeably denote the metric distance function of a metric space $X$.

$M_{K}^{2}$ $\ \ $ denotes the standard two-dimensional, simply-connected model space of constant curvature $K$. 

$M[0,r]$ $\ \ $ for a manifold-with-boundary $M$, denotes the set $\{ x\in M : d(x,\partial M) \leq r \}$.

$\mathcal{M}(n,K^{-},{\lambda}^{\pm},inj(M),i_{\partial},d)$, $\ \ $ for instance, denotes the class of $n$-dimensional manifolds-with-boundary $M$ having lower interior sectional curvature bound $K^{-}$, lower $({\lambda}^{-})$ and upper $({\lambda}^{+})$ bound on the second fundamental form, some (unspecified) \emph{uniform} positive lower bounds to $inj(M)$ and $i_{\partial}(M)$, and an upper diameter bound $d$. A two-sided curvature bound such as occurs in the notation $\mathcal{M}(n,K^{\pm},{\lambda}^{-})$ may occasionally be abbreviated as $\mathcal{M}(n,K,{\lambda}^{-})$.

$\tau(\ldots)$ $\ \ $ indicates a positive constant which tends to $0$ as its arguments in parentheses tend to $0$.\label{indexl1}\\

\newpage

\ifnum \draft=0
  \mainmatter
\fi

\section[Disc Bundle Structures and Metric Structure of Limits]{Disc Bundle Structures and Metric Structure of Limits }
\label{s:2}

\subsection{Statement of Theorems}
\label{ss:2.1}

The structure of nonnegatively-curved open manifolds is described by the celebrated
Soul Theorem of Cheeger and Gromoll, which reduces the study of such manifolds to closed manifolds of nonnegative curvature. Effectively, the topology of the open manifold is entirely contained in the
soul.
After Cheeger and Gromoll's reduction using Busemann functions to initially construct a totally convex subset, the
Soul 
Theorem may be stated in terms of manifolds with locally convex boundary.\\

\begin{theorem}\label{thm:4.1} \bf
Soul
Theorem, \cite{CG:72}).  \rm
Let $(M,\partial M)$ be a complete Riemannian manifold.\\
If $0\leq K_M$, $0\leq II_{\partial M}$ then $M$ is diffeomorphic to a $D^{k}$-bundle (the normal bundle) over a closed totally convex submanifold.
\end{theorem}
\vspace{0.5cm}

Their method of proof consisted of constructing a nested sequence of totally convex subsets, each dropping in dimension, until one arrived at the claimed closed, totally convex submanifold. These subsets are submanifolds (with boundary) of the original manifold under consideration.
 
One would like to see how the hypotheses of nonnegative curvature and locally convex boundary can be weakened, yet the disc-bundle structure conclusion retained, if perhaps some additional hypotheses are a priori made about a supposed limit space.
The
Soul 
Theorem 
may be viewed
not essentially as a theorem about non-negatively curved manifolds, but rather, as its proof via maximal inward equidistant retractions shows, a theorem about the disappearance of boundary. For the
soul is without boundary.

It is possible to go beyond an almost-nonnegative, almost-convex boundary setting. To do this, it is necessary to rely less on convexity and more on subspace methods.\\\\

The following theorems identify manifolds-with-boundary from certain classes as disc bundles over their limit.
In the following, $i$ will be taken sufficiently large, depending on the indicated parameters. \\

\begin{theorem}\label{thm:Yam}
Let $(M_i,\partial M_i)$ be a sequence of Riemannian manifolds-with-boundary, and assume that\\ $M_i \overset{GH}{\longrightarrow} N$, where $N$ is a closed manifold.\\
(i) If $0\leq K_{M_i}$, \hspace{0.7cm} $0\leq II_{\partial M_i}$, then
$F\longrightarrow M_i \longrightarrow N$ is a locally trivial fiber bundle, where the fiber $F$ is
a $D^{k}$-bundle over a closed manifold ($1 \leq k \leq n$),
almost non-negatively curved in the generalized sense, as defined in \cite[section 5]{Y:91}\\
(ii)
If $K^{-}\leq K_{M_i}$, \hspace{0.1cm} $0\leq II_{\partial M_i}$, then $F\longrightarrow M_i \longrightarrow N$ is a locally trivial fiber bundle, where $F$ is a manifold-with-boundary, almost non-negatively curved in the generalized sense.
\end{theorem}

Part (ii) is a direct consequence of
a fibering theorem of Yamaguchi for Alexandrov spaces of curvature bounded below (Theorem~\ref{thm:4.4} given in a following section),

However, part (i) cannot be shown directly with the arguments of Yamaguchi's original fibering theorem~\cite{Y:91}, which assumed geodesic extendibility of the approximating spaces. It potentially could be a direct consequence of
\cite{Y:96}. But it seems that a direct application of \cite{Y:96} would only be able to identify the fiber as a manifold-with-boundary having non-negative curvature in the generalized sense, as defined in
\cite[section 5]{Y:91}.
This is due to the nature of the map produced, which is an almost-Riemannian submersion, instead of a Riemannian submersion.
\\

\vspace{0.3cm}

At least some lower bound $II \geq {\lambda}^{-}$ is necessary in
Theorem~\ref{thm:Yam}
for there to exist a locally trivial fiber bundle structure.
\begin{nonexample}\label{nonexam:8}
$M_i:=S^{2}(1)\setminus B^{2}(\frac{1}{i})$ does not fiber over its limit $S^{2}(1)$.
\end{nonexample}

Here is one of the main theorems of this paper.

\begin{theorem}\label{thm:4.2}
Let $(M_i,\partial M_i)$ be a sequence of Riemannian manifolds-with-boundary.\\
Suppose $M_i$ $\overset{GH}{\longrightarrow} X$,
with either:

(i) $X$ a compact Poincar\'{e} duality space, or

(ii) $X$ geodesically extendible.

\noindent
If $| K_{M_i} | \leq K$, $ | II_{\partial M_i} | \leq {\lambda}$ and $inj(M_i)\geq i_0 >0$, then $D^1\longrightarrow M_i \longrightarrow X$ is a locally trivial fiber bundle.
\end{theorem}

\vspace{0.3cm}

This theorem is fairly sharp, in that all of its hypotheses concerning $M_i$ (except perhaps for the upper curvature bound $K_{M_i} \leq K$) are necessary for the conclusion to hold.

\vspace{0.3cm}

Recall that a topological space $X$ is a \emph{Poincar\'{e} duality space} (over a group $G$ of coefficients) if there exists an integer $k$ (the formal dimension) and a fundamental homology class $[X]\in H_{k}(X; G)$ such that capping $H^{p}(X; G) \overset{\cap [X]}{\underset{\cong}{\longrightarrow}} H_{k-p}(X; G)$ is an isomorphism for all $0\leq p \leq k$. Such spaces include closed topological manifolds, quotients of Poincar\'{e} duality spaces by a group action (in rational coefficients), such as orbifolds,
and more generally any space homotopy equivalent to a Poincar\'{e} duality space.

\vspace{0.4cm}

Recall that a length space $X$ is said to be \emph{geodesically extendible} if for every 
nontrivial geodesic $\gamma : [0,L] \longrightarrow X$ with endpoint $\gamma(L)=x$, there exists a geodesic $\widetilde{\gamma} : [0,L+\epsilon] \longrightarrow X$ with $\epsilon > 0$ which properly extends $\gamma$ past $x$. By convention, a single point is considered a geodesically extendible space.

\vspace{0.5cm}

If the limit space in Theorem~\ref{thm:4.2} is assumed to be a single point, then the upper curvature bounds there can be omitted, and provided some supplementary additional assumptions are made, the result is that the manifold must be a disc:

\begin{theorem}[\cite{Wo:06}]\label{thm:4.5}
Let $(M_i,\partial M_i)$ be a sequence of $n$-dimensional Riemannian manifolds-with-boundary, where
$M_i \overset{GH}{\longrightarrow} X$, and where $X$ is a point.
\begin{itemize}

\item[(i)]
If $K^{-}\leq K_{M_i}$, \hspace{0.2cm} ${\lambda}^{-}\leq II_{\partial M_i}$, $i_{int}(M_i)\geq i_0 >0$, and $\frac{outrad(p_i)}{inrad(M_i)}\longrightarrow 1$, where $p_{i}\in M_i$ realize $inrad(M_i)$ and $outrad(p_i):=\underset{q\in \partial M_i}{\sup}d(p_i,q)$, then $M_i\underset{diffeo}{\approx}~D^n$.

\item[(ii)]
If \hspace{0.2cm} $inj(M_i) \geq i_0 > 0, \pi_{1}(\partial M_i)=1$, and $n\neq 3,4$ then $M_{i} \underset{homeo}{\approx} D^n$.

\end{itemize}
\end{theorem}

\vspace{0.4cm}

Recall Toponogov's splitting theorem: If an open complete Alexandrov space $X$ with $curv X \geq 0$ admits a line, then it splits isometrically as $X=Y \times \mathbb{R}$, for some non-negatively curved Alexandrov space $Y$. The boundary version of this states that if an Alexandrov space with boundary has $curv \geq 0$ and two boundary components, then it splits isometrically as a metric product $X = Y \times D^1$, where $curv Y \geq 0$. With weaker hypotheses in the manifold setting, one still retains the disc bundle structure:\\

\begin{theorem}[Rough Toponogov Splitting]\label{thm:Top}
Let $(M_i,\partial M_i)$ be a sequence of Riemannian manifolds-with-boundary.\\
If $K^{-}\leq K_{M_i}$, \hspace{0.4cm} ${\lambda}^{-}\leq II_{\partial M_i}\leq {\lambda}^{+}$, and $inrad(M_i) \longrightarrow 0$, then
$M_i$ can have at most two boundary components.\\
If in addition $d(M_i) \longrightarrow 0$
and each $M_i$ has exactly two boundary components
then $D^1 \longrightarrow M_i \longrightarrow W$ is a trivial fiber bundle, where $W$ is a
boundary component.
\end{theorem}

By rescaling the sequence to have diameter $d(M_i) \equiv 1$, this is equivalent to the statement that a compact manifold-with-boundary having almost-nonnegative curvature in the interior, almost-convex boundary and two boundary components is topologically a $D^1$-bundle over a boundary component.

\vspace{1cm}

\noindent
The next two theorems yield metric information about limits. In particular, certain limits are characterized as Alexandrov spaces of curvature bounded below.\\

\begin{theorem}\label{thm:4.3}
Let $(M_i, \partial M_i)$ be a sequence of Riemannian manifolds-with-boundary. Suppose $M_i \overset{GH}{\longrightarrow} X$ for some metric space $X$.

If $K^{-}\leq K_{M_i}, \ -\tau(\frac{1}{i}) \leq II_{\partial M_i} \leq {\lambda}^{+}$, 
then $curv X\geq K^{-}$.

\end{theorem}

\vspace{0.3cm}

\begin{theorem}\label{thm:4.6}
Suppose
$M^n$ is a fixed complete manifold without boundary, and $|K_{M}|\leq K$. Suppose $(M_{i}^{n},\partial M_i) \subset M^n$ are immersed submanifolds
of the same dimension,
with $II_{\partial M_i \hookrightarrow M_i} \geq {\lambda}^{-}$.
If $M_i \overset{GH}{\longrightarrow} X$ for some geodesically extendible space $X$, 
then

\qquad (1) $inrad(M_i) \longrightarrow 0$.

If in addition one also has $M_i \overset{H}{\longrightarrow} X$ in $M$, then

\qquad (2) $curv X \geq c(-K,{\lambda}^{-})$, and

\qquad (3) $X$ is a $C^1$-differentiable manifold.

\end{theorem}

\vspace{1cm}

This theorem differs from the previous Theorem~\ref{thm:4.3} in regard that here, the sequence of manifolds-with-boundary are assumed to lie in a fixed, common ambient manifold.
Also, in terms of curvature assumptions, this theorem differs from Theorem~\ref{thm:4.3}
in that the upper bound on $II_{\partial M_i}$ has been dropped, and an upper bound on $K_{M_i}$ imposed in its place.

\vspace{0.2cm}

After providing elementary examples to demonstrate that convergences as in these theorems actually occur, the remainder of the paper will be devoted to the proofs of these theorems.

\vspace{1.5cm}

\subsubsection{Examples}\label{sss:2.1.1}

Examples of convergences as in Theorems~\ref{thm:4.2}, \ref{thm:4.3}, and \ref{thm:4.6} are given below.\\
In examples \ref{exam:1}, \ref{exam:3} and \ref{exam:5}, $M_i$ collapse in dimension, whereas $\partial M_i$, with intrinsic metric, does not.
Examples 1-5 all inradius-collapse.\\

\begin{example}\label{exam:1}
Let $M_i :=$ portion of \, $\mathbb{R}^3$ bounded by the graphs of $f_{0}(x,y)=e^{-(x^2+y^2)}$ and $f_{i}(x,y)=e^{-(x^2+y^2)}+~\frac{1}{i}$.
Then $K_{M_i} \equiv 0$, $|II_{\partial M_i}| \leq {\lambda}$ for some $\lambda < \infty$, and the limit $X$ has $curv X \geq -{\lambda}^2$.

In this example $\partial M_i$, with intrinsic metric, has $curv (\partial M_i) \geq -{\lambda}^2$ and also GH-converges to $X$.

Also, $inj(M_i)\geq i_0$ for some $i_0$, and $M_i$ is diffeomorphically the product $D^{1}(\frac{1}{i}) \times \partial M_i$.
\end{example}

\vspace{0.1cm}

Each of the curvature hypotheses in Theorem~\ref{thm:4.2}, except perhaps for $K^{+}$, are necessary for its conclusion.
The bounds $K^{-}$ and ${\lambda}^{-}$ in Theorem~\ref{thm:4.2} are needed to ensure that the dimension of $X$ is eventually less than that of the $M_i$'s, clearly a necessary condition for a fibering. (See $\S$~\ref{sss:1.2.2}.)
The remark at the beginning of the next section, on p.\pageref{rmk:fiber}, shows why ${\lambda}^{+}$ is necessary for the fiber to be a $D^1$.

The following example demonstrates the sharpness of Theorem~\ref{thm:4.2} with regard to the injectivity radius bound. In particular, the injectivity radius bound is necessary for the fiber to be a $D^1$.

\begin{example}\label{exam:2}
$M_i:= S^{1}(\frac{1}{i}) \times D^{1}(\frac{1}{i}) \overset{GH}{\longrightarrow} pt$, under $K\equiv 0$, $II\equiv 0$, $inj(M_i)=\frac{\pi}{i} \longrightarrow 0$.

The universal Riemannian covers $\mathbb{R} \times D^{1}(\frac{1}{i})$ have $inj\equiv \infty$, so the sequence has bounded covering geometry.
Note that the fiber of the convergence $M_i \overset{GH}{\longrightarrow} pt$ is not a $D^1$, but instead a $D^1$-bundle, since the%
``souls'' of the $M_i$ collapse.
\end{example}

The author suspects that if the injectivity radius bound were to be omitted from Theorem~\ref{thm:4.2}, one would still get a fibering, where the fiber is itself a $D^{1}$-bundle over a closed manifold.

\vspace{0.5cm}

The following example illustrates Theorem~\ref{thm:4.3}.

\begin{example}\label{exam:3}
Let $\epsilon_i > 0$ be any sequence of real numbers tending to zero slower than $\frac{1}{i}$.

Let $M_i := \left( \partial B({[0,1]} \subset x\text{-axis}, \epsilon_i; \mathbb{R}^3) \cap \{ z\geq -\frac{1}{i} \}\right) \times S^{1}(r)$.

Then
$M_i \overset{GH}{\longrightarrow} X=[0,1]\times S^{1}(r)$, \ 
under $\ 0\leq K_{M_i} \longrightarrow \infty$, $-\tau(\frac{1}{i})\leq II_{\partial M_i} \leq 0$.

Also,
$inrad(M_i) = \frac{\pi}{2} \epsilon_i +  {\epsilon_i}{\sin}^{-1}\left( \frac{1}{i\epsilon_i}\right) \longrightarrow 0$
as $i\longrightarrow \infty$.

The limit $X$ has $curv X \geq 0$, but the $M_i$'s do not have Alexandrov curvature bounded below.

\end{example}

In this example, it need not hold that, in the intrinsic metric, $\partial M_i \overset{GH}{\longrightarrow} X$.
In fact $\partial M_i =S^1 \times S^{1}(r) \overset{GH}{\nrightarrow} [0,1]\times S^{1}(r)$.

Note also that $M_i \approx D^2 \times S^1$, a solid torus, does not fiber as locally trivial fiber bundle over its limit, a cylinder.
This phenomenon, and the fact that the sequence $\partial M_i$ with intrinsic metrics does not GH-converge to $X$, is due to the presence of geodesic terminals in the limit.

\vspace{0.5cm}

The next example illustrates Theorem~\ref{thm:4.6}.

\begin{example}\label{exam:4}
There exists a sequence of manifolds-with-boundary $M_i^{n}$ embedded in a fixed ambient manifold $M^{n}$, with $|K_{M}|\leq K$, ${\lambda}^{-}\leq II_{\partial M_i \hookrightarrow M_i}$, $M_i \overset{GH}{\longrightarrow} X$ and $M_i \overset{H}{\longrightarrow} X$, but such that $II_{\partial M_i \hookrightarrow M_i} \longrightarrow +\infty$:

To see this, take $M = \mathbb{R}^n$. Suppose
$X$ is any $l$-dimensional
linear subspace of $M$, or $l$-dimensional embedded round sphere of some fixed radius, where $n-l \geq 2$.

Let $M_i := \overline{B}(X, \epsilon_i ; \mathbb{R}^{n}) = \{ x\in \mathbb{R}^{n} : d(x,X) \leq \epsilon_i \}$ be an $\epsilon_i$-neighborhood of $X$, where $\epsilon_i \longrightarrow 0$.

Then $M_i \overset{H}{\longrightarrow} X$, and, equipped with intrinsic metric induced from $\mathbb{R}^n$, $M_i \overset{GH}{\longrightarrow} X$
under $K_{M_i}\equiv 0$, ${\lambda}^{-} \leq II_{\partial M}\longrightarrow \infty$.
While some principle curvatures remain bounded, at least one must diverge.
\end{example}

\vspace{0.1cm}

\begin{example}\label{exam:5}
There exists a sequence of metrics $g_i$ on $D^3$ such that $(D^3,g_i) \overset{GH}{\longrightarrow} (D^2,std)$, where 

$K\equiv 0, \ 0\leq II \longrightarrow \infty$.

Here, one can take discs $D_{i}^3 \subset \mathbb{R}^3$ with locally convex boundary ``flattening'' to the standard, planar $D^2 \subset \mathbb{R}^3$.

The boundaries $\partial (D^3,g_{i})$, with intrinsic metrics, GH-converge to the metric double $2(D^2,std) = $\\
$D^2 \underset{\partial D^2}{\cup} D^2$.
\end{example}

\vspace{0.2cm}

\noindent
Modifying this example, one has in connection with
Theorem~\ref{thm:4.6}
the following conjecture

\begin{nonexample}\label{nonexam:6}
A nonconvex planar domain $\Omega$, embedded in $\mathbb{R}^3$ via the standard inclusion $\Omega \subset \mathbb{R}^2 \subset \mathbb{R}^3$, cannot be the Hausdorff limit of any sequence of three-dimensional domains in $\mathbb{R}^3$ whose boundaries are smooth and have sectional curvature ($=$ product of principal curvatures) bounded below. (Note that $curv \Omega$ is not bounded below.)
\end{nonexample}

\vspace{0.1cm}
Geodesic extendibility is used in part
(1) of Theorem~\ref{thm:4.6}
to show that the inradius tends to zero. 
It is conceivable that the conclusion in part
(2)
namely,
a lower Alexandrov curvature bound for the limit,
might still be obtained if one did not assume that it was geodesically extendible, but rather, that the sectional curvatures of the boundary were bounded below, and the inradius of the approximating manifolds-with-boundary tended to zero.

\vspace{0.5cm}

\noindent
Finally, to illustrate all parts of Theorem~\ref{thm:4.2},
\begin{nonexample}\label{nonexam:7}
The manifold $S^{2} \setminus \coprod_{1}^{3}B^{2}$, consisting of the two-sphere with three disjoint balls removed, has three boundary components, so cannot admit a sequence of metrics satisfying the hypotheses of Theorem~\ref{thm:4.2}, let alone be a disc bundle. It is not even homotopic to a closed manifold.
\end{nonexample}

\vspace{1.4cm}

Here is a question related to Example~\ref{exam:2} and Theorem~\ref{thm:4.2}:

\begin{question*}
Let $(M_i,\partial M_i)$ be a sequence of Riemannian manifolds-with-boundary, and assume that\\
$M_i \overset{GH}{\longrightarrow} N$, where $N$ is a closed manifold.\\
If $|K_{M_i}|\leq K$, $|II_{\partial M_i}|\leq {\lambda}$, $M_i$ has uniform bounded covering geometry, then is $F \longrightarrow M_i \longrightarrow N$ a locally trivial fiber bundle, where $F$ is a $D^1$-bundle over a closed manifold?
\end{question*}

A sequence of manifolds-with-boundary is said to have uniform bounded covering geometry if the universal covers have sectional curvatures and second fundamental forms bilaterally bounded, in addition to a uniform lower bound to the intrinsic injectivity radius.\\

\subsection{Outline of the Proofs of Theorems}
\label{ss:2.2}

Each part of Theorem~\ref{thm:4.2} is independent from the other, and each involves a different idea.
However, both parts are reduced to a known fibering theorem (Theorem~\ref{thm:thin} below) after proving that inradius tends to zero.
Given the curvature hypotheses of Theorem~\ref{thm:4.2}, the condition $inj(M_i) \geq i_0 > 0$ for a manifold-with-boundary $M_i$ serves to prevent collapse in $(n-1)$ directions, so that the fiber being collapsed must be $1$-dimensional, namely, a $D^{1}$.

Theorem~\ref{thm:4.5}
is fairly direct to obtain.
The conditions in
part (i)
ensure that the manifolds in the sequence are eventually star-like with respect to a point.
Part (ii),
essentially homotopical, follows from facts about local geometric contractibility.
The reader is referred to~\cite{Wo:06} for the proof of Theorem~\ref{thm:4.5}.

Theorem~\ref{thm:Top} relies on the Alexandrov extension procedure from \cite{Wo:07} (see $\S$~\ref{ss:1.2}), and uses a special gluing lemma.
Whereas most 'perturbation' or $\epsilon$-versions of a theorem in Riemannian geometry are obtained from the corresponding usual proof, via an argument involving pasing to the limit, the proof here will proceed somewhat differently. 
Although the main idea still involves passing to the limit, the proof, interestingly enough, does not utilize a reduction to the standard Toponogov splitting theorem, but rather, relies on Yamaguchi's fibration theorem (Theorem~\ref{thm:4.4}).\\

The two fibering theorems referred to above differ not only in their curvature hypotheses, but also in their assumptions of the existence of a limit space.\\

\begin{theorem}\bf (Thin-Manifolds-with-Boundary Theorem, \cite{AB:98}). \label{thm:thin} \rm
\it
There is a dimension-independent constant $c$ ($\geq 0.075$) such that if $M$ is a complete connected manifold satisfying\\
$inrad(M)\cdot\max \{ \sup |K_{M}|, \sup |II|^2 \} < c^2$,
then either $M$ is diffeomorphic to the product of a manifold without boundary and an interval or $M$ can be doubly covered by such a product.
\end{theorem}

\vspace{0.5cm}

\rmk{
The upper bound ${\lambda}^{+}$ in the Thin-Manifolds-with-Boundary Theorem is necessary for its conclusion. Consider a sequence of thinning solid ellipsoids $M_i$ in $\mathbb{R}^{3}$, each homeomorphic to $D^{3}$. Then $K_{M_i}=0$, $0\leq II_{\partial M_i}\longrightarrow \infty$, and the curvature normalized inradius ${inrad}^{2}\cdot \max \{ \sup |K_{M_i}|,{({\lambda}^{-})}^{2} \}$ is identically $0$, yet $M_i$ is not a $D^1$-bundle. The other three curvature bounds ${\lambda}^{-}$, $K^{-}$, and $K^{+}$ are also necessary, as simple examples show:
the unit sphere $S^{2}(1)$ with three disjoint topological balls of increasing areas removed, 
the connected sums $\left(\mathbb{R}\times [0,\frac{1}{i}]\right)\#T^{2}(flat)$ of a thin strip of width $\frac{1}{i}$ with a flat torus of diameter $\frac{1}{i}$, and the diminishing hemispheres $S_{+}^{2}(\frac{1}{i})$, respectively. But observe that the first and the fourth examples given here are nevertheless disc bundles, which suggests that $K^{-}, {\lambda}^{-}$ are the only curvature conditions essential to obtain disc bundle structure.\label{rmk:fiber} \\
}

\begin{definition}
A map $f: Y \longrightarrow X$ between Alexandrov spaces is an \emph{$\epsilon$-almost Lipschitz submersion} if \\
1) $f$ is an $\epsilon$-Hausdorff approximation\\
2) for all $p,q \in Y$
\begin{align*}
\left| \dfrac{|f(p)f(q)|}{|pq|} - \sin \left( \underset{x\in f^{-1}(f(p))}{\inf} {\angle}qpx \right) \right| \leq \epsilon
\end{align*}
\end{definition}

In case $f^{-1}(f(p))$ consists of the single point $p$,
part (2) of the definition is taken to state that\\
$\left| \dfrac{|f(p)f(q)|}{|pq|} - 1 \right| \leq \epsilon$.

\vspace{0.5cm}

\begin{theorem}[\cite{Y:96}]\label{thm:4.4}
Given a positive integer $n$ and a number ${\mu}_0 > 0$, there exist positive numbers $\delta=\delta_n$ and $\epsilon = \epsilon_{\mu_0}$ such that if $X$ is an $n$-dimensional complete Alexandrov space with curvature $\geq -1$ and $\delta$-str.rad$(X) > \mu_0$ and $d_{GH}(Y,X) \leq \epsilon$ for some complete Alexandrov space $Y$ with curvature $\geq -1$, then there exists a $\tau_{n,\mu_0}(\delta,\epsilon)$-almost Lipschitz submersion $f : Y \longrightarrow X$. Furthermore, $f$ is a locally trivial fiber bundle when both $Y$ and $X$ have $C^1$-differentiable structures.
\end{theorem}

\vspace{1cm}

\rmk{
For spaces of curvature bounded below, two major fibering theorems involving the Gromov-Hausdorff distance are known: Yamaguchi's Theorem above, and the Topological Stability Theorem in \cite{P:92}.
The first requires metric completeness and the second compactness. However, 
neither requires geodesic extendibility of the approximating space.

For spaces of curvature bounded above, there seem to be known only two fibering theorems: Theorem~\ref{thm:thin}, specific to manifolds (and not involving Gromov-Hausdorff distance \emph{per se}), and Nagano's Theorem in~\cite{N:02}, concerning geodesically extendible Alexandrov spaces of curvature bounded above (with also a mild condition on geodesic branching, and a uniform $CAT_{k}$ radius bound). The latter theorem requires geodesic extendibility of both the fixed (limit) space and the approximating space, and concludes the existence of a bi-Lipschitz homeomorphism between them.
}

\vspace{0.5cm}

\vspace{1.5cm}
\noindent
With this background, the proof of Theorems~\ref{thm:Yam}, \ref{thm:4.2}, \ref{thm:Top}, \ref{thm:4.3}, \ref{thm:4.6} commences now.\\

\subsection[Proof of Theorem~\ref{thm:Yam}]{Proof of Theorem~\ref{thm:Yam}}\label{ss:2.2.0}

\begin{proof}[Proof of Theorem~\ref{thm:Yam}(ii)]
Part (ii)
is a direct corollary of Yamaguchi's fibering theorem (Theorem~\ref{thm:4.4}).
\end{proof}

The following Proposition is key in the proof of
part (i), as well as in
Theorem~\ref{thm:4.2}.

\newpage

\begin{proposition}[No long segments uniformly transverse to the boundary]\label{prop:inrad}
Let $(M_i,\partial M_i)$ be a sequence of Riemannian manifolds with locally convex boundary and $K^{-}\leq K_{M_i}$.
Assume that $M_i \overset{GH}{\longrightarrow} X$, where $X$ is
a geodesically extendible space (e.g. a closed manifold).
Then for all $R>0$, $\theta > 0$, there exists an $i_{0}$ such that for all $i \geq i_{0}$, $M_i$ admits no minimal segment $[q_{i}p_{i}]_{M_{i}}$ of length $\geq R$ such that
$p_{i} \in \partial M_{i}$
and such that at $p_{i}$, ${\angle}^{M_i}([p_{i}q_{i}]_{M_i}', \partial M_i) \geq \theta > 0$.
\end{proposition}

\rmk{
The constant $i_{0}$ will depend on the extendibility radius of $X$.
}

\rmk{
The Proposition holds if the $M_i$ are Alexandrov spaces with boundary, or pairs $(M_i,E_i)$, consisting of closed Alexandrov spaces together with an extremal subset.
}

\begin{proof}
$ \ $

\noindent
\textcircled{1} For each $i$, let $p_i \in \partial M_i$ be an arbitrary point.\\
\quad $\textcircled{2}$ Suppose there exist minimal segments
$[q_{i}p_{i}]$
with ${\angle}^{M_i}([p_{i}q_{i}]',\partial M_i) \geq \theta > 0$
such that $|q_{i}p_{i}| = R$ for some fixed $R > 0$ which is assumed independent of $i$.\\
\quad $\textcircled{3}$ $[q_{i}p_{i}]$ converges to a minimal segment $[qp] \subset X$ (after possibly taking a subsequence).\\
\quad $\textcircled{4}$ Extend $[qp]$ to a geodesic $[qpr]$ with $|pr|=\delta > 0$.

$[qpr]$ may be assumed minimal by restricting $R$ and/or $\delta$, if necessary\\
\quad $\textcircled{5}$ There exist points $r_{i}$ in $M_i$ converging to $r$.\\\\
Then\\

(1) ${\angle}^{M_i}q_{i}p_{i}r_{i} \leq \pi - \theta$ for all $i$\\

(2) ${\overline{{\angle}}}_{k}^{X}qpr = {\angle}^{X}qpr = \pi$\\

(3) $ {{\angle}}^{M_i}q_{i}p_{i}r_{i}\geq {\overline{{\angle}}}_{k}^{M_i}q_{i}p_{i}r_{i} $ since $curv M_i\geq k=K^{-}$\\
Now
\begin{align*}
|q_{i}p_{i}| &\longrightarrow |qp| = R,\\
|p_{i}r_{i}| &\longrightarrow |pr| = \delta,\\
\text{ and } \qquad |q_{i}r_{i}| &\longrightarrow |qr| = R+\delta
\end{align*}
imply

(4) ${\overline{{\angle}}}_{k}^{M_i}q_{i}p_{i}r_{i} \longrightarrow {\overline{{\angle}}}_{k}^{X}qpr$\\
by the law of cosines in the model space.
But (2),(3),(4) together contradict (1). So $R$ cannot be chosen uniformly positive and independent of $i$ after all.
\end{proof}

\vspace{0.5cm}

\begin{proof}[Proof of Theorem~\ref{thm:Yam}(i)]

Suppose $M_i \overset{GH}{\longrightarrow} N$.

By the Soul Theorem, 
each $M_i$ is a $D^{k_i}$ bundle over the totally convex closed submanifold $S_i$:
\begin{align}\label{eq:1}
D^{k_i} \longrightarrow M_i \longrightarrow S_i
\end{align}

Since $S_{i} \subset M_i$ is totally convex,
\begin{align}\label{eq:2}
d_{GH}(S_{i},M_i) \leq d_{H}^{M_i}(S_{i},M_i).
\end{align}

Now Proposition~\ref{prop:inrad} implies
in particular that the Sharafutdinov retraction $M_i \longrightarrow S_{i}$ takes less and less time, as $i \longrightarrow \infty$, when the generalized gradient flow it follows is normalized to unit speed.

Therefore
$d_{H}^{M_i}(S_{i},M_i) \longrightarrow 0$
and so by (\ref{eq:2}),
$d_{GH}(S_{i},M_i) \longrightarrow 0$
as $i \longrightarrow \infty$.

Hence $d_{GH}(S_{i},N) \longrightarrow 0$.

Since $S_i$ is a closed manifold with non-negative curvature, 
one can apply Yamaguchi's fibration theorem for manifolds \cite{Y:91} (or \cite{Y:96}) to obtain a locally trivial fiber bundle
\begin{align}
F \longrightarrow S_i \longrightarrow N
\end{align}
for sufficiently large $i$, for some closed manifold fiber $F$.
Up to diffeomorphism, composing this fiber bundle map with the fiber bundle map
(\ref{eq:1}) above yields a bundle
$D^{k_i}\widetilde{\times} F \longrightarrow M_i \longrightarrow N$ for sufficiently large $i$, where $D^{k_i}\widetilde{\times} F$ denotes a (possibly twisted) product.
\end{proof}

\vspace{0.5cm}

\subsection[Proof of Theorem~\ref{thm:4.2} (i)]{Proof of Theorem~\ref{thm:4.2} (i)}\label{sss:2.2.1}

\noindent
Part (i)
utilizes the following two lemmata. The reader might want to scan Appendices $\S$\ref{s:1}, $\S$\ref{ss:1.2} and $\S$\ref{s:A3} before reading the proofs.\\

\begin{lemma}\label{lem:4.8}
Suppose $(M_{i}, \partial M_i)$ is a sequence of $n$-dimensional Riemannian manifolds-with-boundary such that $K^{-} \leq K_{M_i}, \ |II| \leq \lambda, \ d(M_i) \leq d, \ \ M_{i}^{n} \overset{GH}{\longrightarrow} N$, with $N$ a closed topological manifold (or more generally, a Poincar\'{e} duality space). Then $dim_{\mathcal{H}} N < n$, i.e. the sequence $\{ M_i \} $ volume collapses.
\end{lemma}

\begin{proof}
Assume $dim_{\mathcal{H}} N \geq n$. Then by the dimension estimate (Corollary~\ref{cor:2.2}), $dim_{\mathcal{H}} N \leq n$, so $dim_{\mathcal{H}} N = n$. By Proposition~\ref{prop:2.1}, there is an embedding
$M_i \overset{\approx}{\hookrightarrow} \widetilde{M_i}$, where $curv (\widetilde{M_i})\geq k(K^{-},\lambda)$.
Note that $d(\widetilde{M_i}) \leq d + 2t_0$.
Extracting a subsequence by precompactness,
$\widetilde{M_i} \overset{GH}{\longrightarrow} Y$ for some compact Alexandrov space $Y$ of dimension $k\leq n$. 

The projection maps $\pi_i: \widetilde{M_i} \longrightarrow M_i$ defined by
\begin{align*}
\pi_i = \begin{cases}
id_{M_i} &\text{ on } M_i\\
\text{orthogonal projection onto base } \partial C_{M_i}=\partial M_i &\text{ on } C_{M_i}
\end{cases}
\end{align*}
are surjective. They are also
Lipschitz: $|\pi_{i}(x)\pi_{i}(y)| \leq L |xy|$ for all $x,y\in \widetilde{M_i}$, for some constant $L$. In fact, using a common warping function for all $M_i$, one may take $L=\frac{1}{\epsilon}$, by Lemma~\ref{lem:2.4}.
So by Proposition~\ref{prop:8.1} in $\S$\ref{s:A3}, there exists a surjective, $L$-Lipschitz map $\pi : Y \longrightarrow N$. In particular, $n = dim N \leq dim Y = k$. Hence $k=n$.

But then, by the Topological Stability Theorem for Alexandrov spaces of curvature bounded below, $\widetilde{M_i} \underset{homeo}{\approx} Y$.
By Proposition~\ref{prop:2.8}, there is a homotopy equivalence $Y \overset{\simeq}{\longrightarrow} N$ given by the deformation retraction arising from the natural deformation retractions $\widetilde{M_i} \overset{\simeq}{\longrightarrow} M_i$.
So $M_{i}^{n} \approx \widetilde{M_i} \approx Y \simeq N^n$ for all sufficiently large $i$. Fix such $i$ and let $M=M_i$.

\begin{alignat*}{2}
H^{p}(M,\partial M) &\cong H_{n-p}(M)
&& \quad \text{for any $0 \leq p \leq n$, by duality,}
\text{ since $M$ compact and $dim M = n$}\\
&\cong H_{n-p}(N) && \quad \text{since } M \simeq N \\
&\cong H^{p}(N) && \quad \text{by duality, since $N$ closed and $dim N = n$}\\
&\cong H^{p}(M) && \quad \text{since } M \simeq N
\end{alignat*}
For the first duality isomorphism given above, see, e.g., corollary 9.3, p.351 of \cite{Br:97}.
(Note that compactness of $M$ (and thence of $\partial M$) is essential for these isomorphisms, e.g., consider $M = {\mathbb{R}}_{+}^{n} \simeq N = \mathbb{R}^n$.)\\

The long exact cohomology sequence of the pair $(M,\partial M)$
\begin{align*}
\cdots \longrightarrow
H^{n-1}(M,\partial M) \overset{\cong}{\longrightarrow} H^{n-1}(M) 
\longrightarrow H^{n-1}(\partial M) \longrightarrow  H^{n}(M,\partial M) \overset{\cong}{\longrightarrow}H^{n}(M)
\end{align*}
then yields that $H^{p}(\partial M) = 0$ for all $0 \leq p \leq n-1$, which
contradicts the fact that (in ${\mathbb{Z}}_{2}$ coefficients) $H^{n-1}(\partial M) = H_{c}^{n-1}(\partial M) \geq {\mathbb{Z}}_{2}$ since $\partial M$ is compact and non-empty.
\end{proof}

\vspace{0.4cm}

\noindent
The second lemma which will be used in the proof of Theorem~\ref{thm:4.2} (i) is a version of Berger's isoembolic volume inequality.\\

\begin{lemma}\label{lem:4.5}
Let $(M,\partial M)$ be a complete manifold of dimension $n$. Then for some constant $c(n)>0$, \
$vol(M) \geq c(n) { \min \{ inrad(M), i_{int}(M) \} }^{n}$
\end{lemma}

\begin{proof}
If $B(p,r)$ is any metric ball lying entirely in the interior of $M$, and $i_{int}(x) \geq 2r$ for $x\in B(p,r)$, then all geodesics emanating from $x$ are minimizing at least until they hit $\partial B(p,r)$. Therefore, if $r\leq \frac{i_{int}(M)}{2}$, any geodesic emanating from an interior point of $B(p,r)$ is minimizing at least until it hits $\partial B(p,r)$.
(In the terminology of \cite{C:80}, this means $\widetilde{\omega}=1$ on $B(p,r)$, or equivalently, that the cut locus of any interior point of $B(p,r)$ w.r.t. the usual exponential map lies outside $B(p,r)$.) Now by
\cite[Theorem 11]{C:80}, 
we have that for $0 < t \leq r \leq \frac{i_{int}(M)}{2}$,
\begin{align*}
\dfrac{vol(\partial B(p,t))}{{vol(B(p,t))}^{(n-1)/n}} \geq c^{1/n}
\end{align*}
where $c=2^{n-1}{vol(S^{n-1}(1))}^n / {vol(S^{n}(1))}^{n-1}$.
Integrating both sides w.r.t. $t$ (from $0$ to $r$) yields
\begin{align*}
n \cdot vol(B(p,r))^{1/n} \geq c^{1/n}\cdot r
\end{align*}
so
\begin{align} \label{ineq:v1}
vol(B(p,r)) \geq c \cdot r^{n}
\end{align}
for some new constant $c$.

Now choose $p \in M$ to realize $inrad(M)$, so that $inrad(p)=inrad(M)$.

Let $r:= \min \{ inrad(M),\frac{i_{int}(M)}{2} \}$. Then $B(p,r)$ is an interior ball, and
\begin{align*}
vol(M) &\geq vol(B(p,inrad(M)))\\
&\geq vol(B(p, r))\\
&\geq c(n)\cdot r^{n} \quad \text{ by (\ref{ineq:v1})}.
\end{align*}
So for a less sharp constant,
\begin{align*}
vol(M) &\geq c(n) \cdot {\min \{ inrad(M), i_{int}(M) \}}^{n}. \qedhere
\end{align*}
\end{proof}

\vspace{0.6cm}

Combining Lemma~\ref{lem:4.8} and Lemma~\ref{lem:4.5}, and recalling that $i_{int}(M) \geq inj(M)$, one obtains

\begin{proposition}
If $K^{-} \leq K_{M_i}, \ |II_{\partial M_i}| \leq \lambda, \ d(M_i) \leq d, \ inj(M_i)\geq i_0 > 0, \ M_{i}^{n} \overset{GH}{\longrightarrow} N$, with $N$ a topological manifold without boundary, then $inrad(M_i) \longrightarrow 0$.
\end{proposition}

\vspace{0.7cm}

\begin{proof}[Proof of Theorem~\ref{thm:4.2}(i)]
Immediate from this proposition and Theorem~\ref{thm:thin}.
\end{proof}

\vspace{1cm}

\subsection{Proof of Theorem~\ref{thm:4.2} (ii)}
\label{s:A1}

This section provides a geometric proof to Theorem~\ref{thm:4.2}(ii), 
employing lemmata which have independent interest. Although this proof, which uses only local arguments, is longer, it is in a way more general than the proof of Theorem~\ref{thm:4.2}(i) given in \S\ref{sss:2.2.1}, in that the limit $X$ need not be assumed a
Poincar\'{e} duality space, but merely a geodesically extendible length space with positive injectivity radius. For instance, the one-point union $S^1 \vee S^1$ is such a space, which is not a Poincar\'{e} duality space [$H^{1}(S^1 \vee S^1)=\mathbb{Z} \oplus \mathbb{Z}$ whereas $H_{0}(S^1 \vee S^1)=\mathbb{Z}$ ]. However, that proof in \S\ref{sss:2.2.1} covers cases that the present one does not, e.g., metric cones or suspensions over manifolds.\\

In fact, the proof given here
in \S\ref{s:A1}
essentially shows that one only needs to assume
that $X$ is a geodesic metric space which is 'weakly geodesically-extendible' in the sense that (i) Alexandrov angles exist\footnote{
i.e., $\underset{s,t \rightarrow 0}{\lim}{\angle}^{X}\alpha(s)p\beta(t)$ exists for any point $p\in X$ and any (shortest) curves $\alpha$ and $\beta$ with $\alpha(0)=\beta(0)=p$} and (ii) for any nontrivial geodesic $\gamma : [0,L] \longrightarrow X$, there exists a nontrivial geodesic $\sigma : [0,\epsilon] \longrightarrow X$, $(\epsilon>0)$ for which $\sigma(0)=\gamma(0)$ and ${\angle}(\gamma'(0), \sigma'(0) ) \geq \theta$, where $\theta > \frac{\pi}{2}$. (Here the length $\epsilon$ of $\sigma$ and the angle $\theta$ are allowed to depend on the point $\gamma(0) \in X$.) Such a class of spaces includes Alexandrov spaces of curvature bounded below which have a uniform lower $(n,\delta)$-strain radius bound, where $\delta < \frac{\pi}{2}$. But it also \emph{a priori} includes spaces which have no finite lower or upper Alexandrov curvature bound, such as the prism-block, and certain other cell complexes.\\

Some of these spaces may be ruled out of consideration by the fact that the hypotheses on the $M_i$ alone force their limit $X$ to be an Alexandrov space of curvature bounded above (by $k=k(K^{+},{\lambda}^{-})$), with uniformly bounded lower $CAT_{k}$-radius.
In light of this, and also in light of the proof of Theorem~\ref{thm:4.2}(ii) given in this section, Theorem~\ref{thm:4.2}(ii) may be viewed a fibering theorem in the category of CBA spaces. Even with this data about the limit space, however, there are still many geodesically extendible CBA spaces which are not manifolds or Poincar\'{e} duality spaces (a large class of examples arises from gluing, for instance). \\

The main ingredients of the proof of Theorem~\ref{thm:4.2}(ii) are: the extension produced in $\S$\ref{ss:1.2}, the injectivity estimates from
$\S$\ref{s:1},
arc/chord comparison ($\S$\ref{s:A2}), and a lemma on angles, obtained for instance from Proposition~\ref{prop:4.2} below.
The main idea of the proof is fairly well expressed already in the special case when the boundary is locally convex.

\vspace{1cm}

\subsubsection{Locally convex boundary case}
\label{ss:A1.1}

By Proposition~\ref{prop:inrad}, $inrad(M_i) \longrightarrow 0$.
Now invoke Theorem~\ref{thm:thin} to get a disc bundle structure $D^{1} \longrightarrow M \longrightarrow X$.

\vspace{1cm}

\subsubsection{General case}
\label{ss:A1.2}

\begin{proof}[Proof of Theorem~\ref{thm:4.2}(ii)]
In conjunction with the Thin-Manifolds-with-Boundary Theorem, the proof is completed by the following

\begin{lemma}\label{lem:4.3}
Suppose $(M_i,\partial M_i)\overset{GH}{\longrightarrow}X$, with $X$ a geodesically extendible length space,
under $K^{-}\leq K_{M_i} \leq K^{+}$, ${\lambda}^{-}\leq II_{\partial M_i}\leq {\lambda}^{+}$, $inj(M_i)\geq i_0$. Then $inrad(M_i)\longrightarrow 0$.
\end{lemma}

\end{proof}

The proofs of the following three lemmata, which prepare for Lemma~\ref{lem:4.3}, will be given later.

Although $M$ need not have curvature bounded below in the Alexandrov sense, it is possible to obtain an estimate to replace step (3) in Proposition~\ref{prop:inrad}.
The first lemma, the main angle estimate, establishes that angles of certain sufficiently small triangles in $M$ are comparable to the angles considered in the extension of $M$.

\begin{lemma}\label{sublem:2}
Let $\widetilde{M}$ be an Alexandrov extension of $M$ as in Proposition~\ref{prop:2.1}. 
For all $x,y,z\in M$ with $|xy|_{M},|yz|_{M}\leq R$, $y\in \partial M$ and $[xy]_{M}\perp \partial M$,
\begin{align*}
|{\angle}^{M}xyz - {\angle}^{\widetilde{M}}xyz|\leq {\tau}_{K^{\pm},{\lambda}^{\pm},i_0}(R)
\end{align*}
\end{lemma}

\vspace{0.5cm}

Now we introduce a notion which will be useful in the sequel.
As discussed in \cite{Lyt:01} and \cite{AB:03}, it is a synthetic metric space surrogate for the second-fundamental form of subspaces, when the spaces involved might not be smooth. Such a notion, formulated in terms of the metric distance functions, is thus relevant when considering Gromov-Hausdorff convergence. \\

\begin{definition}[\cite{Lyt:01}]\label{def:cconvex} 
A subspace $Z$ of a length space $X$ is said to be \emph{$(C,2,\rho)$-convex} (for some $C\geq 0$ and positive function $\rho$ on $X$) if for all $w\in Z$
and any $x,y\in B(w,\rho(w); X) \cap Z$, the metrics satisfy $\ d_{Z}(x,y) \leq d_{X}(x,y)+Cd_{X}^{3}(x,y)$.
\end{definition}

\vspace{0.1cm}

In applications, one typically takes $\rho$ to be a positive constant.
If one does not want to emphasize the function $\rho$, or if it is understood, $Z$ is called merely \emph{$(C,2)$-convex}. 
For instance, a $(0,2)$-convex subset is simply a locally convex subset.
Note that in the definition, the subspace $Z$ is not a-priori required to be a length space itself.\\

Under the hypotheses of Theorem~\ref{thm:4.2}, both $M$ and $\partial M$ are $(C,2)$-convexly embedded in $\widetilde{M}$:

\begin{lemma}\label{sublem:3}
$|xy|_{\widetilde{M}} \leq |xy|_{M} \leq |xy|_{\widetilde{M}} + C{|xy|_{\widetilde{M}}}^{3}$ for all $x,y\in M$ sufficiently close, where $C=C(K^{\pm},{\lambda}^{\pm},i_0)\geq 0$\\
Furthermore, $|xy|_{\widetilde{M}} \leq |xy|_{\partial M} \leq |xy|_{\widetilde{M}} + C{|xy|_{\widetilde{M}}}^{3}$ for all $x,y\in \partial M$ sufficiently close, where $C=C(K^{\pm},{\lambda}^{\pm},i_0)\geq 0$\\
\end{lemma}

The last lemma required is a perturbation of the fact that if the normalized excess at a vertex of a model triangle is small, the angle at that vertex is close to $\pi$:

\begin{lemma}\label{sublem:4}
Fix a constant $C\geq 0$. For a triangle in a model space with side lengths $a,b,c$, $c$ sufficiently small ($Cc^3 << c,a,b$), $\dfrac{a+b-(c+Cc^3)}{\min \{ a,b \}}\leq \tau_{a,b,c}(Cc^3)$ implies $\theta\geq \pi - \tau_{a,b,c}(Cc^3)$, where $\theta$ denotes the angle opposite side $c$.
\end{lemma}

\vspace{0.8cm}

Given these three lemmata, the proof of Lemma~\ref{lem:4.3}, hence the proof Theorem~\ref{thm:4.2}(ii), proceeds as follows.\\

\begin{proof}[Proof of Lemma~\ref{lem:4.3}]
The lemma is trivial if $X=pt$ since then $d(M_i)\longrightarrow 0$, which forces $inrad(M_i)\longrightarrow 0$. It may be supposed $X \neq pt$.
A lower bound $inj(X) \geq i_{0} > 0$ is inherited from the assumed bound for the $M_i$'s.

Assume to the contrary that $inrad(M_i)\geq i_2 > 0$ for some $i_2$ and all $i$. By Lemma~\ref{lem:4.1}, there exists a sequence $p_i\in \partial M_i$ such that $p_i\longrightarrow p\in X$ and $i_{\partial}(p_i)\geq i_2$ for all $i$.
Let $[q_{i}p_{i}]_{M_i}$ be a geodesic segment normal to $\partial M_i$. By definition of $i_{\partial}(p_i)$, a lower bound $i_{\partial}(p_i)\geq i_2$ implies that
$q_i$ may be chosen such that the minimal segment $[q_{i}p_{i}]_{M_i}$ has length
\begin{align}\label{ref:4.1}
|q_{i}p_{i}|_{M_i}\geq i_2 > 0
\end{align}

The minimal segments $\{ [q_{i}p_{i}]_{M_i} \}$ have as limit a minimal segment $[qp]_{X}\subseteq X$ (of length $\geq i_2$), by Lemma~\ref{lem:4.4}. Since $X$ is geodesically extendible and $inj(X)\geq i_0>0$, there exists a point $z\in X$ such that $[qp]_{X}\cup [pz]_{X}$ is also a segment, where say, $|pz|_{X}= \delta > 0$ for some fixed $0<2\delta<\min\{i_0,i_2 \}\leq i_2\leq i_{\partial}$. Let $f_{i} : X \longrightarrow M_i$ be an $\epsilon_i$-Hausdorff-approximation.
Set $z_{i}:=f_{i}(z)$.\\

Restricting if necessary, it may be assumed that $|q_{i}p_{i}|_{M_i}=R$ for each sufficiently large $i$. ($R$ will later be chosen smaller).

Let $\widetilde{M_i}$ be an Alexandrov extension of $M_i$ as in Proposition~\ref{prop:2.1}.\\

Since $[q_{i}p_{i}]_{M_i} \perp \partial M_i$, it is obvious that
\begin{align}
{\angle}^{M_i}q_{i}p_{i}z_{i} &\leq \pi/2. \label{ineq:5.1}
\end{align}

Now
\begin{align}\label{ineq:5.2}
|{\angle}^{M_i}q_{i}p_{i}z_{i} - {\angle}^{\widetilde{M_i}}q_{i}p_{i}z_{i}| &\leq {\tau}_{K^{\pm},{\lambda}^{\pm},i_0}(R)
\end{align}
by Lemma~\ref{sublem:2},
and
\begin{align}\label{ineq:5.3}
{\angle}^{\widetilde{M_i}}q_{i}p_{i}z_{i} &\geq {\overline{\angle}}^{\widetilde{M_i}}q_{i}p_{i}z_{i}
\end{align}
since $curv \widetilde{M_i}$ is bounded below.

Furthermore, by Lemma~\ref{sublem:3},
\begin{align*}
|q_{i}p_{i}|_{\widetilde{M_i}}+|p_{i}z_{i}|_{\widetilde{M_i}}-\left( |q_{i}z_{i}|_{\widetilde{M_i}} + C |q_{i}z_{i}|_{\widetilde{M_i}}^{3} \right)
&\leq |q_{i}p_{i}|_{M_i} + |p_{i}z_{i}|_{M_i} - |q_{i}z_{i}|_{M_i}\\
&\leq (|qp|_{X}+\epsilon_i) + (|pz|_{X}+\epsilon_i) - (|qz|_{X}-\epsilon_i)\\
&= 3{\epsilon}_{i} 
\end{align*}
since $[qpz]$ is a segment in $X$.

If $|q_{i}p_{i}|_{\widetilde{M_i}}=|q_{i}p_{i}|_{M_i}=R$,
then by Lemma~\ref{sublem:3} again, 
\begin{align*}
|p_{i}z_{i}|_{\widetilde{M_i}} \geq |p_{i}z_{i}|_{M_i} - C |p_{i}z_{i}|_{M_i}^3  \geq (\delta-\epsilon_i) - C(\delta+\epsilon_i)^3.
\end{align*}

\emph{Assuming} that $R$ is bounded below by a fixed positive constant, then eventually $R>\epsilon_i$, and the normalized excess of a model triangle with sidelengths $a=|q_{i}p_{i}|_{\widetilde{M_i}}$, $b=|p_{i}z_{i}|_{\widetilde{M_i}}$, and $c+Cc^3=\left( |q_{i}z_{i}|_{\widetilde{M_i}} + C |q_{i}z_{i}|_{\widetilde{M_i}}^{3} \right)$ 
is bounded above by
$\frac{3\epsilon_i}{\min\{ R,(\delta-\epsilon_i)-C{(\delta+\epsilon_i)}^3\}}$,
hence tends to $0$ as $\epsilon_i \longrightarrow 0$.
Therefore
\begin{align}\label{ineq:5.4}
{\overline{\angle}}^{\widetilde{M_i}}q_{i}p_{i}z_{i}\geq \pi - \tau(|q_{i}z_{i}|_{\widetilde{M_i}})
\end{align}
by Lemma~\ref{sublem:4}.

\noindent
If
$R$ is such that ${\tau}_{K^{\pm},{\lambda}^{\pm},i_0}(R) \leq \pi /12$\\
and $i$ such that $\tau(|q_{i}z_{i}|_{\widetilde{M_i}})\leq \pi /12$ and $\epsilon_{i} < R$\\
then by (\ref{ineq:5.2}), (\ref{ineq:5.3}), (\ref{ineq:5.4}),
\begin{align*}
{\angle}^{M_i}q_{i}p_{i}z_{i}
&\geq \pi - \tau(|q_{i}z_{i}|_{\widetilde{M_i}})  - {\tau}_{K^{\pm},{\lambda}^{\pm},i_0}(R)\\
&\geq \pi - \frac{2\pi}{12} = \frac{5\pi}{6}
\end{align*}
contradiction to (\ref{ineq:5.1}).\\

Therefore
$R\longrightarrow 0$ as $i$ increases, contradiction to (\ref{ref:4.1}). So we must have $inrad(M_i)\longrightarrow 0$.
\end{proof}

\vspace{1cm}

\vspace{1cm}
No arc/chord comparison is available in $\widetilde{M}$ in general, since
$inj(\widetilde{M})=0$ can occur.\\
However,
$M \in \mathcal{M}(n,K^{\pm},{\lambda}^{\pm})$ implies $C_{M} \in \mathcal{M}(n,K_{C}^{\pm},{\lambda}_{C}^{\pm})$ for some constants $K_{C}^{\pm},{\lambda}_{C}^{\pm}$.\\
Explicitly, one may take
\begin{align*}
K_{C}^{-} &= \min \{ K_{C,radial}^{-}, K_{C,tangential}^{-} \}\\
&= \min \{ \underset{t}{\inf} \left( -\frac{\phi''(t)}{\phi(t)} \right), \underset{t}{\inf} \frac{1}{{\phi}^2(t)} \left[ K_{\partial M}^{-} - {|\phi'(t)|}^2 \right]  \}\\
&\geq \min \{ r.h.s. \text{ of } (\ref{functineq:4}),  r.h.s. \text{ of } (\ref{functineq:6}) \}\\
&= c(\phi,K_{M}^{-},{\lambda}^{\pm}) > -\infty\\
K_{C}^{+}
&= \max \{  \underset{t}{\sup} \left| \frac{\phi''(t)}{\phi(t)} \right|,
\underset{t}{\sup} \frac{1}{{\phi}^2(t)} \left[ K_{\partial M}^{+} - {|\phi'(t)|}^2 \right] \}\\
&\leq \max \{  \underset{t}{\sup} \left| \frac{\phi''(t)}{\phi(t)} \right|, \underset{t}{\sup} \frac{1}{{\phi}^2(t)} \left[ \left( K_{M}^{+} + \max \{  {|{\lambda}^{-}|}^2,{|{\lambda}^{+}|}^2 \} \right) - 0 \right] \}\\
&= c(\phi,K_{M}^{+},{\lambda}^{\pm}) < \infty\\
{\lambda}_{C}^{-} &= 0\\
{\lambda}_{C}^{+} &= |\min \{ 0, {\lambda}^{-} \}|
\end{align*}

In particular, $C_{M}$ has an upper Alexandrov curvature bound
\begin{align}
curv C_{M} \leq K_{C}^{+}
\end{align}
And
\begin{align*}
inj(\partial C_{M}) &= inj(\partial M)\\
&\geq c(n,K^{+},{\lambda}^{\pm},inj(M)\geq i_0)>0
\end{align*}
by Proposition~\ref{prop:1.2}(iv),
so
\begin{align*}
inj(C_{M})\geq c(n,K^{\pm},{\lambda}^{\pm},inj(\partial C_{M}), i_{\partial}=t_0,d)>0
\end{align*}
by Proposition~\ref{prop:1.2}(iii).
Hence
\begin{align}
CAT_{K_{C}^{+}}Rad(C_{M})
\geq \min \{ \frac{\pi}{2\sqrt{K_{C}^{+}}}, inj(C_{M}) \}
\geq c(n,K^{\pm},{\lambda}^{\pm},i_{0},d) > 0.
\end{align}
by (\cite{AB:96}, Theorem 4.3, p.78).
Therefore arc/chord comparison is available in $C_{M}$ within this radius.

\vspace{1cm}

Arc/chord comparison then yields $(C,2)$-convexity of $M$ and $\partial M$ in $\widetilde{M}$. More precisely, it yields\\
$(C,2,\min \{ CAT_{k_{0}}Rad(C_{M}), 1 \})$- and $(C,2,\min \{ CAT_{k_{0}}Rad(M), CAT_{k_{0}}Rad(C_{M}), 1 \})$-convexity, respectively, if $curv M \leq k_{0}$ and $curv C_{M} \leq k_{0}$:

\begin{proof}[Proof of Lemma~\ref{sublem:3}]
$|xy|_{\widetilde{M}} \leq |xy|_{M}$ for all $x,y\in M$ since $M\subseteq \widetilde{M}$.\\
Suppose $|xy|_{\widetilde{M}} \leq 1$, and also $|xy|_{\widetilde{M}} < CAT_{k_{0}}Rad(C_{M})$, where $curv C_{M} \leq k_{0}$.\\
Let $[xy]_{\widetilde{M}}=\cup [x_{i_j}x_{i_{j+1}}]_{M} \bigcup \cup [x_{i_k}x_{i_{k+1}}]_{C_{M}}$, where $x_{i_k} \in \partial C_{M}$.\\
Note that $\cup [x_{i_j}x_{i_{j+1}}]_{M} \bigcup \cup [x_{i_k}x_{i_{k+1}}]_{\partial C_{M}}$ is a path in $M$ from $x$ to $y$.

Therefore
\begin{align*}
|xy|_M &\leq \sum |x_{i_j}x_{i_{j+1}}|_M + \sum |x_{i_k}x_{i_{k+1}}|_{\partial C_{M}}\\
&\leq \sum |x_{i_j}x_{i_{j+1}}|_{M} + \sum \left( |x_{i_k}x_{i_{k+1}}|_{C_{M}}+C |x_{i_k}x_{i_{k+1}}|_{C_{M}}^{3}  \right) \quad \text{ by arc/chord comparison in }C_{M},\\
& \quad \text{ where $C=C(K_{C}^{+},{\lambda}_{C}^{+}) = C(K^{+},{\lambda}^{\pm})$,
since $|x_{i_{k}}x_{i_{k+1}}|_{C_{M}} = |x_{i_{k}}x_{i_{k+1}}|_{\widetilde{M}} \leq |xy|_{\widetilde{M}} \leq 1$
}\\
&= \left( \sum |x_{i_j}x_{i_{j+1}}|_M + \sum |x_{i_k}x_{i_{k+1}}|_{C_{M}} \right) + \sum C |x_{i_k}x_{i_{k+1}}|_{C_{M}}^{3}\\
&\leq |xy|_{\widetilde{M}}+C|xy|_{\widetilde{M}}^{3}
\end{align*}

The proof of the second statement of the lemma is analogous:\\

$|xy|_{\widetilde{M}} \leq |xy|_{\partial M}$ for all $x,y\in \partial M$ since $\partial M\subseteq \widetilde{M}$.\\

Suppose $|xy|_{\widetilde{M}} \leq 1$, and also $|xy|_{\widetilde{M}} < CAT_{k_{0}}Rad(M), CAT_{k_{0}}Rad(C_{M})$, where $curv M, curv C_{M} \leq k_{0}$.\\
Let $[xy]_{\widetilde{M}}=\cup [x_{i_j}x_{i_{j+1}}]_{M} \bigcup \cup [x_{i_k}x_{i_{k+1}}]_{C_{M}}$, where $x_{i_k} \in \partial C_{M}$. Note that $\cup [x_{i_j}x_{i_{j+1}}]_{\partial M} \bigcup \cup [x_{i_k}x_{i_{k+1}}]_{\partial C_{M}}$ is a path in $\partial M$ from $x$ to $y$. Therefore
\begin{align*}
|xy|_{\partial M} &\leq \sum {|x_{i_j}x_{i_{j+1}}|}_{\partial M} + \sum |x_{i_k}x_{i_{k+1}}|_{\partial C_{M}}\\
&\leq \sum \left( |x_{i_j}x_{i_{j+1}}|_{M} +  C|x_{i_j}x_{i_{j+1}}|_{M}^{3} \right) + \sum \left( |x_{i_k}x_{i_{k+1}}|_{C_{M}}+C |x_{i_k}x_{i_{k+1}}|_{C_{M}}^{3}  \right)\\
& \quad \text{ by arc/chord comparison in }M \text{ and } C_{M},
\text{ where
$C=C(K^{\pm},{\lambda}^{\pm},i_{0})$
 }\\
&= \left( \sum |x_{i_j}x_{i_{j+1}}|_M + \sum |x_{i_k}x_{i_{k+1}}|_{C_{M}} \right) + C \left( \sum |x_{i_j}x_{i_{j+1}}|_{M}^{3}
+ \sum |x_{i_k}x_{i_{k+1}}|_{C_{M}}^{3} \right)\\
&\leq |xy|_{\widetilde{M}}+C|xy|_{\widetilde{M}}^{3} \qedhere
\end{align*}
\end{proof}

\vspace{1cm}

\begin{proof}[Proof of Lemma~\ref{sublem:4}]
Let $\theta$ denote the angle opposite side $c$ in a model triangle $\Delta$ with sidelengths $a,b,c$.\\
Let $\theta_C$ denote the angle opposite side $c+Cc^3$ in a model triangle $\Delta_C$ (in the same model space) with sidelengths $a,b, c+Cc^3$.
For simplicity suppose the model space is the Euclidean plane $\mathbb{R}^2$. By the law of cosines,
\begin{align*}
|\cos(\theta)-\cos({\theta}_{C})| = \frac{a^2+b^2-c^2}{2ab} - \frac{a^2+b^2-{(c+Cc^3)}^2}{2ab} = Cc^3\frac{2c+Cc^3}{2ab}
\end{align*}
which implies
\begin{align*}
|\theta-{\theta}_{C}| \leq \tau_{a,b}(Cc^3)%
\end{align*}

\noindent
(and similarly for the other models $S^2, H^2$).

Obviously,
$normalized \ excess(\Delta_C):=\dfrac{a+b-(c+Cc^3)}{\min \{ a, b \} }\longrightarrow 0$ implies $\theta_C \uparrow \pi$.
So 
if $\dfrac{a+b-(c+Cc^3)}{\min \{ a, b \} } \leq \tau_{a,b,c}(Cc^3)$ then
$\theta_C \geq \pi - \tau(normalized \ excess(\Delta_C)) = \pi - \tau(\tau_{a,b,c}(Cc^3)) = \pi - \tau_{a,b,c}(Cc^3)$.

Therefore
$\theta \geq \theta_C - \tau_{a,b}(Cc^3)
\geq (\pi -\tau_{a,b,c}(Cc^3)) - \tau_{a,b}(Cc^3)
= \pi - \tau_{a,b,c}(Cc^3)$.
\end{proof}

\vspace{1cm}

\begin{lemma}[\protect{\cite[lemma 10.8.13]{BBI:01}}]\label{lem:3}
Suppose $curv X \geq k$, $x \in X$ and $a_{-1},a_{1} \in X$ with ${\angle}a_{-1}xa_{1} \geq \pi - \delta$.  If $b\in X$ is any point satisfying $|xb| < \frac{\epsilon}{4} \min \{ |xa_{-1}|, |xa_{1}| \}$ then
\begin{align*}
0 < {\angle}a_{1}xb - {\overline{\angle}}a_{1}xb < 2 \max \{ \epsilon, \delta \}.
\end{align*}
\end{lemma}

\vspace{1cm}

The following proposition is a CBB dual to \cite[Section 6.1]{Lyt:01}.
Recall that a subspace $Z$ of a length space $X$ is $(C,2)$-convex if their intrinsic metrics satisfy $d_{Z} \leq d_{X}+Cd_{X}^3$, at least locally.\\

If $Z$ is any subspace of $X$, with an intrinsic metric, then
it is usually the case that
${\angle}^{Z}(v,w)\geq {\angle}^{X}(v,w)$ for any $v,w\in {\Sigma}_{z}Z$, $(z\in Z)$, whenever ${\angle}^{X}$ and ${\angle}^{Z}$ exist.

Conversely,\\

\noindent
\begin{proposition}\label{prop:4.2}
Suppose $curv X\geq k$. Let $Z\subset X$ be $(C,2)-$convex, for some constant $C\geq 0$. Then ${\angle}^{Z}(\gamma_1,\gamma_2)={\angle}^{X}(\gamma_1,\gamma_2)$ for all $Z-$geodesics $\gamma_1$ and $\gamma_2$.
\end{proposition}

\begin{proof}
Since a $Z-$geodesic will again be $(C,2)-$convex, it may be supposed $Z=\gamma$. Set $x_{t}=\gamma(t)$.
Let $n$ denote the dimension of $X$.

First suppose that $x_{0}$ is a non-singular point of $X$ in the sense that $(n,\delta)$-$str.rad(x_{0}) > 0$ for all $\delta > 0$.
For the moment, let $\delta > 0$ be fixed but arbitrary. Then $(n,\delta)$-$str.rad(x_0) > 0$ by definition.
There exists (for $t$ sufficiently small, depending on $\delta$, e.g., $t < (n,\delta)$-$str.rad(x_{0})$) a point $y_t\in X$ such that $|x_{0}y_{t}|_X = |x_{0}x_{t}|_{X}$ and ${\angle}^{X}x_{t}x_{0}y_{t} \geq \pi - \delta$.
Consider the triangle $\Delta x_{0}x_{f(t)}x_{t}$, where we choose $f(t)=t^2$.
Let $\epsilon:=8t$. For sufficiently small $t$,
\begin{align}\label{ineq:4.1}
|x_{0}x_{f(t)}|_{X}<\frac{\epsilon}{4}\min \{ |x_{0}x_{t}|_{X}, |x_{0}y_{t}|_{X}\}
\end{align}
since $|x_{0}x_{f(t)}|_{X}\leq |x_{0}x_{f(t)}|_{Z}=f(t)=t^2$ and 
$|x_{0}x_{t}|_{X}=|x_{0}y_{t}|_{X}\geq |x_{0}x_{t}|_{Z} - C |x_{0}x_{t}|_{Z}^{3} = t - Ct^3$ and
$t^2<\frac{\epsilon}{4}(t-Ct^3)$ when $\epsilon=8t$ and $t$ sufficiently small ($t< \sqrt{\frac{1}{2C}}$).

Inequality (\ref{ineq:4.1}) implies
\begin{align}\label{ineq:4.2}
0<{\angle}^{X}x_{t}x_{0}x_{f(t)} - {\overline{\angle}}^{X}x_{t}x_{0}x_{f(t)} < 2\max \{ \epsilon, \delta\}
\end{align}

by (\cite{BBI:01}, lemma 10.8.13).

Note that $d_{X}\leq d_{Z} \leq d_{X}+Cd_{X}^3$ implies $d_{Z}-Cd_{Z}^{3} \leq d_{X}$.

Letting $a=|x_{0}x_{f(t)}|_{X}, b=|x_{0}x_{t}|_{X}, c=|x_{f(t)}x_{t}|_{X}$,
\begin{align*}
f(t)-Cf(t)^3 &\leq a \leq f(t)\\
t-Ct^3 &\leq b \leq t\\
(t-f(t))-C{(t-f(t))}^3 &\leq c \leq t-f(t)
\end{align*}

Then, assuming $k<0$,
\begin{align}\label{ineq:4.3}
{\overline{\angle}}^{X}x_{t}x_{0}x_{f(t)} &\leq {\angle}^{X}x_{t}x_{0}x_{f(t)}\notag\\
&\leq {\overline{\angle}}^{X}x_{t}x_{0}x_{f(t)} + 2\max \{ \epsilon, \delta\} \quad \text{ by (\ref{ineq:4.2})}\notag\\
&= {\cos}^{-1}\left( \dfrac{\cosh(\sqrt{|k|}a)\cosh(\sqrt{|k|}b)-\cosh(\sqrt{|k|}c)}{\sinh(\sqrt{|k|}a)\sinh(\sqrt{|k|}b)} \right) +2\max \{ \epsilon, \delta\}  \notag \\
&\leq {\cos}^{-1}\left(  \dfrac{\cosh(\sqrt{|k|}(f(t)-C{f(t)}^3))\cosh(\sqrt{|k|}(t-Ct^3))-\cosh(\sqrt{|k|}(t-f(t)))}{\sinh(\sqrt{|k|}f(t))\sinh(\sqrt{|k|}t)} \right) + 2\max \{ \epsilon, \delta\}\notag\\
&\leq \sqrt{2Ct}+O(t^{3/2})+2\max \{ 8t, \delta\} 
\end{align}
whenever $t < \min \{ (n,\delta)$-$str.rad(x_{0}),  \sqrt{\frac{1}{2C}}  \}$. $ \ $
(Note that $(n,\delta)$-$str.rad(x_{0})$ is non-increasing as $\delta$ tends to $0$.)

This implies that $\gamma$ has a unique initial direction $v\in \Sigma_{x_0}X$ and that the angle between $v$ and $[x_{0}x_{t}]_{X}^{'}$ is at most a constant (depending only on $C$ and $k$) times a power of $t$, plus a constant involving $\delta$. This can be taken arbitrarily small, since $\delta > 0$ was arbitrary.

 By the triangle inequality for angles, this implies that the angle between two $Z-$geodesics equals the angle between these two curves in $X$.\\

Now suppose $x_{0}$ is a singular point of $X$.

Let
\begin{align*}
x_{0}^{i} &\mapsto x_{0}\\
x_{t}^{i} &\mapsto x_{t}\\
x_{f(t)}^{i} &\mapsto x_{f(t)}
\end{align*}
where each $x_{0}^{i} \in X$ is a non-singular point of $X$. Such a sequence exists since the set $S_{X}$ of singular points of $X$ has Hausdorff dimension $dim_{\mathcal{H}}(S_{X}) \leq n-1$ and hence is nowhere dense (\cite{OS:94}).

Since $curv X \geq k$, angles are lower semi-continuous:
\begin{align*}
{\angle}^{X}x_{t}x_{0}x_{f(t)} \leq \underset{i \longrightarrow \infty}{\liminf} {\angle}^{X} x_{t}^{i}x_{0}^{i}x_{f(t)}^{i}
\end{align*}
and the right-hand side may be bounded above by (\ref{ineq:4.3}), just as in the previous case. Again we conclude that the angle between two $Z-$geodesics equals the angle between these two curves in $X$.%
\end{proof}

\vspace{0.5cm}

\begin{corollary}
Let $X$ and $Z$ be as in the proposition, with $(n,\delta)$-$str.rad(X)\geq i_0 > 0$ (for some $\delta < \pi / 4$, say; in particular, $X$ is weakly geodesically extendible). Then $Z$ has no $C^1$-smoothly closed geodesic of length less than $i_1=i_1(C,k,i_0)>0$.
\end{corollary}

\begin{proof}
This follows similarly as in~\cite{Lyt:01}. 
Let $\gamma : [0,2t] \longrightarrow Z$ be a $C^{1}$-smoothly closed geodesic loop in $Z$. Let $\eta$ be a minimal $X$-geodesic from $\gamma(0)$ to $\gamma(t)$. Then
\begin{align*}
\pi &= {\angle}^{Z}({\gamma}'(0),{\gamma}'(2t))\\
&= {\angle}^{X}({\gamma}'(0),{\gamma}'(2t)) \quad \text{ by Proposition~\ref{prop:4.2}}\\
&\leq {\angle}^{X}({\gamma}'(0), \eta') + {\angle}^{X}(\eta',{\gamma}'(2t))\\
&\leq (c t^{\alpha} + 2\delta) + (c t^{\alpha} + 2\delta) \quad \text{ for some constants } c=c(k,C) \text{ and } \alpha, \text{ by (\ref{ineq:4.3})}
\end{align*}
which implies
\begin{align*}
t \geq {\left( (\pi - 4\delta ) / 2c \right)}^{1/\alpha} > 0
\end{align*}
if $t < \min \{ i_{0}, \sqrt{\frac{1}{2C}}, \frac{\delta}{8} \}$, as desired.

It is easy to check that these constants $c$ and $\alpha$ may be chosen uniformly, if $t$ is sufficiently small relative to $k$ and $C$. Take $\alpha:=\frac{1}{2}$. 
Using the crude estimates
\begin{align*}
1 + \frac{x^2}{2} \leq \cosh(x) &\leq 1 + \frac{x^2}{2} + x^4\\
\sinh(x) &\leq x + x^3\\
\cos(x) &\leq 1 -  \frac{x^2}{2} + \frac{x^4}{24},
\end{align*}
which are valid for all $x$,
one has
\begin{align*}
\cos({\overline{\angle}}^{X}x_{t}x_{0}x_{f(t)})
&=  \dfrac{\cosh(\sqrt{|k|}(f(t)-C{f(t)}^3))\cosh(\sqrt{|k|}(t-Ct^3))-\cosh(\sqrt{|k|}(t-f(t)))}{\sinh(\sqrt{|k|}f(t))\sinh(\sqrt{|k|}t)} \notag\\
&\geq \cos(const\cdot t^{1/2}) \notag
\end{align*}
or
\begin{align*}
{\overline{\angle}}^{X}x_{t}x_{0}x_{f(t)} \leq const \cdot t^{1/2}
\end{align*}

\noindent
if $k$ (which may be assumed negative) is normalized to $k=-1$, 
$const$ is chosen as $const:=\sqrt{2(C+2)}$, and $t$ is taken sufficiently small (depending only on $C$, and implicitly, on $k$) such that the following
eight
inequalities hold:
\begin{align*}
\left[ {3-\frac{1}{6}{(C+2)}^2} \right] t \geq -\frac{1}{20}\\
\left[ {\frac{C^2}{2}+C-\frac{15}{4}} \right] t^2 \geq -\frac{1}{20}\\
\left[ {3-\frac{1}{6}{(C+2)}^2} \right] t^3 \geq -\frac{1}{20}\\
\left[ 1-\frac{C}{2} \right] t^4 \geq -\frac{1}{20}\\
\left[ -1-\frac{1}{6}{(C+2)}^2 \right] t^5 \geq -\frac{1}{20}\\
\left[ -\frac{1}{6}{(C+2)}^2 \right]t^7 \geq -\frac{1}{20}\\
\left[ \frac{C^2}{4}-\frac{C^3}{2} \right] t^{10} \geq -\frac{1}{20}\\
\left[ -\frac{C^3}{2} \right] t^{12} \geq -\frac{1}{20}
\end{align*}\qedhere
\end{proof}

\vspace{0.4cm}

\begin{proposition}\label{prop:1}
Let $Z \subseteq X$. Suppose $curv X\geq k$, $curv Z \geq k$, $X$ is geodesically extendible, and $d_{Z} \leq d_{X} + Cd_{X}^{3}$ locally. Then
\begin{align*}
{\angle}^{X}([xy]_{X},[xy]_{Z}) \leq \tau_{C}(|xy|_{X})
\end{align*}
for all $x,y\in Z$.
\end{proposition}

\begin{proof}
This follows from Proposition~\ref{prop:4.2}.
Since $X$ is assumed geodesically extendible, the argument there involving $(n,\delta)$-strain radius simplifies.
\end{proof}

\vspace{0.8cm}

\begin{corollary}\label{cor:1}
Let $y,z$ be as in Lemma~\ref{sublem:2}.
For any $u,v \in [yz]_{\widetilde{M}} \cap \partial M$,
\begin{align*}
{\angle}^{\widetilde{M}} ( [uv]_{\widetilde{M}}, [uv]_{\partial M}) \leq \tau(|uv|_{\widetilde{M}})
\end{align*}
\end{corollary}

\begin{proof}
Immediate from Lemma~\ref{sublem:3} and Proposition~\ref{prop:1}.
\end{proof}

\vspace{1cm}

\subsubsection{Proof of Lemma~\ref{sublem:2}}
\label{ss:A1.3}

Note that
$curv \widetilde{M} \geq const$ by construction, and $M \subset \widetilde{M}$ is $(C,2)$-convexly embedded by Lemma~\ref{sublem:3}, so Proposition~\ref{prop:4.2} is available for employment in the proof of Lemma~\ref{sublem:2} which follows (although only the one-way inequality of the remark immediately preceding Proposition~\ref{prop:4.2} is needed, and even then only in the special situation where the spaces are smooth manifolds and the submanifold (with smooth boundary) is of the same dimension as its ambient manifold.) \\

\begin{proof}[Proof of Lemma~\ref{sublem:2}]
Let $x,y,z \in M$, with $y\in \partial M$, $[xy]_{M}$ a minimizing segment orthogonal to $\partial M$, and $x,z\in B(y,\frac{1}{2}CAT_{k}Rad(M); M)$.\\
Note that $[xy]_M = [xy]_{\widetilde{M}}$.\\
It may be assumed that $[yz]_{\widetilde{M}}\nsubseteq M$ as otherwise $[yz]_{M}=[yz]_{\widetilde{M}}$, and ${\angle}^{M} xyz = {\angle}^{\widetilde{M}} xyz$ would hold trivially.

\begin{align}\label{ineq:13}
|{\angle}^{M} xyz - {\angle}^{\widetilde{M}} xyz|
&= |{\angle}^{M}([xy]_M,[yz]_M) - {\angle}^{\widetilde{M}}([xy]_{\widetilde{M}},[yz]_{\widetilde{M}})| \quad \text{by definition} \notag \\
&= |{\angle}^{\widetilde{M}}([xy]_M,[yz]_M) - {\angle}^{\widetilde{M}}([xy]_M,[yz]_{\widetilde{M}})| \quad \text{ by Proposition~\ref{prop:4.2}} \notag \\
&\leq {\angle}^{\widetilde{M}}([yz]_M,[yz]_{\widetilde{M}}) \quad \text{ by triangle inequality}
\end{align}

One needs to relate quantities in $M$ and $\widetilde{M}$, such as angles, via quantities in $\partial M$. One obstacle is that, while $z$ lies in $M$, $z$ does not necessarily lie in $\partial M$. To overcome this,
let $w :=$ last point of $[yz]_{\widetilde{M}} \cap \partial M$ such that $[wz]_{\widetilde{M}} \subseteq M$ (so $[wz]_{\widetilde{M}} = [wz]_{M}$).

We may assume $w \neq y$.

Possibly $w=z$.\\

If $w = z$, then $z\in \partial M$. So for the angle at $y$
\begin{align}
{\angle}^{\widetilde{M}} ([yz]_{M},[yz]_{\widetilde{M}})
&\leq {\angle}^{\widetilde{M}} ([yz]_{M},[yz]_{\partial M})
+ {\angle}^{\widetilde{M}} ([yz]_{\partial M},[yz]_{\widetilde{M}})
\quad \text{ by triangle inequality} \notag \\
&= {\angle}^{M} ([yz]_{M},[yz]_{\partial M})
+ {\angle}^{\widetilde{M}} ([yz]_{\partial M},[yz]_{\widetilde{M}}) \notag \\
&\leq \tau(|yz|_{M})
+ {\angle}^{\widetilde{M}} ([yz]_{\partial M},[yz]_{\widetilde{M}})
\quad \text{ by arc/chord comparison in $M$} \notag \\
&\leq \tau(|yz|_{M})
+ \tau(|yz|_{\widetilde{M}})
\quad \text{ by Corollary~\ref{cor:1}} \notag \\
&= \tau(|yz|_{\widetilde{M}}) \quad \text{ by Lemma~\ref{sublem:3} }\label{ineq:18} 
\end{align}

\vspace{0.3cm}

If $w \neq z$, then
at $y$
\begin{align}\label{ineq:12}
{\angle}^{\widetilde{M}}([yz]_{M},[yz]_{\widetilde{M}}) 
&\leq {\angle}^{\widetilde{M}}([yz]_{M},[yw]_{M}) 
    + {\angle}^{\widetilde{M}}([yw]_{M},[yz]_{\widetilde{M}}) \quad \text{ by triangle inequality} \notag\\
&= {\angle}^{M} ([yz]_{M},[yw]_{M}) 
 + {\angle}^{\widetilde{M}}([yw]_{M},[yw]_{\widetilde{M}}) \quad \text{ since $[yw]_{\widetilde{M}} \subseteq [yz]_{\widetilde{M}}$} \notag \\
&\leq {\angle}^{M}([yz]_{M},[yw]_{M}) 
    + {\angle}^{M}([yw]_{M},[yw]_{\partial M}) 
    + {\angle}^{\widetilde{M}}([yw]_{\partial M},[yw]_{\widetilde{M}})
\quad \text{ by triangle inequality } \notag\\
&\leq (*) + \tau(|yw|_{M}) + {\angle}^{\widetilde{M}}([yw]_{\partial M},[yw]_{\widetilde{M}})
\quad \text{ by arc/chord comparison in $M$} \notag \\
&\leq (*) + \tau(|yw|_{M}) + \tau(|yw|_{\widetilde{M}})
\quad \text{by Corollary~\ref{cor:1}} \notag \\
&\leq (*) + \tau(|yz|_{\widetilde{M}})
\quad \text{}
\end{align}

In order to estimate (\ref{ineq:13}) from above, it suffices to prove that $(*):={\angle}^{M}zyw \leq \tau(|yz|_{\widetilde{M}})$. First we will show that ${\angle}^{M}ywz \longrightarrow \pi$ as $|yz|_{\widetilde{M}} \longrightarrow 0$.

\begin{align*}
\pi &= {\angle}^{\widetilde{M}}ywz
\quad \text{ since $w \in [yz]_{\widetilde{M}}$}\\
&= {\angle}^{\widetilde{M}}([yw]_{\widetilde{M}},[wz]_{\widetilde{M}})
\quad \text{ by definition }\\
&\leq {\angle}^{\widetilde{M}}([yw]_{\widetilde{M}},[wy]_{M}) + 
{\angle}^{\widetilde{M}}([wy]_{M},[wz]_{\widetilde{M}})\\
&\leq \tau(|yw|_{\widetilde{M}}) + {\angle}^{\widetilde{M}}([wy]_{M},[wz]_{M})
\quad \text{ by the same argument as in (\ref{ineq:18}) above}\\
&= \tau(|yw|_{\widetilde{M}}) + {\angle}^{M}([wy]_{M},[wz]_{M})
\quad \text{ by Lemma~\ref{sublem:3} and Proposition~\ref{prop:4.2}}
\end{align*}
Thus
\begin{align}\label{ineq:11}
{\angle}^{M}ywz \geq \pi - \tau(|yw|_{\widetilde{M}}).
\end{align}

In $M$ consider the triangle ${\Delta}^{M} ywz$.
Since $curv M \leq k_{M}^{+}$ (wlog $k_{M}^{+} > 0$) and (\ref{ineq:11}) holds,
\begin{align*}
{\overline{\angle}}^{M} ywz = {\overline{\angle}}^{M}([yw]_{M},[wz]_{M}) \geq {\angle}^{M}([yw]_{M},[wz]_{M}) \longrightarrow \pi
\end{align*}
as $|yz|_{\widetilde{M}} \longrightarrow 0$.

By the law of sines in the model space $M_{k_{M}^{+}}^{2}$,
\begin{align*}
\sin {\overline{\angle}}^{M} zyw
&= \frac{\sin\left( \sqrt{k_{M}^{+}} \, |wz|_{M}\right) }{ \sin\left( \sqrt{k_{M}^{+}} \, |yz|_{M} \right) } \sin {\overline{\angle}}^{M} ywz\\
&\leq
\frac{\sin\left( \sqrt{k_{M}^{+}} \, |wz|_{\widetilde{M}}\right) }{ \sin\left( \sqrt{k_{M}^{+}} \, |yz|_{\widetilde{M}} \right) }
\sin {\overline{\angle}}^{M} ywz \quad \text{ since $|wz|_{M}=|wz|_{\widetilde{M}}$, $|yz|_{\widetilde{M}} \leq |yz|_{M}$ }\\
&\leq \sin {\overline{\angle}}^{M} ywz \quad \text{ since $|wz|_{\widetilde{M}} \leq |yz|_{\widetilde{M}}$}\\
& \quad \longrightarrow 0
\end{align*}
which implies that
at $y$,
${\overline{\angle}}^{M}zyw = {\overline{\angle}}^{M}([yz]_{M},[yw]_{M}) \longrightarrow 0$. (sine can also tend to zero when its argument tends to $\pi$, but it is impossible for ${\overline{\angle}}^{M}zyw$ and ${\overline{\angle}}^{M}ywz$ to both tend to $\pi$, since the corresponding comparison triangle for ${\Delta}^{M}ywz$ will be small compared to $\frac{\pi}{\sqrt{k_{M}^{+}}}$.)

Again, since $curv M \leq k_{M}^{+}$, this forces
\begin{align}\label{ineq:16}
(*) = {\angle}^{M}([yz]_{M},[yw]_{M}) \longrightarrow 0.
\end{align}

This completes the proof of Lemma~\ref{sublem:2}.
\end{proof}

\vspace{1.5cm}

\subsection[Proof of Theorem~\ref{thm:Top}]{Proof of Theorem~\ref{thm:Top}}\label{s:5iii}

\vspace{0.2cm}
The proof of Theorem~\ref{thm:Top} will use the following proposition on limits of gluings.

\newpage

\begin{proposition}\label{prop:8.4}
Assume $Z_i := X_{i} \underset{A_i}{\cup} Y_{i}$, where

$X_i$ are path-connected length spaces, 

$Y_{i} = \coprod_{k} {\left( Y_i \right)}_{k}$ is a disjoint union of path-connected length spaces $Y_i$

$X_i \overset{GH}{\longrightarrow} pt$

$Y_i \overset{GH}{\longrightarrow} Y$ via $\epsilon_i$-Hausdorff approximation $f : Y_i \longrightarrow Y$

$Y = \coprod_{k} {\left( Y \right)}_{k}$

\quad \quad (so ${\left( Y_{i} \right)}_{k} \overset{GH}{\longrightarrow} Y_{k}$ via $\epsilon_i$-Hausdorff approximation given by restriction, $f|_{{(Y_i)}_{k}}$)

$A_i \subset X_i, Y_i$ is a closed subset of both $X_i$ and $Y_i$, where

$A_{i} = \coprod_{k} {\left( A_i \right)}_{k}$ with each $ {\left( A_i \right)}_{k}=A_i \cap (Y_i)_k$ non-empty and path-connected

$A_i \longrightarrow A$ as subsets (under $Y_i \overset{GH}{\longrightarrow} Y$)

\quad \quad (so in particular, $A \subseteq Y$, where $A \subseteq Y$ is closed)

\quad \quad where $A = \coprod_{k} {\left( A \right)}_{k}$, each ${\left( A \right)}_{k}$ path-connected

Then $Z_i \overset{GH}{\longrightarrow} Z:=Y/A$
where $Z$ carries the quotient metric.
\end{proposition}

\vspace{0.5cm}

\begin{figure}[h]
\begin{center}
\begin{overpic}[scale=.6]{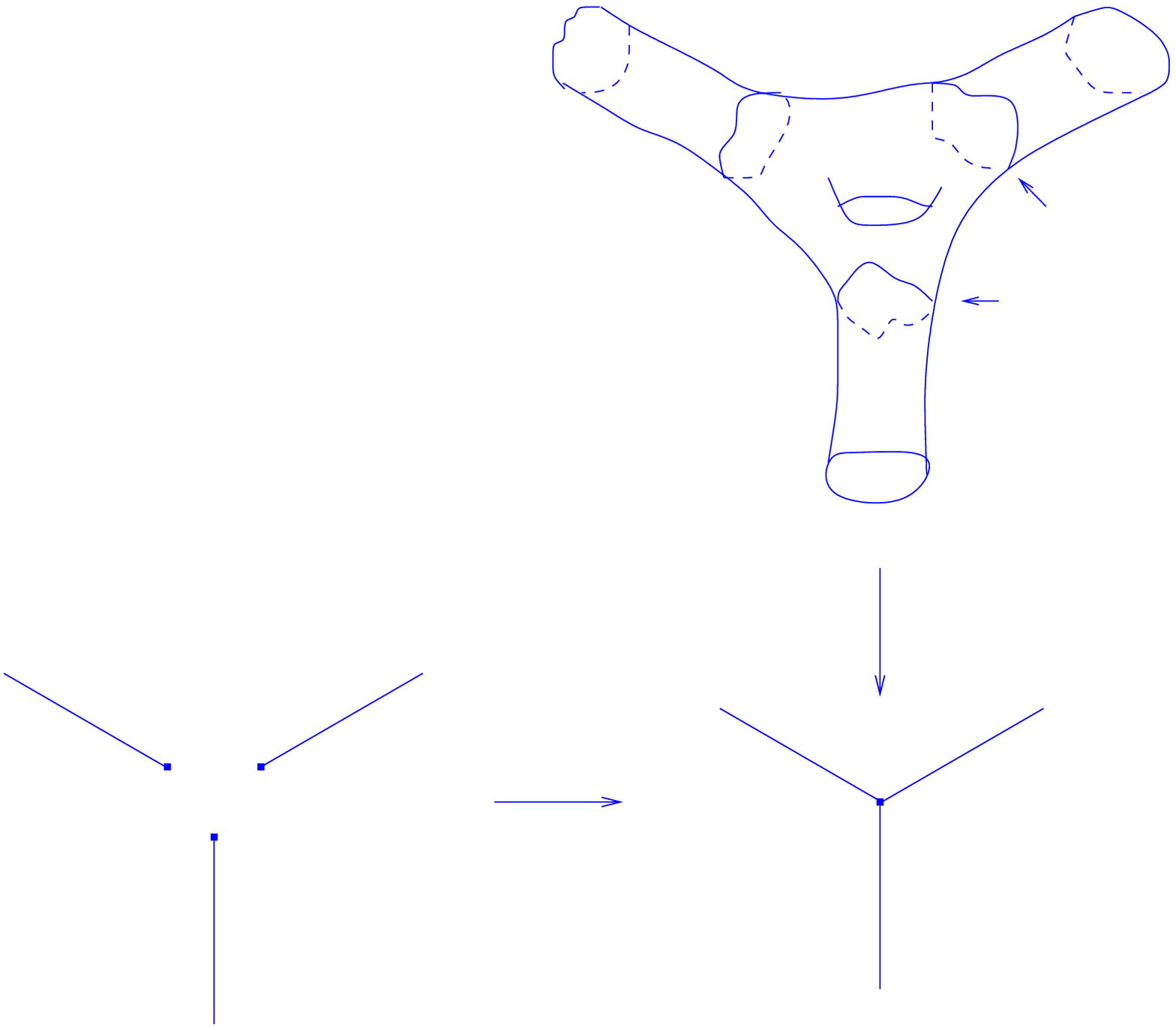}
\put(73,75){$X_i$}
\put(85,60){${(A_i)}_j$}
\put(89,66){${(A_i)}_k$}
\put(82,51){${(Y_i)}_j$}
\put(94,71){${(Y_i)}_k$}
\put(96,16){$Z$}
\put(19,14){$A_j$}
\put(22,18){$A_k$}
\put(19,8){$Y_j$}
\put(29,23){$Y_k$}
\put(103,63){$Z_i$}
\put(46,22){$\pi$}
\end{overpic}
\caption[]{}
\label{fig:figname}
\end{center}
\end{figure}

\begin{proof}

Let $\pi : Y \longrightarrow Y/A$ be the quotient map. 
Note that
$\pi(A)$, the image of $A$ under $\pi$, is a single point $\{ pt \}$ in $Z$.

Define a map 
$F : Z_{i} \longrightarrow Z$ by
$F(z):=\begin{cases}
\pi f(z) & z\in Y_{i}\\
\pi(A) & z\in X_{i} \setminus A_i
\end{cases}$

Will prove that $F$ is a Hausdorff approximation.

By definition of quotient metric,
since ${(A)}_{j}$ and ${(A)}_{k}$ are path-connected by assumption,
\begin{align}\label{eq:8.2}
|F(z)F(w)|_{Y/A} = \min \{ |f(z)f(w)|_{Y}, |f(z)A|_{Y}+|Af(w)|_{Y} \}
\end{align}

for any $z,w\in Y_i$.

Here, the distance from a point to a set is defined as
$|pA| = \underset{q \in A}{\inf}|pq|$.

Since $A_i \longrightarrow A$ as subsets, $d_{H}^{Y}(f(A_i),A) \leq \tau(\frac{1}{i})$.
So for $p\in Y$,
\begin{align}\label{ineq:8.4.subs}
||pA|_{Y} - |pf(A_i)|_{Y}| \leq \tau(\frac{1}{i})
\end{align}

If $z \in Y_i$, then
\begin{align}
|zA_{i}|_{Y_i} &\leq |zw|_{Y_i} \notag \\
&\leq |f(z)f(w)|_{Y}+\epsilon_i \notag \\
&= |f(z)f(A_i)|_{Y}+\epsilon_i \notag \\
&\leq |f(z)A|_{Y}+\tau(\frac{1}{i}) + \epsilon_i \label{ineq:8.4.1}
\quad \text{ by (\ref{ineq:8.4.subs})}
\end{align}
where $w\in A_i$ is chosen to make the equality in the third line above hold (possible since $f(A_i)$ may be assumed closed)

On the other hand,
\begin{align}
|f(z)A| &\leq |f(z)f(A_i)|+\tau(\frac{1}{i})
\quad \text{ by (\ref{ineq:8.4.subs}) } \notag \\
&\leq |f(z)f(w)|+\tau(\frac{1}{i}) \notag \\
&\leq |zw|+\epsilon_i+\tau(\frac{1}{i}) \notag \\
&= |zA_i|+\epsilon_i+\tau(\frac{1}{i})\label{ineq:8.4.2}
\end{align}
where the point $w\in A_i$ here is chosen to make the equality in the fourth line above hold (possible since $A_i$ is closed)

Combining (\ref{ineq:8.4.1}) and (\ref{ineq:8.4.2}) yields
\begin{align}\label{ineq:8.4.subsii}
||zA_i|_{Y_i} - |f(z)A|_{Y} | \leq \epsilon_i + \tau(\frac{1}{i})
\end{align}

\vspace{0.5cm}

There are four cases to consider.

\underline{Case $z,w\in X_{i} \setminus A_{i}$} :

Note that
$|zw|_{Z_i} \leq |zw|_{X_i} \leq d(X_i) \longrightarrow 0$.

Hence
\begin{align}\label{ineq:8.5}
|F(z)F(w)|_{Z} - |zw|_{Z_i} &= |\pi(A)\pi(A)|_{Z} - |zw|_{Z_i} \notag\\
&\leq 0
\end{align}

and
\begin{align}
|zw|_{Z_i} - |F(z)F(w)|_{Z} &\leq |zw|_{Z_i} \notag\\
&\leq d(X_i) \notag\\
& \quad \longrightarrow 0
\end{align}

\underline{Case $z\in {\left( Y_i \right)}_{j}, w\in  {\left( Y_i \right)}_{k}$,  
$j \neq k$} :

Observe that since $z$ and $w$ lie in different path-components of $Y_i$, any $Z_i$ geodesic from $z$ to $w$ must pass through $X_i$ via $A_i$, so
\begin{align*}
|zw|_{Z_i}
&\geq |zA_i|_{ {\left( Y_i \right)}_{j}}  + |A_{i}w|_{ {\left( Y_i \right)}_{k}}\\
&\geq \left( |f(z)A|_{ {\left( Y \right)}_{j} } - \epsilon_{i} -\tau(\frac{1}{i})\right)  +  \left( |Af(w)|_{ {\left( Y \right)}_{k} } - \epsilon_{i} -\tau(\frac{1}{i}) \right) \\
& \quad \text{ by (\ref{ineq:8.4.subsii}), since $f$ an $\epsilon_{i}$-Hausdorff approximation }
\end{align*}

Hence
\begin{align}
|F(z)F(w)|_{Z} - |zw|_{Z_i}
&= |F(z)F(w)|_{Y/A} - |zw|_{Z_i} \notag\\
&\leq |f(z)A|_{Y} + |Af(w)|_{Y} - |zw|_{Z_i} \quad \text{ by (\ref{eq:8.2})}\notag\\
&\leq |f(z)A|_{Y} + |Af(w)|_{Y} - \left( \left( |f(z)A|_{ {\left( Y \right)}_{j} }  - \epsilon_{i} -\tau(\frac{1}{i}) \right) + \left( |Af(w)|_{ {\left( Y \right)}_{k}} -\epsilon_{i} -\tau(\frac{1}{i}) \right) \right) \notag\\
& \quad \text{ by the observation above } \notag\\
&= 2\epsilon_{i} + 2 \tau(\frac{1}{i}) \longrightarrow 0
\end{align}

(Note that $|f(z)A|_{Y}=|f(z)A|_{{(Y)}_{j}}$ since $z\in {(Y_i)}_{j}$ implies $f(z)\in {(Y)}_{j}$.)

For the converse,
\begin{align}
|zw|_{Z_i} - |F(z)F(w)|_{Z}
&\leq \left( d(X_i) + d^{ {\left( Y_i \right)}_{j} } (z,A_i) + d^{ { {\left( Y_i \right)}_{k}} } (A_i,w) \right) - |F(z)F(w)|_{Y/A} \quad \text{ by triangle inequality } \notag\\
&\leq d(X_i) + \left( d^{  { {\left( Y \right)}_{j}} } (f(z),A) + \epsilon_{i} + \tau(\frac{1}{i}) \right)  + \left( d^{  { {\left( Y \right)}_{k}} }(A,f(w)) + \epsilon_{i} + \tau(\frac{1}{i}) \right) - |F(z)F(w)|_{Y/A} \notag\\
& \quad \text{ by (\ref{ineq:8.4.subsii}), since $f$ an $\epsilon_i$-Hausdorff approximation} \notag\\
&= d(X_i) + \left( |f(z)A|_{{(Y)}_{j}} + \epsilon_{i} + \tau(\frac{1}{i}) \right)  + \left( |Af(w)|_{{(Y)}_{k}} + \epsilon_{i} + \tau(\frac{1}{i}) \right) - \left( |f(z)A|_{Y} + |Af(w)|_{Y} \right)  \notag\\
& \quad \text{ by (\ref{eq:8.2}), since $f(z)$ and $f(w)$ in different path-components of $Y$,} \notag\\
& \quad \text{ hence $|f(z)f(w)|_{Y}=\infty$ } \notag\\
&= d(X_i) + 2\epsilon_{i} + 2\tau(\frac{1}{i}) \notag\\
& \quad \longrightarrow 0
\end{align}

\underline{Case $z\in  {\left( Y_i \right)}_{j}, w\in {\left( Y_i \right)}_{j}$} : 

\begin{align}
|F(z)F(w)|_{Z} - |zw|_{Z_i}
&= |F(z)F(w)|_{Y/A} - |zw|_{Z_i} \notag\\
&= \min \{ |f(z)f(w)|_{Y}, |f(z)A|_{Y}+|Af(w)|_{Y} \} - |zw|_{Z_i} \quad \text{ by (\ref{eq:8.2})}\notag\\
&\leq \min \{ |zw|_{Y_i} + \epsilon_i, \left( |zA_i|_{Y_i} + \epsilon_i +\tau(\frac{1}{i}) \right) + \left( |A_{i}w|_{Y_i} + \epsilon_i + \tau(\frac{1}{i}) \right) \} - |zw|_{Z_i}\notag\\
& \quad \text{ by (\ref{ineq:8.4.subsii}), since $f$ an $\epsilon_i$-Hausdorff approximation and $A_i \longrightarrow A$ as subsets} \notag\\
&\leq \min \{ |zw|_{Y_i} + \epsilon_i, \left( |zA_i|_{Y_i} + \epsilon_i + \tau(\frac{1}{i}) \right) + \left( |A_{i}w|_{Y_i} + \epsilon_i + \tau(\frac{1}{i}) \right) \} \notag\\
& \quad - \min \{ |zw|_{Y_i} \quad , \quad |zA_i|_{Y_i}+|A_{i}w|_{Y_i} \} \notag\\
\intertext{
\quad \quad If $|zw|_{Y_i} \leq |zA_i|_{Y_i} + |A_{i}w|_{Y_i}$ then the right-hand side of this equals
}
r.h.s. &= \left( |zw|_{Y_i} + \epsilon_i \right) - |zw|_{Y_i} \notag\\
&= \epsilon_i \notag
\intertext{ 
\quad \quad Otherwise, if
$|zw|_{Y_i} \geq |zA_i|_{Y_i} + |A_{i}w|_{Y_i}$ then
}
r.h.s. &= \epsilon_{i} + \min \{ |zw|_{Y_i}, |zA_i|_{Y_i} + |A_{i}w|_{Y_i} + \epsilon_i + 2\tau(\frac{1}{i}) \} - \left( |zA_i|_{Y_i} + |A_{i}w|_{Y_i} \right) \notag\\
&= \epsilon_{i} + \min \{ \underset{\geq 0}{\underbrace{|zw|_{Y_i} - \left( |zA_i|_{Y_i} + |A_{i}w|_{Y_i} \right)}} , \, \epsilon_i +2\tau(\frac{1}{i}) \} \notag\\
&\leq \epsilon_i + \epsilon_i +2\tau(\frac{1}{i}) \notag
\intertext{ \quad \quad so that, either way, }
|F(z)F(w)|_{Z} - |zw|_{Z_i}
&\leq 2\epsilon_{i} + 2\tau(\frac{1}{i}) \notag\\
& \quad \longrightarrow 0
\end{align}

Conversely,
\begin{align}
|zw|_{Z_i} - |F(z)F(w)|_{Z}
&= |zw|_{Z_i} -  \min \{ |f(z)f(w)|_{Y}, |f(z)A|_{Y}+|Af(w)|_{Y} \} \quad \text{ by (\ref{eq:8.2}), since $z,w\in Y_i$ } \notag\\
&= \max \{  |zw|_{Z_i} - |f(z)f(w)|_{Y}, |zw|_{Z_i} - (|f(z)A|_{Y} + |Af(w)|_{Y})  \} \notag\\
&\leq \max \{ |zw|_{Y_i} - |f(z)f(w)|_{Y}, |zA_i|_{Y_i} + d(X_i) + |A_{i}w|_{Y_i} - \left( |f(z)A|_{Y} + |Af(w)|_{Y}  \right)  \} \notag\\
& \quad \text{ by triangle inequality, since ${(A_i)}_j$ path-connected }  \notag\\
& \quad \text{ (note that $|zA_i|_{Y_i} = |z {(A_i)}_j |_{Y_i}$ and likewise 
$|A_{i}w|_{Y_i} = |{(A_i)}_{j} w|_{Y_i}$ ) } \notag\\
&= \max \{ \epsilon_{i}, d(X_i) + \left( |zA_i|_{Y_i} -|f(z)A|_{Y} \right)   + \left(  |A_{i}w|_{Y_i} - |Af(w)|_{Y}  \right)  \} \notag\\
&\leq \max \{ \epsilon_{i}, d(X_i) + \epsilon_{i} + \epsilon_{i} +2\tau(\frac{1}{i}) \} \quad \text{ by (\ref{ineq:8.4.subsii}), since $f$ is an $\epsilon_i$-Hausdorff approximation } \notag\\
& \quad \text{ and $A_i \longrightarrow A$ as subsets } \notag\\
& \quad \longrightarrow 0
\end{align}

\underline{Case $z\in X_i \setminus A_i, w\in Y_i$} :

First note that $A_i \cap [wz]_{Z_i} \neq \emptyset$.
Let $v\in A_i \cap [wz]_{Z_i}$ be the last point on the $Z_i-$segment from $w$ to $z$ such that $[vz]_{Z_i} \subseteq X_i$. Then
\begin{align}\label{ineq:8.4.3}
|A_{i}w|_{Y_i} &= |v'w|_{Y_i} \quad \text{ for some $v'\in A_i$, since $A_i \subseteq Y_i$ closed} \notag\\
&= |v'w|_{Z_i} \quad \text{ since $[v'w]_{Z_i} \subseteq Y_i$} \notag\\
&\leq |zv|_{Z_i} + |v'w|_{Z_i} \notag\\
&\leq |zv|_{Z_i} + |vw|_{Z_i} \quad \text{ by choice of $v'$} \notag\\
&= |zw|_{Z_i} \quad \text{ by choice of $v$}
\end{align}

Thus
\begin{align}
|F(z)F(w)|_{Z} - |zw|_{Z_i}
&= |\pi(A) \pi f(w)|_{Y/A} - |zw|_{Z_i} \quad \text{ by definition of $F$} \notag\\
&\leq |\pi(A)\pi f(w)|_{Y/A} - |A_{i}w|_{Y_i} \quad \text{ by (\ref{ineq:8.4.3})} \notag\\
&\leq |Af(w)|_{Y} - |A_{i}w|_{Y_i} \notag\\
&\leq \epsilon_i + \tau(\frac{1}{i}) 
\quad \text{ by (\ref{ineq:8.4.subsii}) }\notag\\
& \quad \longrightarrow 0
\end{align}

Conversely,
let $v'\in A_i$ be chosen as before such that $|A_{i}w|_{Y_i} = |v'w|_{Y_i}$.

Then
\begin{align}\label{ineq:8.4.4}
|zw|_{Z_i} &\leq |zv'|_{Z_i} + |v'w|_{Z_i} \quad \text{ by the triangle inequality} \notag\\
&\leq |zv'|_{X_i} + |v'w|_{Y_i} \quad \text{ since $X_i \subseteq Z_i$ and $Y_i \subseteq Z_i$} \notag\\
&\leq d(X_i) + |v'w|_{Y_i} \notag\\
&= d(X_i) + |A_{i}w|_{Y_i}
\end{align}

Hence
\begin{align}\label{ineq:8.6}
|zw|_{Z_i} - |F(z)F(w)|_{Z}
&= |zw|_{Z_i} - |\pi(A)\pi f(w)|_{Z} \quad \text{ by definition of $F$} \notag\\
&= |zw|_{Z_i} - |\pi(A)\pi f(w)|_{Y/A} \notag\\
&= |zw|_{Z_i} - |Af(w)|_{Y} \quad \text{ since $f(w)\in Y$ and ${\pi}^{-1}(\pi(A))=A$} \notag\\
&\leq d(X_i) + |A_{i}w|_{Y_i} - |Af(w)|_{Y} \quad \text{ by (\ref{ineq:8.4.4})} \notag\\
&\leq d(X_i) + \epsilon_i + \tau(\frac{1}{i}) \quad \text{ by (\ref{ineq:8.4.subsii})} \notag\\
& \quad \longrightarrow 0
\end{align}

(\ref{ineq:8.5})-(\ref{ineq:8.6}) together prove that $F : Z_i \longrightarrow Z$ is a $(d({X_i})+2\epsilon_i+2\tau(\frac{1}{i}))$-Hausdorff approximation. 
Therefore $Z_{i} \overset{GH}{\longrightarrow} Z=Y/A$.\footnote{observe that all the inequalities of the form $|zw|_{Z_i} - |F(z)F(w)|_{Z} \leq \ldots$ require the upper diameter bound $d(X_i)$.}
\end{proof}

\vspace{1cm}

\begin{proof}[Proof of Theorem~\ref{thm:Top}]

Suppose to the contrary that $M_i$ has at least three boundary components, for all $i$.

For each $i$ there exists a point $p_{i} \in M_i$ such that $inrad(p_i) = inrad(M_i)$ and for which the closure of $B(p_i, inrad(p_i); M_i)$ intersects $\partial M_i$ in at least $3$ distinct points, call them ${\{ 0 \}}_{i,k}$ $(k=1,\ldots 3)$, where ${\{ 0 \}}_{i,k}$ lies in the $k$-th component of $\partial M_i$. $p_i$ can be chosen as a cut point of the boundary, of order $3$.

Form the Alexandrov extension $\widetilde{M_i}$ of $M_i$.
There exist $3$ distinct, disjoint segments $\{ I_{i,k}:=[0,t_{0}]_{i,k} \}_{k=1}^{3}$ such that $I_{i,k} \subset C_{M_i,k}:={(\partial M_i)}_{k} \times_{\phi} [0,t_{0}]$.

By Lemma \ref{lem:4.4}, there exist $3$ limit segments ${ \{ I_{k} \} }_{k=1}^{3}$ in $\underset{GH}{\lim}\widetilde{M_i}$.

For any $t \geq t_{0}/2$, if ${\{ t \}}_{i,j}$ and ${\{ t \}}_{i,k}$ lie in ${ [0,t_{0}] }_{i,j}$ and  ${ [0,t_{0}] }_{i,k}$ respectively, then
$d_{\widetilde{M_i}}( {\{ t \}}_{i,j}, {\{ t \}}_{i,k}) \geq 2t$ for all $i$ and all $j\neq k$
since ${ \{ I_{i,k} \} }_{k=1}^{3}$ are pairwise distinct and lie in different components $C_{M_i,k}$. On the other hand,
\begin{align*}
d_{\widetilde{M_i}}( {\{ 0 \} }_{i,j},{ \{ 0 \}  }_{i,k}) \leq d_{M_i}({ \{ 0 \} }_{i,j},{ \{ 0\} }_{i,k}) \leq 2 inrad(M_i) \longrightarrow 0
\end{align*}
for any $j,k=1,\ldots,3$.

Thus all $I_{k}$ intersect precisely at their left endpoints. Their union would yield a branch point in the limit of the $\widetilde{M_i}$'s. But
$curv \left( \underset{GH}{\lim}\widetilde{M_i} \right) \geq k$.
Therefore $\partial M_i$ can have no more than two components.\\

Now we prove the second statement of Theorem~\ref{thm:Top}.
\footnote{in the course of the proof, the conclusion of the first statement will be reproved in a different manner}

Write $\partial M_i = \coprod_{k=1}^{m}{(\partial M_i)}_{k}$ as a union of components. 
So ${\partial M_i} {\times}_{\phi} [0,t_{0}] = \coprod_{k=1}^{m}\left(  {(\partial M_i)}_{k}{\times}_{\phi} [0,t_{0}] \right)$ and $\widetilde{M_i}=M_i \cup_{\partial M_i} \left( \partial M_i {\times}_{\phi} [0,t_{0}] \right)$ as in Proposition~\ref{prop:2.1}.

The first goal is to show that
$M_i \overset{GH}{\longrightarrow} pt$ implies $(\partial M_i)_{k} \overset{GH}{\longrightarrow} pt$ for every $k$.
It suffices to show that $d(\partial M_i) \longrightarrow 0$.

This will be carried out using
Proposition~\ref{lem:2.8}
and a scaling trick, as follows:

Let $d_{i}:=d(M_i) \longrightarrow 0$. 

Consider the rescaled manifolds
$\overline{M_i}:=\frac{1}{d_i} M_i$

Then
\begin{align*}
d(\overline{M_i}) &\equiv 1,\\
K_{\overline{M_i}} &\geq -\tau(\frac{1}{i}),\\
\text{ and } |II_{\partial \overline{M_i}}| &\leq \tau(\frac{1}{i}).
\end{align*}

By
Proposition~\ref{lem:2.8},
each component of $\partial \overline{M_i}$ therefore has intrinsic diameter bounded above by
$d(\partial \overline{M_i}) \leq D$
for some constant $D$ independent of $i$.

Scaling down,
\begin{align}\label{eq:4.1}
d(\partial M_i) = d_{i}d(\partial \overline{M_i}) \leq d_{i} D \longrightarrow 0
\end{align}
as desired.

Now apply Proposition~\ref{prop:8.4} with
\begin{align*}
X_i &= M_i\\
Y_i &= \partial M_i {\times}_{\phi} [0,t_{0}]\\
A_i &= \partial M_i\\
Z_i &= \widetilde{M_i}
\end{align*}

By (\ref{eq:4.1}),  $ { (\partial M_i) }_{k} \overset{GH}{\longrightarrow} {pt}_{k}$ for each boundary component $ { (\partial M_i) }_{k} $, so ${(A_i)}_{k} =  {( \partial M_i)}_{k} \longrightarrow {pt}_{k} = {(A)}_{k}$ as subsets, under the convergence $Y_i \overset{GH}{\longrightarrow} Y$. 
In particular, each ${(A)}_{k}\equiv { pt  }_{k}$ is certainly path-connected, and Proposition~\ref{prop:8.4} is applicable as stated.

Note that
\begin{align*}
A = \underset{k}{\coprod}{\{ pt \} }_{k}
\end{align*}

Then
\begin{align*}
Y &= \underset{GH}{\lim} \left( \partial M_i {\times}_{\phi} [0,t_{0}] \right) \\
&= \left( \underset{GH}{\lim} \, \partial M_i \right) {\times}_{\phi} [0,t_{0}] \quad \text{by Proposition~\ref{lem:4.7}}\\ 
&=  \underset{k}{\coprod} {\{ pt \}}_{k} {\times}_{\phi} [0,t_{0}] \quad \text{ by (\ref{eq:4.1})}
\end{align*}

Therefore by Proposition~\ref{prop:8.4}, $Z=Y/A$ isometrically, or
\begin{align*}
\underset{GH}{\lim} \widetilde{M_i}
&= \left( \underset{k}{\coprod} {\{ pt \}}_{k} {\times}_{\phi} [0,t_{0}] \right)
 / \left( \underset{k}{\coprod}{\{ pt \} }_{k} \right) \\
& = \left( \underset{k}{\coprod} {\{ pt \}}_{k} {\times}_{\phi} [0,t_{0}] \right) / \left( pt_{j} \sim pt_{k} \right)\\
& = \text{ wedge of $m$ intervals $[0,t_{0}]$ }
\end{align*}

Observe that a boquet of $m$ intervals does not admit $curv\geq k$ if $m \geq 3$.

However, $\underset{GH}{\lim} \widetilde{M_i}$ has $curv\geq k$.
Therefore $m\leq 2$.\\

Suppose $m=2$. Then $Z:=\underset{GH}{\lim} \widetilde{M_i} = D^1$
is the wedge of two intervals $\{pt_1 \} {\times}_{\phi} [0,t_{0}]$ and $\{pt_2 \} {\times}_{\phi} [0,t_{0}]$.

Both $\widetilde{M_i}$ and $Z$ admit natural codimension-$0$ isometric extensions $\widetilde{\widetilde{M_i}}$ and $\widetilde{Z}$, respectively, since their boundaries are totally geodesic. Clearly one still has $\widetilde{\widetilde{M_i}} \overset{GH}{\longrightarrow}\widetilde{Z}$.  $\widetilde{M_i}$ and $Z$ are contained in open subsets $U$ and $V$ of $\widetilde{\widetilde{M_i}}$ and $\widetilde{Z}$, respectively. By Theorem~\ref{thm:4.4}, there is a $\tau(\frac{1}{i})$-almost Lipschitz submersion $f: U \longrightarrow V$, which is simultaneously a $C^1$-locally trivial fiber bundle, since $\widetilde{M_i}$ and $Z$, hence $U$ and $V$, admit $C^1$ differentiable structures. Since the base is contractible, $M_i \underset{homeo}{\approx} \widetilde{M_i}$ is actually a trivial fiber bundle over $D^1$.
\end{proof}

\vspace{2cm}

\subsection{Proof of Theorem~\ref{thm:4.3}}

\begin{proof}[Proof of Theorem~\ref{thm:4.3}]

Consider, as in Proposition~\ref{prop:2.1}, (Alexandrov) extensions of $M_i$, but where the warping function $\phi$ is not fixed for all $i$, but rather varies with $i$ (i.e. is ``optimally adjusted'' for $M_i$). In other words, consider 
\begin{align*}
\widetilde{M_{i,{\phi_{i}}}} = M_{i} \cup \left( \partial M_{i} {\times}_{\phi_{i}} [0,t_{0,i}] \right)
\end{align*}
where the functions $\phi_{i}$ satisfy
\begin{align*}
{\phi_{i}}'(0) &= {\overline{\lambda}}_{i}:=\min \{ 0, {\lambda}^{-}(M_i) \} = -\tau(\frac{1}{i})\\
{\phi_{i}}(t_{0,i}) &= {\epsilon}_{i}
\end{align*}
in addition to the other conditions in Lemma~\ref{lem:2.1}, where $t_{0,i}$ and $\epsilon_{i}$ are chosen as
\begin{align}
t_{0,i} &:= \frac{10}{\sqrt{i}}\label{const:4.1}\\
\epsilon_{i} &:= 1-\frac{1}{i^{3/2}}.\label{const:4.2}
\end{align}

(This essentially means that the $\phi_{i}$ are approaching the function which is identically equal to 1 at $t=0$.)

Then just as in Proposition~\ref{prop:2.1}
\begin{align*}
curv \widetilde{M_{i,\phi_{i}}} \geq k_{i}:=\min \{ K^{-}, K_{C_{i}}^{-} \}
\end{align*}
where the right-hand side may depend on $i$, but remains bounded below (in terms of $K^{-}, {\lambda}^{\pm}$) as $i \longrightarrow \infty$.

To see this, begin by noting that
given ${\overline{\lambda}}_{i} = -\tau(\frac{1}{i})$, it may be assumed by taking a subsequence that 
\begin{align}
{\overline{\lambda}}_{i} = -\tau(\frac{1}{i}) = -\frac{1}{i}.
\end{align}

Then $\frac{ | {\overline{\lambda}}_{i}| t_{0,i}}{1-{\epsilon}_{i}} = 10$
so
by (\cite[(2.8)]{Wo:07}),
the radial curvatures
are bounded below by
$\underset{t}{\inf} \left( -\frac{{\phi}_{i}''(t)}{{\phi}_{i}(t)} \right)$,
which is no less than 
the least of the quantities 0 and
\begin{align}\label{functineq:7}
-\frac{1}{ {\epsilon}_{i}}
\left[
2\frac{{\overline{\lambda}}_{i}}{t_{0,i}}
+ \frac{{( {\overline{\lambda}}_{i} )}^{2}}{1-{\epsilon}_{i} }
\right]
= -\frac{1}{1-\frac{1}{i^{3/2}}} \left[
-\frac{2}{i}\left( \frac{i^{1/2}}{10} \right) + \frac{1/i^2}{1/i^{3/2}}
\right]
= -\frac{1}{1-\frac{1}{i^{3/2}}} \left[
-\frac{1}{5 i^{1/2}} + \frac{1}{i^{1/2}}
\right]
= -\frac{4}{5}\frac{1}{1-\frac{1}{i^{3/2}}}\cdot \frac{1}{i^{1/2}} \, ,
\end{align}
i.e., the radial curvatures are bounded below by the last quantity.

Again, since $\frac{|{\overline{\lambda}}_{i}| t_{0,i}}{1-{\epsilon}_{i}} = 10$,
(\cite[(2.11)]{Wo:07}) implies that
the tangential curvatures
are bounded below by
\begin{align}\label{functineq:8}
\frac{1}{{\phi}_{i}^2(t)}
 \left[ K_{\partial M_i}^{-} - {|{\phi}_{i}'(t)|}^2 \right]
&\geq \frac{1}{{\epsilon}_{i}^2}\left[ \min \{ K_{\partial M_i}^{-},0 \} -
{| {\overline{\lambda}}_{i} |}^2
\right] \notag\\
&= \frac{1}{{(1-1/i^{3/2})}^2}
\left[ \min \{ K_{\partial M_i}^{-},0 \} -
\frac{1}{i^2}
\right]
\end{align}
where $K_{\partial M_i}^{-} \geq K_{M_i}^{-} -
\tau(\frac{1}{i}){\lambda}^{+}
= K^{-} - \tau(\frac{1}{i})$

It follows that the right-hand side of (\ref{functineq:8}) tends to $K^{-}$ as $i\longrightarrow \infty$.\\

If in an orthonormal frame the sectional curvatures are bounded (from below) on all coordinate two-planes, then the sectional curvatures are bounded (from below) on arbitrary two-planes. Therefore (\ref{functineq:7}) and (\ref{functineq:8}) prove that $K_{C_i}$, hence $k_{i}$, is uniformly bounded from below, independent of $i$.\\

Let $x_i \in M_i$ be arbitrary.
Suppose there exists a (pointed) limit
\[
(X,x) = \underset{p-GH}{\lim} (M_i,x_i)
\]
where $X$ is a length space. It may be assumed that $X$ is metrically (Cauchy) complete.

Let $R>0$ be arbitrary.

$(M_i,x_i) \overset{p-GH}{\longrightarrow} (X,x)$ (and $X$ complete) implies that
\begin{align}
B(x_i, R; M_i) \overset{GH}{\longrightarrow} B(x,R;X),
\end{align}
where the balls are equipped with the respective restricted metrics.
(Here, for simplicity of notation, in contrast to the usual conventions made, $B(x,R;X)$ etc. denotes a closed ball of radius $R$ in $X$ with center $x$.)

For any $x_i \in M_i$, $ B(x_i,R; M_i) \subseteq B(x_i,R; \widetilde{M_{i,{\phi_{i}}}})$

Consider the maps $\pi_{i} : \widetilde{M_{i,{\phi_{i}}}} \longrightarrow M_i$ from Lemma~\ref{lem:2.4}.
The restriction
\[
\pi_{i}| : B(x_i,R; \widetilde{M_{i,{\phi_{i}}}}  ) \longrightarrow B(x_i,R; M_i)
\]
is surjective and is, by Lemma~\ref{lem:2.5}, a $\tau\left(\frac{1}{i}\right)$-Hausdorff approximation, where
$\tau\left(\frac{1}{i}\right)=\max \{ 2t_{0,i}, \left( \dfrac{1}{\epsilon_i}-1 \right) (R+t_{0,i})  \}
= \max \{ \frac{20}{\sqrt{i}}, \frac{1}{i^{3/2}-1}(R+\frac{20}{\sqrt{i}}) \}$
(Here again the balls are equipped with the restricted metrics of their ambient spaces.)

By precompactness
of the class of Alexandrov spaces with lower curvature bound, with respected to the pointed Gromov-Hausdorff topology,
there exists a (pointed) limit
\[
(Y,y) = \underset{p-GH}{\lim} (\widetilde{M_{i,{\phi_{i}}}}, x_i)
\]
for some $y\in Y$ with $x_i \mapsto y$, say. This entails that
\begin{align}
B(x_i,R; \widetilde{M_{i,{\phi_{i}}}}) \overset{GH}{\longrightarrow} B(y ,R;Y)
\end{align}
and so by the triangle inequality
\begin{align*}
d_{GH}(B(y, R; Y), B(x,R; X))
&\leq d_{GH}(B(y,R; Y), B(x_i, R; \widetilde{M_{i,{\phi_{i}}}}))\\
&+ d_{GH}(B(x_i, R; \widetilde{M_{i,{\phi_{i}}}}), B(x_i,R; M_i))\\
&+ d_{GH}(B(x_i,R; M_i), B(x, R; X))
\end{align*}

Letting $i\longrightarrow \infty$, one has $d_{GH}(B(y,R;Y),B(x,R;X) )=0$.
Since the balls are closed, $B(y,R;Y) \underset{isom}{\equiv} B(x,R;X)$. Since $Y$ has lower curvature bound, so does $B(y,R;Y)$ and hence so does $B(x,R;X)$. But $x$ was arbitrary, so by Toponogov's theorem, $curv X$ is bounded below.
\end{proof}

\vspace{1cm}

\vspace{1cm}

\subsection{Proof of Theorem~\ref{thm:4.6}}

The first part (1) will be proved independently from the other parts (2)-(3), and (2)-(3) also do not rely on (1).
The proof of parts (2) and (3) uses work of Alexander, Berg, Bishop, Lytchak, and Nikolaev.\\

\begin{proof}[Proof of Theorem~\ref{thm:4.6} (1)]

In order to show that $inrad(M_i) \longrightarrow 0$, it suffices, by Lemma~\ref{lem:4.1}, to show that $i_{\partial, M_i}(p_i) \longrightarrow 0$ for any choice of points $p_i \in \partial M_i$, $(i=1,2,\ldots)$.
Suppose by way of contradiction that for each $i$, a point $p_i \in \partial M_i$ is chosen arbitrarily,
and $i_{\partial, M_i} (p_i) \geq R$, for some fixed $R>0$ which is independent of $i$.

By definition of $i_{\partial, M_i}(p_i)$, there exists points
$q_i \in M_i$, with $[q_{i}p_{i}]_{M_i}$ a minimizing segment hitting $\partial M_i$ orthogonally at $p_i$, and with $|q_{i}p_{i}|_{M_i} = R > 0$ for all $i$.

Suppose $q_i \mapsto q$ and $p_i \mapsto p$ for some $q, p\in X$. There exists a minimal geodesic segment $[qp]_{X}$ since $X$ is a geodesic metric space. Since $X$ is assumed geodesically extendible, one can extend $[qp]_{X}$ to $[qpr]_{X}$ where $|pr|_{X}=\delta$ for some $\delta > 0$. There exists a sequence of points $r_i \in M_i$
with $r_i \mapsto r$. \\
Then
\begin{align*}
|q_{i}p_{i}|_{M_i} &= R\\
|p_{i}r_{i}|_{M_i} &\leq |pr|_{X} + \tau(\frac{1}{i}) = \delta + \tau(\frac{1}{i})\\
|q_{i}r_{i}|_{M_i} &\leq R + \delta + \tau(\frac{1}{i}).
\end{align*}
where the last inequality is the triangle inequality.

Since int$(M_i)$ is open in $M^n$, $K_{M_i} \equiv K_{M}$ and so $|K_{M_i}| \leq K$. Then $curv M_i \leq k = k(K,{\lambda}^{-})$ (\cite{ABB:93}).

Although $M$ was not necessarily assumed compact, the argument in what follows will be local, so consider a closed ball $\overline{B}(p,100; M)$ in $M$, of some fixed radius, e.g. $100$.
Since the ball is compact, one has $inj_{M}(x) \geq i_0 > 0$ for all $x\in \overline{B}(p,100; M)$, for some constant $i_{0}$ depending on $\overline{B}(p,100; M)$.
Hereafter assume that all distances and lengths are less than $100$.
Effectively then, one may as well assume
$inj(M) \geq i_{0} > 0$.

By \cite{AB:03}, the lower bound $inj(M) \geq i_{0} > 0$ implies a lower bound $inj(M_{i}) \geq i_{1} > 0$ for some constant $i_{1}=i_{1}(K,{\lambda}^{-},i_{0})$.

The conditions $curv(M_i) \leq k$ and $inj(M_i)\geq i_{1} > 0$ imply that each $M_i$ is $CAT_{k}$ for a uniform lower $CAT_{k}Rad$ bound independent of $i$. Since $M_i \overset{GH}{\longrightarrow} X$, this is known to imply that $X$ is also $CAT_{k}$. In particular, $curv(X)\leq k$ and $inj(X)\geq i_{1} > 0$.

Given that $inj(X)\geq i_{1} > 0$ and $e_{qp}(r):=|qp|_{X}+|pr|_{X}-|qr|_{X} = 0$, if $[qr]_{X}$ and $[qp]_{X}\cup [pr]_{X}$ were distinct this would force $|qr|\geq i_{1}$.
However,
$|qr|_{X}\leq |qp|_{X}+|pr|_{X} = R + \delta < i_{1}$ if $R,\delta > 0$ are chosen sufficiently small relative to $i_{1}$ (e.g., $R,\delta < \frac{i_{1}}{2}$). Therefore $[qr]_{X}$ and $[qp]_{X}\cup [pr]_{X}$ coincide.

Furthermore,
\begin{align}\label{ineq:c.1}
d_{\infty}(\sigma_{i},\gamma_{i}) := 
\underset{0\leq t \leq L_i}{\sup} d(\sigma_i(t),\gamma_i(t)) \longrightarrow 0
\end{align}

where 
\begin{align*}
\sigma_i &= [q_{i}p_{i}]_{M_i} \cup [p_{i}r_{i}]_{M_i}\\
\gamma_i &= [q_{i}r_{i}]_{M_i}
\end{align*}
are parametrized on a common interval $[0,L_i]$. 
For otherwise, one would get nontrivial minimal geodesics $\sigma=\lim \sigma_i$, $\gamma = \lim \sigma_i$ in $X$ by Lemma~\ref{lem:4.4}, and $d_{\infty}(\sigma_i, \gamma_i)\geq r > 0$ would imply that $d_{\infty}(\sigma, \gamma)\geq r > 0$, i.e., $\sigma, \gamma $ would be distinct.

If $R,\delta < R_{0}:= \min \{ \frac{\pi}{2 \sqrt{K}}, i_{0}/4 \} $, then
\begin{align*}
{\Delta}^{M_i}q_{i}p_{i}r_{i} \subset B(p_{i},R_{0} ; M)
\end{align*}
for all sufficiently large $i$.
The latter is a convex ball in
$M$ since the convexity radius has the lower bound $conv (M)\geq \min\{ \frac{\pi}{2\sqrt{K}}, \frac{1}{2}inj(M) \} > 0$.  In particular, this ball is a $CAT(K)$-domain.

All of the geodesic sides $[p_{i}r_{i}]_{M_i}$, $[q_{i}p_{i}]_{M_i}$, and
$[q_{i}r_{i}]_{M_i}$ of the triangle ${\Delta}^{M_i}q_{i}p_{i}r_{i}$ have pointwise arc/chord curvature $\kappa_{g} \leq |{\lambda}^{-}| $ when considered as curves in $M$.

By arc/chord comparison in $B(p_{i}, R_{0} ; M)$ (see $\S$\ref{s:A2}), there are base angle estimates
\begin{align*}
{\angle}^{M}([p_{i}r_{i}]_{M_i},[p_{i}r_{i}]_{M})
&\leq \frac{|{\lambda}^{-}|}{2}|p_{i}r_{i}|_{M} + O(|p_{i}r_{i}|_{M}^{3})\\
&= \frac{|{\lambda}^{-}|}{2}|p_{i}r_{i}|_{M_i} + O(|p_{i}r_{i}|_{M_i}^{3})\\
&\leq \frac{|{\lambda}^{-}|}{2}(\delta + \tau(\frac{1}{i}) ) + O\left( {(\delta + \tau(\frac{1}{i}) )  }^{3} \right)
\end{align*}
and
\begin{align*}
{\angle}^{M}([q_{i}p_{i}]_{M_i},[q_{i}p_{i}]_{M})
&\leq \frac{|{\lambda}^{-}|}{2}|q_{i}p_{i}|_{M} + O(|q_{i}p_{i}|_{M}^{3})\\
&= \frac{|{\lambda}^{-}|}{2}|q_{i}p_{i}|_{M_i} + O(|q_{i}p_{i}|_{M_i}^{3})\\
&\leq \frac{|{\lambda}^{-}|}{2}(R + \tau(\frac{1}{i}) ) + O\left( {(R + \tau(\frac{1}{i}) )  }^{3} \right)
\end{align*}
(Here the $O's$ depend implicitly on the fixed values $K$ and ${\lambda}^{-}$.)

By the triangle inequality,
\begin{align}
{\angle}^{M}([q_{i}p_{i}]_{M},[p_{i}r_{i}]_{M})
&\leq {\angle}^{M}([q_{i}p_{i}]_{M_i},[p_{i}r_{i}]_{M_i})
+ {\angle}^{M}([p_{i}r_{i}]_{M_i},[p_{i}r_{i}]_{M})
+{\angle}^{M}([q_{i}p_{i}]_{M_i},[q_{i}p_{i}]_{M})\\
&\leq \frac{\pi}{2} 
+ \frac{|{\lambda}^{-}|}{2}(\delta + \tau(\frac{1}{i}) ) + O\left( {(\delta + \tau(\frac{1}{i}) )  }^{3} \right)
+ \frac{|{\lambda}^{-}|}{2}(R + \tau(\frac{1}{i}) ) + O\left( {(R + \tau(\frac{1}{i}) )  }^{3} \right)\\
&= \frac{\pi}{2}+O\left( R+\delta+2\tau(\frac{1}{i}) \right)\label{ineq:c.3}.
\end{align}
In other words, the angle formed by the segments $ [q_{i}p_{i}]_{M}$ and $[p_{i}r_{i}]_{M}$ is close to a right angle, if $R, \delta$ are sufficiently small and $i$ sufficiently large.

Furthermore,
if $\Delta^{{M_i}}q_{i}p_{i}r_{i}$ has sufficiently small perimeter, then
\begin{align}
d_{H}([q_{i}r_{i}]_{M_i}, [q_{i}r_{i}]_{M})
&\leq \frac{|{\lambda}^{-}|}{8}|q_{i}r_{i}|_{M}^{2} + O(|q_{i}r_{i}|_{M}^{3})\\
&= \frac{|{\lambda}^{-}|}{8}|q_{i}r_{i}|_{M_i}^{2} + O(|q_{i}r_{i}|_{M_i}^{3})\\
&\leq \frac{|{\lambda}^{-}|}{8} {\left( R+\delta+\tau(\frac{1}{i}) \right) }^{2} + O\left(  { \left( R+\delta+\tau(\frac{1}{i}) \right) }^{3} \right)\\
&= O\left( {\left( R+\delta+\tau(\frac{1}{i})\right)}^{2} \right) \label{ineq:c.4}
\end{align}
by width estimates (cf. $\S$\ref{s:A2}).\\

Now the goal is to show that 
the uniform closeness  (\ref{ineq:c.1}) of $\gamma_{i}$ to $\sigma_{i}$
implies, for sufficiently small $R,\delta$ and large $i$, that $d(p_{i},[q_{i}r_{i}]_{M})$ would be smaller than is allowed by a Toponogov excess estimate.
On the one hand,
\begin{align}
d(p_{i},[q_{i}r_{i}]_{M})
&\leq d(p_{i},[q_{i}r_{i}]_{M_i}) + d_{H}([q_{i}r_{i}]_{M_i},[q_{i}r_{i}]_{M})\\
&\leq d_{\infty}(\sigma_{i},\gamma_{i})
+ d_{H}([q_{i}r_{i}]_{M_i},[q_{i}r_{i}]_{M})\\
&\leq \tau(\frac{1}{i}) + d_{H}([q_{i}r_{i}]_{M_i},[q_{i}r_{i}]_{M})
\quad \text{ by (\ref{ineq:c.1})}\\
&\leq \tau(\frac{1}{i}) + O\left( {\left( R+\delta+\tau(\frac{1}{i})\right)}^{2} \right)
\quad \text{ by (\ref{ineq:c.4}).}
\label{ineq:c.5}
\end{align}

However, by the triangle inequality,
\begin{align}
d(p_{i},[q_{i}r_{i}]_{M})
&\geq \frac{1}{2}\left( |q_{i}p_{i}|_{M}+|p_{i}r_{i}|_{M}-|q_{i}r_{i}|_{M} \right)\\
&\geq \frac{1}{2}\left( R-\tau(\frac{1}{i})+\delta-\tau(\frac{1}{i})-|q_{i}r_{i}|_{M} \right). \label{ineq:c.6}
\end{align}

Since $K_{M} \geq -K$ and $B(p_{i}, R_{0} ; M)$ is a convex ball in
$M$, one may apply Toponogov's comparison theorem (hinge version) inside it.
\begin{align*}
|q_{i}r_{i}|_{M} &\leq \frac{1}{\sqrt{K}}{\cosh}^{-1} \left[
\cosh(\sqrt{K} |q_{i}p_{i}|_{M}) \cosh(\sqrt{K} |p_{i}r_{i}|_{M})
\right.\\
& \hspace{2.2cm} - \left. 
\sinh(\sqrt{K} |q_{i}p_{i}|_{M}) \sinh(\sqrt{K} |p_{i}r_{i}|_{M})
\cos({\angle}^{M}q_{i}p_{i}r_{i})
\right]\\
&\leq \frac{1}{\sqrt{K}}{\cosh}^{-1} \left[
\cosh(\sqrt{K} (R+\tau(\frac{1}{i}))) \cosh(\sqrt{K} (\delta+\tau(\frac{1}{i}))) \right.\\
& \hspace{2.2cm} + \left. \sinh(\sqrt{K} (R+\tau(\frac{1}{i}))) \sinh(\sqrt{K} (\delta+\tau(\frac{1}{i}))) 
\sin \left(
O\left( R+\delta+2\tau(\frac{1}{i}) \right)
\right)
\right]
\end{align*}
by (\ref{ineq:c.3}).

Substituting this into (\ref{ineq:c.6}), one obtains a lower bound
\begin{align}
d(p_{i},[q_{i}r_{i}]_{M})
&\geq \frac{1}{2}\left[ R-\tau(\frac{1}{i})+\delta-\tau(\frac{1}{i})
- 
\frac{1}{\sqrt{K}}{\cosh}^{-1} \left[
\cosh(\sqrt{K} (R+\tau(\frac{1}{i}))) \cosh(\sqrt{K} (\delta+\tau(\frac{1}{i}))) \right. \right.\notag \\
& \hspace{2.2cm} + \left. \left. \sinh(\sqrt{K} (R+\tau(\frac{1}{i}))) \sinh(\sqrt{K} (\delta+\tau(\frac{1}{i}))) 
\sin \left(
O\left( R+\delta+2\tau(\frac{1}{i}) \right)
\right)
\right]
\right] \notag \\
&\geq O\left( \min\{ R, \delta \} -\tau(\frac{1}{i}) \right)
\end{align}
which contradicts the upper bound (\ref{ineq:c.5}) if $R$ and $\delta$ are restricted to be sufficiently small (yet positive and independent of $i$), and $i$ is sufficiently large (e.g. so that $\tau(\frac{1}{i}) << R, \delta$).
One concludes from this contradiction that $R \longrightarrow 0$ as $i \longrightarrow \infty$.
This completes the proof of
Theorem~\ref{thm:4.6} (1)
that $inrad(M_i)\longrightarrow 0$.
\end{proof}

\footnote{\it Side Remark \rm :  
In the absence of an ambient space, an analogue of Theorem~\ref{thm:4.6} part (1) can still be proven, just assuming $M_{i}\in \mathcal{M}(n,K^{\pm},{\lambda}^{-},inj(M_i)\geq i_1)$ and the following 

Extension Hypothesis: For any $M \in \mathcal{M}(K^{\pm},{\lambda}^{-})$ there exists $\rho=\rho(K^{\pm},{\lambda}^{-})>0$ and a smooth extension $\widetilde{M}$ of $M$, of the same dimension, such that
\begin{align*}
K^{-}-C(\rho) &\leq K_{\widetilde{M}} \leq K^{+}+C(\rho),\\
\rho &\leq d(\partial M, \partial \widetilde{M})
\end{align*}
for some constant $C(\rho)<\infty$.\\

Under the stated assumptions of Theorem~\ref{thm:4.6}, this hypothesis follows since $M_{i}^{n}$ admits a uniform immersed tubular neighborhood in $M^n$: the outward focal distance of $\partial M_i$ is $\geq c(K^{+}, {\lambda}^{-})=\frac{1}{\sqrt{K^{+}}} \arctan \left( \frac{\sqrt{K^{+}}}{{\lambda}^{-}} \right) > 0$,
so there exists a uniform extension of $M_i$ (via the outward normal bundle of $\partial M_i$ in $M$) satisfying the conditions of the hypothesis.\\

By the extension hypothesis, each $M_i$ may be extended to a smooth manifold-with-boundary $\widetilde{M_i}$ with
 $ K^{-}-C(\rho) \leq K_{\widetilde{M_i}} \leq K^{+}+C(\rho)$. And $inj_{\widetilde{M_i}}(p) \geq i_{0}(i_{1},K^{\pm},\rho) > 0$ for some $i_{0}$ and all $p \in \widetilde{M_i}$ with $d(p,\partial {\widetilde{M_i}}) \geq \rho/4$, by standard injectivity radius decay estimates.
Then one can repeat the same proof as just given above, replacing the fixed $M$ there with $\widetilde{M_i}$. 
}

\vspace{0.4cm}

\vspace{1.4cm}

The remaining parts (2), (3)
of Theorem~\ref{thm:4.6}
are independent of (1), and may be derived from the following proposition.

\begin{proposition}\label{prop:4.1}(from earlier preprint version of \cite{Lyt:01}
\footnote{Note: From the older preprint version of \cite{Lyt:01}, Lytchak deleted (for editorial reasons rather than lack of correctness) the third section
which was called "Riemannian manifolds" and which essentially contained
proposition~\ref{prop:4.1} above.})
Let $M$ be a Riemannian manifold with bi-laterally
bounded Alexandrov curvature and with intrinsic injectivity radius bounded below by a positive
constant. If $N \subset M$ is a $C^{1,1}$ (immersed) submanifold without
boundary, then $N$ itself has (Alexandrov) curvature bi-laterally
bounded.
\end{proposition}

\begin{proof}(from earlier preprint version of \cite{Lyt:01}) 
Approximating the metric on $M$ by $C^{\infty}$-smooth Riemannian metrics, in the $C^{1,1}$ topology, it may be assumed that $M$ is $C^{\infty}$-smooth.
Approximating $N$ by $C^{\infty}$-smooth Riemannian submanifolds with uniformly bounded $C^{1,1}$ norms, it may be assumed that $N$ is $C^{\infty}$-smooth.

The $C^{1,1}$ norm is controlled, by Theorem 1.2 of \cite{Lyt:01}.

It suffices to show that the curvature and the injectivity radius of a smooth submanifold $N$ of a smooth Riemannian manifold $M$ can be bounded by the curvature of $M$ and the $C^{1,1}$ norm of $N$.

The curvature bound follows from the Gauss Equations.

The lower bound on injectivity radius follows from Lemma 6.1 of \cite{Lyt:01}.
\end{proof}

\vspace{0.5cm}

\begin{corollary}\label{cor:4.2}
Let $M$ be a Riemannian manifold with $K^{-}\leq K_{M} \leq K^{+}$, $Z\subset M$ a geodesically extendible $(C,2)$-convex subset.
Then $c^{-}(K^{-},C) \leq curv Z \leq c^{+}(K^{+},C)$.
\end{corollary}

\begin{proof}
By
\cite[Theorem 1.3]{Lyt:01}
if $M$ is a ($C^{0}$-) Riemannian manifold
and $Z \subset M$ is a compact $(C, 2)$-convex subset, then $Z$ has positive
reach. Since $Z$ is geodesically extendible by assumption, $Z$ is a
$C^{1,1}$ submanifold
\cite[proposition 1.4]{Lyt:01}.
Then Corollary~\ref{cor:4.2} follows from Proposition~\ref{prop:4.1}.
\end{proof}

\vspace{1cm}

\begin{proof} [Proof of Theorem~\ref{thm:4.6} (2)]
$K_{M_i} \equiv K_{M} \leq K^{+}$ since int$(M_i)$ is open in $M^n$.
By the same reasoning as in the proof of part (ii)(1), it may be assumed, by working locally on a fixed compact ball of large radius, that
$inj(M) \geq i_{0}$.
For all $i$, $inj(M_i) \geq i_{1}(K^{+},{\lambda}^{-},i_{0}) > 0 \ $ (\cite{AB:03}).

By \cite{ABB:93}, $curv M_i \leq k=k(K^{+},{\lambda}^{-})$.
The hypothesis on ${\lambda}^{-}$ implies that $M_i$ is $(C,2,i_{0})$-convex in $M^n$, where $C=C({\lambda}^{-}) = \frac{{({\lambda}^{-})}^2}{24}$.

Therefore, if $M_i \overset{GH}{\longrightarrow} X$ and $M_i \overset{H}{\longrightarrow} X$, then $X\subset M^n$ is $(C,2,i_{0}/2)$-convex, by Lemma~\ref{lem:C.3}. By Corollary~\ref{cor:4.2}, $curv X \geq c(K^{-},{\lambda}^{-})$ for some constant $c$.
\end{proof}

\vspace{1cm}

\begin{proof}[Proof of Theorem~\ref{thm:4.6} (3)]
This is immediate from part (2) and the fact that $curv X \leq c(K^{+},{\lambda}^{-})$ (see Proposition 1.4 of \cite{Lyt:01}, or \cite{Nik:83}).
\end{proof}

\appendix

\section[Injectivity Radii]{Injectivity Radii}
\label{s:1}

\subsection[Definitions]{Injectivity Radii---Definitions}
\label{ss:1.1}

In a manifold-with-boundary, the usual Riemannian exponential map is not well-defined because geodesics may bifurcate. Nevertheless, one may define a natural notion of conjugate radius and several notions of injectivity radius, all of which are useful. 

\vspace{1cm}

\begin{definition}[\cite{ABB:93}]\label{def:conj}
$M$ is said to have no \emph{conjugate point} along a geodesic $\gamma$ from a point $p\in M$ if the right-hand endpoint map on the space of geodesics from $p$, in the uniform metric, acts homeomorphically from a neighborhood of $\gamma$ onto a neighborhood of its right-hand endpoint. Define the \emph{conjugate radius} $conj_{M}(p)$ of $p$ to be the infimum of distances between $p$ and any distict point $q$ which is conjugate to it. Define $conj(M)=\underset{p\in M}{\inf} \{ conj_{M}(p)\}$.
\end{definition}

This coincides with the usual definition of conjugacy in terms of Jacobi fields, when the space under consideration is a closed manifold.

\begin{lemma}[\cite{ABB:93}, Cor.3, p.711]\label{lem:1.1} 
For $(M,\partial M)$ with $K_{M}\leq K$, and
$K_{\partial M}(U,V)\leq K$ for all those 2-planes $U\wedge V$ for which $II(U,V)<0$,
there exists a lower bound $conj(M)\geq \frac{\pi}{\sqrt{K}}$.
\end{lemma}

\label{def:bigon}
Define a \emph{geodesic bigon} (in a general geodesic metric space $X$) to the union of two distinct minimal geodesic segments having the same length and common initial and terminal endpoints. \\

\begin{figure}[h]
\begin{center}
\begin{overpic}[scale=.6]{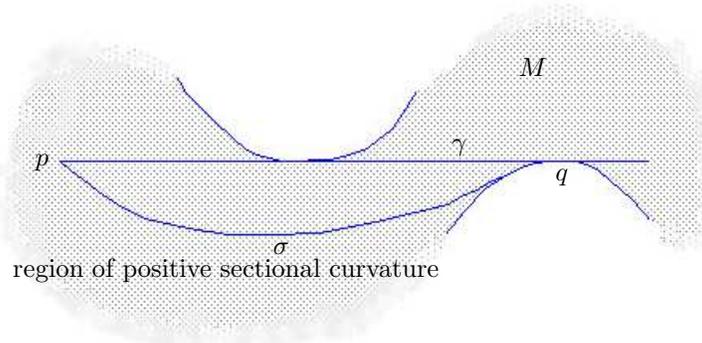}
\put(77,24){$q$}
\put(72,38){$M$}
\put(4,11){region of positive sectional curvature}
\put(7,26){$p$}
\put(39,14){$\sigma$}
\put(63,28){$\gamma$}
\end{overpic}
\caption[minimizing beyond a bigon endpoint]{minimizing beyond a bigon endpoint}\label{fig:1}
\end{center}
\end{figure}

Interestingly enough, in a manifold-with-boundary geodesics may minimize beyond a conjugate point, as well as beyond a nontrivial bigon endpoint, contrary to the case with closed manifolds. The essential reason resides in the fact that tangent vectors to distinct geodesics may form an angle equal to $0$ (as at $q$ in the following example).
See Figure \ref{fig:1}. $\gamma \cup \sigma$ represents a geodesic bigon, yet $\gamma$ may be extended as a minimizer beyond the point $q$.

(Another example was found independently in~\cite{ABB:93b})

\begin{definition}\label{def:i_int}
For a Riemannian manifold-with-boundary $(M,\partial M)$,
and $p\in M$, 
let $i_{int}(p)=\sup\{ r>0 :$ any unit-speed geodesic  $\gamma: \left[0,t_{\gamma}\right]\longrightarrow M$ issuing from $p$ is distance minimizing up to the distance $\min(t_{\gamma},r) \}$, where $t_{\gamma}$ is the first time $\gamma$ intersects $\partial M$ (so $t_{\gamma}=\infty$ if $\gamma\cap \partial M = \emptyset$).
Define $i_{int}(M)=\underset{p\in M}{\inf} \{ i_{int}(p) \}$, the \emph{interior injectivity radius}
\end{definition}

\vspace{0.3cm}

\begin{definition}\label{def:inj}
For a geodesic metric space $X$, and $x\in X$, let $inj(x)=\sup \{ r>0 :$ any minimizing geodesic issuing from $x$ is unique up to distance $r \}$. Define $inj(X)=\underset{x\in X}{\inf} \{ inj(x) \}$,  the \emph{intrinsic injectivity radius of $X$}.
\end{definition}

\vspace{0.3cm}

\begin{definition}\label{def:i_partial}
For a Riemannian manifold-with-boundary $(M,\partial M)$, and $p\in \partial M$, let $i_{\partial}(p)=\sup \{ r>0 :$ any minimizing geodesic $\gamma$ issuing from $p$ normally to $\partial M$ uniquely minimizes distance to $\partial M$, up to distance $r$ (i.e., $\gamma(0)=p$ and $d(\gamma(r),\partial M)=r$) $ \}$. 
Define $i_{\partial}(M)=\underset{p\in \partial M}{\inf} \{ i_{\partial}(p) \}$,  the \emph{boundary injectivity radius} of $(M,\partial M)$.
\end{definition}

\noindent
This is essentially the length of the longest vector for which the normal exponential map of the boundary is nonsingular.

\vspace{0.6cm}

\rmk{
$i_{\partial}(M)=\min \{ Foc(\partial M), \frac{1}{2}L \}$, where $Foc(\partial M)$ denotes the minimum focal distance for the normal exponential map of the boundary, and $L$ represents the length of a shortest segment, meeting $\partial M$ at right angles at both its endpoints. It is known that $Foc(\partial M) \geq \frac{1}{\sqrt{K^{+}}} \arctan \left( \frac{\sqrt{K^{+}}}{{\lambda}^{+}} \right)$, if $K_{M} \leq K^{+}$ and $II_{\partial M} \leq {\lambda}^{+}$. \\\\
}

\noindent
\bf Example: \rm
Euclidean space with a ball of radius $r$ removed, $M=(\mathbb{R}^n\setminus B^{n}(r),g_{std})$, has
\begin{center}
$i_{int}(M)=\infty$, \quad $inj(M)=\pi r$ \quad and $\ \ i_{\partial}(M)=\infty$.\\
\end{center}
More precisely, for any point $p\in M$,
\begin{center}
$i_{int}(p)=\infty$, \quad $inj_{M}(p) = \sqrt{{(R+r)}^2-r^2} + r\left( \frac{\pi}{2} + {\sin}^{-1}\left(\frac{r}{R+r}\right) \right)$\\
\end{center}
where $R=d(p,\partial M)$.

\vspace{0.5cm}

\noindent
\bf Example: \rm
The standard sphere of radius $1$, with a ball of radius $r$ removed,
$M=(S^{n}(1)\setminus B^{n}(r),g_{std})$, has, for $0<r<\pi$

$i_{int}(M)=\begin{cases}
\pi & r\leq \pi/2\\
\infty & r> \pi/2
\end{cases}$, $\ \ \ \ $
$inj(M)=\begin{cases}
\pi \sin(r) & r\leq \pi/2\\
\infty & r> \pi/2
\end{cases}$, $\ \ \ \ $
and $\ \ i_{\partial}(M)=\pi-r$.\\\\\\

It is immediate from the definitions that for any $p\in M$,
\begin{align}\label{ineq:1.13}
i_{int}(p) &\geq inj_{M}(p)
\intertext{and}inj_{M}(p) &\geq \min \{ i_{int}(p), R \}
\end{align}
whenever $R=d(p,\partial M)$. Thus, away from the boundary, $i_{int}(p)$ and $inj_{M}(p)$ are comparable quantities.
Relation (\ref{ineq:1.13}) implies, again by definition, that for any manifold-with-boundary $M$,
\begin{align*}
i_{int}(M)\geq inj(M)
\end{align*}
always holds.
Moreover,
for a closed manifold $M$, or a manifold with locally convex boundary, one has equality:
$i_{int}(M)= inj(M)$.\\

\vspace{1cm}

Lastly, there are two more invariants that needs to be introduced.

\vspace{0.5cm}

\begin{definition}\label{def:inrad}
For a manifold-with-boundary M, the \emph{inradius} is $inrad(M)=\underset{p\in M}{\sup} \{ inrad(p)\}$, where $inrad(p)=\sup \{ r>0 : B(p,r)\subseteq \text{int}M \} = d(p,\partial M)$.
\end{definition}

This quantity, giving the radius of the largest metric ball which fits entirely in the interior of $M$, is defined even for manifolds with nonsmooth boundary.

\begin{definition}\label{def:reach}
For a manifold $M$ with (not necessarily smooth) boundary, define the \emph{reach} of a boundary point $p\in \partial M$ as
$reach(p) = \sup \{ r \geq 0$ : there exists an  $x\in int M$ such that  $|px|_{M}=r$ and $[px]_{M}$ is the unique minimizing segment realizing the distance to  $\partial M \}$
\end{definition}

\rmk{
For manifolds whose boundary is at least $C^2$, one has
$i_{\partial}(p) = reach(p)$
}

\vspace{0.4cm}

\begin{lemma}\label{lem:4.1}
$inrad(M) = \underset{p\in \partial M}{\sup} reach(p)$.
\end{lemma}

As a corollary, for any sequence $(M_i,\partial M_i)$ of manifolds-with-boundary, $inrad(M_i)\longrightarrow 0$ if and only if $reach(p_i) \longrightarrow 0$ for all $p_i\in \partial M_i$.

\begin{proof}
 ( $\geq$ ): Suppose $reach(p)\geq i_0 > 0$ for some $p\in \partial M$. Let $[px]_{M}$ be a minimizing segment realizing $reach(p)$. Then $inrad(M)\geq inrad(x)\geq i_0$,
For if $\partial M$ intersected the interior of $B(x,i_0)$,
with $d(x,\partial M)=d(x,q)<i_0$ say, for some point $q\in \partial M$, $q\neq p$, then the geodesic
$[px]_{M}$ would not minimize the distance to $\partial M$ after all. This contradicts the definition of $reach(p)$ and choice of $x$.\\

\noindent
( $\leq$ ):
We need to produce a sequence of points $p_{k} \in \partial M$ for which $reach(p_{k}) \longrightarrow inrad(M)$. 

By definition of $inrad(M)$, there exists a sequence of points $x_{k} \in M$ such that
$inrad(x_{k}) \longrightarrow inrad(M)$.
It can be assumed that the $x_{k}$ are chosen as the centers of interior balls
whose closure intersects $\partial M$ in at least
one point, say $p_{k}$.
Then we claim that
$reach(p_{k})\geq inrad(x_{k})$.

To see this, consider a minimal geodesic segment $[x_{k}p_{k}]$.
Next, let $y_{k}$ belong to the interior of $[x_{k}p_{k}]$. Note that $y_{k}$ belongs to the interior of the manifold, since it is contained in the interior ball $B(x_{k}, inrad(x_{k}); M)$. To get the stated lower bound on $reach(p_{k})$, it only remains to show that there can be no geodesic segment (say $\sigma$) from $y_{k}$ to $\partial M$, distinct from $[y_{k}p_{k}] \subset [x_{k}p_{k}]$, which has length equal or smaller than $|y_{k}p_{k}|$.

Suppose $|\sigma| \leq |y_{k}p_{k}|$. Then 
$[x_{k}y_{k}]\cup \sigma$, having length no greater than $|x_{k}p_{k}|$, would be a minimal segment from $x_{k}$ to $\partial M$ distinct from the original segment $[x_{k}p_{k}]$. This contradicts the fact that there is no branching in the interior of $M$ (i.e., there exists a convex ball centered at $y_{k}$ and entirely contained in the interior of $M$, which has curvature bounded below by some finite constant).
\end{proof}

\subsection[InterRelations]{Injectivity Radii---InterRelations}
\label{ss:1.3}

The main technical result in this section is Proposition~\ref{prop:1.2}, which relates the injectivity radii introduced in $\S$\ref{ss:1.1} and gives an exponential decay rate for the intrinsic injectivity radius of a manifold-with-boundary. 
The estimates may be viewed as giving more easily verifiable conditions under which the intrinsic injectivity radius is bounded below.
Interesting in their own right, the estimates are used in 
Theorem~\ref{thm:4.2}(ii) in $\S$\ref{s:2}.\\

The boundary may be considered a type of generalized point, namely, what one would get by puncturing a closed manifold.
In this sense, $inj(\partial M)$ together with $i_{\partial}(M)$ (concretely, the quantity
$\min\{ inj(\partial M), i_{\partial}(M) \}$)
function as the injectivity radius of this generalized point. If $M$ is $n$-dimensional, $inj(\partial M)$ accounts for $n-1$ directions, and $i_{\partial}(M)$ accounts for $1$ direction.\\

With this observation, the following proposition may be anticipated.
To state it, let $l_{M}(p)$ denote the length of the shortest nontrivial (not-necessarily smoothly closed) geodesic bigon based at $p$, the sides of which are allowed to contact the boundary in their interior or endpoints. Recall that by definition (page~\pageref{def:bigon}), a geodesic bigon is the union of two distinct minimal segments having common initial and terminal endpoints.
By results of \cite{ABB:93}, a given manifold-with-boundary $M$ is an Alexandrov space of curvature bounded above, so is locally $CAT(k)$ for some $k$.
This implies
in particular that $l_{M}(p) > 0$ is positive for any point $p \in M$.
It turns out that a geodesic bigon which realizes $l_{M}(p)$ will automatically be $C^1$-smoothly closed at its endpoint.\\

\begin{proposition}[\cite{Wo:06}]\label{prop:1.2}
For any complete Riemannian manifold $M$ with boundary, and any point $p\in M$,
\begin{itemize}
\item[(i)]
$inj_M(p)\geq \min\{conj_M(p),\frac{1}{2}\mathit{l}_{M}(p)\}$.
\item[(ii)] There exists a constant $c$ such that if $p\in M[0,i_{\partial}/2]$,

$\mathit{l}_{M}(p)\geq c(n, K^{\pm},{\lambda}^{\pm},inj(\partial M),i_{\partial})>0$.
\item[(iii)] There exists a constant $c$ such that for all $p\in M$,

$inj_M(p)\geq c(n,K^{\pm},{\lambda}^{\pm},inj(\partial M),i_{\partial},d)>0$.
\item[(iv)] There exists a constant $c$ such that

$inj(\partial M) \geq c(K^{+},{\lambda}^{\pm},inj(M))>0$.
\end{itemize}
\end{proposition}

\vspace{0.5cm}

\rmk{
Unlike the situation in the closed manifold case, one can have inequality in part (i) of Proposition~\ref{prop:1.2}. See the beginning of section \ref{ss:1.1} for an example. However, one is invariably interested in lower bounds to injectivity radius, as opposed to upper bounds.
}

\vspace{0.2cm}

\rmk{
The estimates (i) and (iv) in Proposition~\ref{prop:1.2} are sharp. The constants in the lower estimate (iii) may be sharpened, but exponential decay is inevitable, as examples show.\\
}

For part (ii), let $\gamma$ denote a geodesic realizing ${\mathit{l}}_{M}(p)$, where $p \in M[0,i_{\partial}/2]$. 
The idea is to project $\gamma$ to $\partial M$, obtaining a curve whose geodesic curvature in $\partial M$ is uniformly bounded, use arc/chord comparison in $\partial M$ (see $\S$\ref{s:A2}, Theorem~\ref{thm:arc}(c) ) to extract a lower bound for the length of the projected curve, and then use the Lipschitzness of the original projection map to obtain a lower bound for the length of $\gamma$.

\vspace{0.3cm}

For the third part (iii) of Proposition~\ref{prop:1.2}, note that it is automatic from (i), (ii) and Lemma~\ref{lem:1.1} that the intrinsic injectivity radius $inj_{M}(p)$ of $M$ at $p$ is bounded below, for all $p\in M[0,i_{\partial}/2]$. The point is that (iii) provides a lower bound for all $p\in M$.

\vspace{0.5cm}

   \section{(Alexandrov) Extension}\label{ss:1.2}

Beginning with a Riemannian manifold-with-boundary $(M,\partial M)$ one may manufacture a collar, which, when isometrically glued to the boundary, yields an Alexandrov space of curvature bounded below. Outside the gluing locus $\partial M$, the resulting extension $\widetilde{M}$ is $C^{\infty}$ smooth. Actually, $\widetilde{M}$ is a $C^{0}$ Riemannian manifold with a $C^{1,\alpha}$ differentiable manifold structure. The following lemma constructs the collar, and the proposition after it constructs $\widetilde{M}$. The rest of the section details consequences. See~\cite{Wo:07} for more details of the results summarized in this section.

\subsection{Construction}\label{sss:1.2.1}

\begin{lemma}\label{lem:2.1}
Suppose
$M$ is any manifold-with-boundary having
$K_{M}\geq K^{-}$ and ${\lambda}^{-}\leq II_{\partial M}\leq {\lambda}^{+}$.
Then for any $t_{0}>0$, there exists an intrinsic metric on ${\partial M}\times [0,t_0]$, such that $II_{{\partial M}\times \{0 \}}\geq |\min \{0,{\lambda}^{-}\}|$ and $II_{{\partial M}\times \{t_0 \}}\equiv 0$ and the sectional curvature of ${\partial M}\times [0,t_0]$ is bounded below by a constant $c(K^{-},{\lambda}^{\pm},t_0)$.
\end{lemma}

\begin{figure}[!h]
\begin{center}
\begin{overpic}[scale=.8]{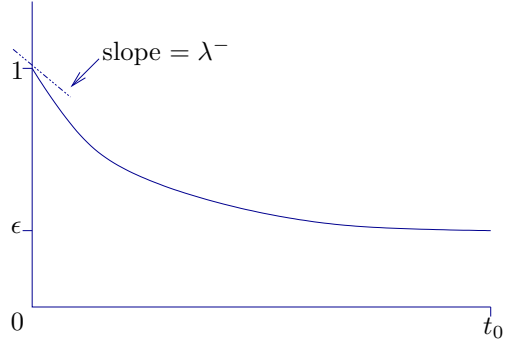}
\put(0,-2){$0$}
\put(19,54){slope $= {\lambda}^{-}$}
\put(0,50){$1$}
\put(0,18){$\epsilon$}
\put(98,-3){$t_0$}
\end{overpic}
\caption[warping function $\phi$]{warping function $\phi$}
\label{fig:phi}
\end{center}
\end{figure}

\begin{proof}
Let $\overline{\lambda}=\min \{0,{\lambda}^{-}\}$. Fix some $t_0>0$ and $0<\epsilon<1$.
For some $K=K(\lambda,\epsilon,t_0)\in \mathbb{R}$, there exists a $C^{\infty}$ monotone non-increasing function $\phi(t)$, defined on $[0,t_0]$, which satisfies
\begin{align*}
{\phi}''+K\phi &\leq 0\\
\phi(0) &= 1\\
-\infty<{\phi}'(0) &\leq \overline{\lambda}\\
\phi(t_0) &= \epsilon\\
{\phi}'(t_0) &= 0.
\end{align*}
See Figure~\ref{fig:phi}.

Consider the warped-product metric on $\partial M \times [0,t_0]$ given by $g_{1}(x,t)={dt}^2+{\phi(t)}^{2}g_{\partial M}(x)$.
\end{proof}

\vspace{0.3cm}

In the lemma, $t_0$ and $\epsilon$ are independent free parameters which may be chosen according to one's purpose. The optimal (i.e., the greatest) lower bound $K_0$ achievable for some $\phi$ (satisfying the above requirements) decreases to $-\infty$ as $t_0$ decreases to $0$ (when $\epsilon$ fixed). It also decreases to $-\infty$ as $\epsilon$ tends to $1$ (when $t_0$ fixed), provided ${\lambda}^{-}<0$ is fixed too. \\

\vspace{0.5cm}
\noindent
\bf An Explicit warping function $\phi$: \rm\\
Here is an explicit construction of a warping function $\phi$ which satisfies the condition of the lemma.

Assume that
\begin{align*}
0 &< \epsilon < 1\\
0 &< t_0\\
\text{ and } \quad \overline{\lambda} &\leq 0.
\end{align*}

It may be further assumed that $\overline{\lambda} < 0$, since otherwise ${\lambda}^{-} \geq 0$, and then the boundary, being locally convex, would not require an extension (in this case one could just take $\phi \equiv 1$). \\

For $0 \leq t < t_0$, define
\begin{align}\label{functionphi}
\phi(t)=(1-\epsilon) \exp{\left[ \frac{\overline{\lambda}t_{0}^{2}}{1-\epsilon} \left( \frac{1}{t_{0}-t} - \frac{1}{t_{0}} \right) \right] } + \epsilon.
\end{align}

Extend $\phi$ to be defined on $[0,t_{0}]$ by requiring continuity of $\phi$ and all its derivatives:
$\phi(t_0):=\underset{t \uparrow t_{0}}{\lim} \ \phi(t)$, 
$\phi'(t_0):=\underset{t \uparrow t_{0}}{\lim} \ \phi'(t)$, etc.

If $\frac{|\overline{\lambda}|t_{0}}{1-\epsilon} > \frac{6}{3-\sqrt{3}}$
then
the radial curvatures
are bounded below by
\begin{align}\label{functineq:4}
-\frac{\phi''(t)}{\phi(t)}
&\geq 
-\frac{1}{\epsilon} \max \{ 0,
\left[
2\frac{\overline{\lambda}}{t_{0}}
+ \frac{{(\overline{\lambda})}^{2}}{1-\epsilon }
\right]
\}. 
\end{align}

If $\frac{|\overline{\lambda}|t_{0}}{1-\epsilon} > 2$ then
the tangential sectional curvatures are bounded below by
\begin{align}\label{functineq:6}
\frac{1}{{\phi}^2(t)} \left[K_{\partial M}^{-} - {|\phi'(t)|}^2 \right]
&\geq  \min \{
0, \frac{K_{\partial M}^{-}}{{\phi}_{\min}^2} \} - \frac{ \max {|\phi'(t)|}^2 }{\min {\phi}^{2}(t)} \notag\\
&\geq \frac{1}{{\epsilon}^2}\left[ \min \{ K_{\partial M}^{-},0 \} -   {|\overline{\lambda}|}^2  \right].
\end{align}

If in an orthonormal frame the sectional curvatures are bounded (from below) on all coordinate two-planes, then the sectional curvatures are bounded (from below) on arbitrary two-planes.
So (\ref{functineq:4}) and (\ref{functineq:6}) together prove that the sectional curvatures of $({\partial M}\times [0,t_0], g_{1})$ are bounded below by a constant $c(K^{-},{\lambda}^{\pm},t_0)$, as stated in
Lemma~\ref{lem:2.1}.

\vspace{1.7cm}

\begin{proposition}[Construction of Extension]\label{prop:2.1}
Fix $n\geq 2$ and $K^{-},\lambda^{\pm}\in \mathbb{R}$. For any\\
$M\in \mathcal{M}(n,K^{-},{\lambda}^{\pm})$ $\equiv$ $\{M \ \text{ Riemannian n-manifold } : K^{-}\leq K_{M}, \ {\lambda}^{-}\leq II_{\partial M}\leq {\lambda}^{+} \}$ there exists an isometric, uniform extension $\widetilde{M}$ of $M$ which is an Alexandrov space of curvature bounded below by a constant $k=k(K^{-},{\lambda}^{\pm})$. The extension is uniform in the sense that the distance in $\widetilde{M}$ between $\partial M$ and $\partial \widetilde{M}$ is no smaller than a constant which may be chosen arbitrarily.
\end{proposition}

\vspace{0.5cm}

\vspace{0.3cm}

\rmk{
The upper bound ${\lambda}^{+}$ is not needed in
Proposition~\ref{prop:2.1}
when dim$M=2$, since in this situation there are no tangential two-planes of which to speak of curvature.\\
}

\vspace{0.7cm}

\noindent
\bf Properties of the extension $\widetilde{M}$: \rm

\begin{align}\label{property:}
i_{\partial}(\widetilde{M}) &\geq t_0 \tag{i}\\
d(\widetilde{M}) &\leq d(M) + 2t_0 \tag{ii}\\
|xy|_{M} &\leq \frac{1}{\epsilon} |xy|_{\widetilde{M}} \text{ for all } x,y\in M.
\quad \text{In particular, } d(M) \leq \frac{1}{\epsilon} d(\widetilde{M})
\tag{iii}
\\
\partial M &= \frac{1}{\epsilon} \partial \widetilde{M}.
\quad \text{In particular, } d(\partial M) = \frac{1}{\epsilon} d(\partial \widetilde{M})
\tag{iv}
\end{align}

\rmk{
Under only the hypotheses of Proposition~\ref{prop:2.1}, neither $M$ nor $C_{M}$ need be a locally convex subset of $\widetilde{M}$. For instance, to see that $C_{M}$ need not be locally convex in $\widetilde{M}$, take $M$ to be the result of cutting lengthwise (through the apex) a rounded-off cone with small cone-angle, so that the resulting boundary of $M$ near the apex is totally geodesic, but elsewhere has some concavity. \\
}

\vspace{1cm}

There is a projection map from $\widetilde{M}$ to $M$, whose Lipschitz continuity (as well as Lipschitz constant)
is important in
applications.

\begin{lemma}[\cite{Wo:07}]\label{lem:2.4}
Let $(M,\partial M)$ be fixed.
Define a map $\pi : \widetilde{M} \longrightarrow M$ by
\[
\pi(x) = \begin{cases}
\text{ x }& \text{ if } x\in M\\
\text{ orthogonal projection of } x \text { onto } \partial C_{M}=\partial M & \text{ if } x\in C_{M}
\end{cases}
\]

Then for all $x,y\in \widetilde{M}$, $|\pi(x)\pi(y)|_{M} \leq \frac{1}{\epsilon}|xy|_{\widetilde{M}}$.
\end{lemma}

\vspace{0.5cm}

\vspace{0.9cm}

Even though $d_{H}^{\widetilde{M}}(M,\widetilde{M})=t_{0}$, it need not necessarily hold that $d_{GH}(M,\widetilde{M})\leq t_{0}$. (This would be true if $M$ were convexly embedded in $\widetilde{M}$, in the sense that $|xy|_{M} = |xy|_{\widetilde{M}}$ for all $x,y\in M$.) However, if $d(M)\leq d$, the map $\pi$ defined above is a $c(t_{0},\epsilon,d)$-Hausdorff approximation for a computable constant $c$.\\

\begin{lemma}\label{lem:2.5}
$\pi : \widetilde{M} \longrightarrow M$ is a $\max \{ 2t_{0}, (\frac{1}{\epsilon} - 1)(d+2t_{0}) \}$-Hausdorff approximation.
\end{lemma}

\begin{proof}
For any $x,y\in \widetilde{M}$, 
\begin{align*}
|xy|_{\widetilde{M}}
&\leq |x\pi(x)|_{\widetilde{M}} + |\pi(x)\pi(y)|_{\widetilde{M}} + |\pi(y)y|_{\widetilde{M}}\\
&\leq 2t_{0} + |\pi(x)\pi(y)|_{\widetilde{M}}\\
&\leq 2t_{0} + |\pi(x)\pi(y)|_{M}
\end{align*}
Hence
\begin{align*}
|xy|_{\widetilde{M}} - |\pi(x)\pi(y)|_{M} \leq 2t_{0}.
\end{align*}
Conversely,
by Lemma~\ref{lem:2.4},
\begin{align*}
|\pi(x)\pi(y)|_{M} - |xy|_{\widetilde{M}}
&\leq \frac{1}{\epsilon}|xy|_{\widetilde{M}} - |xy|_{\widetilde{M}}\\
&\leq \left( \frac{1}{\epsilon} - 1 \right) d(\widetilde{M})\\
&\leq \left( \frac{1}{\epsilon} - 1 \right) (d+2t_{0}) \qedhere
\end{align*}
\end{proof}

\vspace{1cm}

\subsection{Corollaries of Construction}\label{sss:1.2.2}

\begin{corollary}[\cite{Wo:07}]\label{cor:2.2}
If a sequence $M_{i}\in \mathcal{M}(n,K^{-},{\lambda}^{\pm},d)$ GH-converges to a metric space X (necessarily compact and geodesic) then $dim_{\mathcal{H}}X\leq n$.
\end{corollary}

\begin{proposition}[\cite{Wo:07}]\label{lem:2.8}
For any $M\in \mathcal{M}(n,K^{-},{\lambda}^{\pm},d)$, the intrinsic diameter of any boundary component of $M$ is uniformly bounded above by
$\ d(\partial M) \leq D(n,K^{-},{\lambda}^{\pm},d)$.
\end{proposition}

\begin{proposition}[\cite{Wo:07}]\label{prop:2.8}
Suppose $\{ M_{i} \}$ is a sequence of $n$-dimensional Riemannian manifolds-with-boundary such that $K^{-} \leq K_{M_i}, \ |II| \leq \lambda, \ d(M_i) \leq d$,
and $M_i$ GH-converge to a limit space.
Then there is a homotopy equivalence
\[
\underset{GH}{\lim} M_i \simeq \underset{GH}{\lim} \widetilde{M_i} \, ,
\]
where $\widetilde{M_i}$ are the (Alexandrov) extensions of $M_i$, as in
Proposition~\ref{prop:2.1}.
\end{proposition}

\newpage

\section{Arc/chord Comparison}
\label{s:A2}

Arc/chord comparison relates geometric quantities associated with a curve (in a space of curvature bounded above by $K$) having a particular geodesic curvature to corresponding quantities of a model curve in a space form of curvature $K$. Such quantities include arclength, chordlength, base-angles, and width. See~\cite{AB:96} for more discussion and proofs of the results summarized below. \\

The setting for arc/chord comparison is in $CAT(k)$-domains.

A metric space is called $CAT(k)$, for some $k\in \mathbb{R}$, if any two points are joined by a minimal geodesic segment, and the space has Alexandrov curvature bounded from above, i.e., each triangle of perimeter $< \frac{2\pi}{\sqrt{k}}$ satisfies the triangle comparison condition.
Equivalently, it is $CAT(k)$ if minimizers of length $< \frac{\pi}{\sqrt{k}}$ exist, are unique, and vary continuously with their endpoints (\cite{AB:96}, p.68).\label{def:CAT_k}

Spaces of curvature bounded above are locally $CAT(k)$, where the value of $k$ might vary from one region to another.
A typical example of a $CAT(k)$ space is a metric ball $B(x,r; X)$ of radius $r$, in a space $X$ with curvature bounded above by $k$, where $r$ is less than the so-called $CAT_{k}$-radius, which may be estimated from below
by
\begin{align*}
CAT_{k}Rad(x) = \min \{ \frac{\pi}{2\sqrt{k}}, inj(x) \}
\end{align*}
(See \cite{AB:96}, Theorem 4.3, p.78 for a proof of this estimate.)

By definition, $CAT_{k}Rad(X) := \underset{x\in X} \inf CAT_{k}Rad(x)$.\\\\

Let $U$ be a $CAT(K)$-domain.\\

\noindent
\bf Arclength, chordlength \rm

\begin{theorem}[\protect{\cite[Corollary 1.2]{AB:96}}]\label{thm:arc}
Let $\gamma$ be a curve in $U$ with pointwise arc/chord curvature $\kappa$, where $\kappa \leq k$, and assume the sum of its arclength $s$ and chordlength $r$ is less than $2\pi / \sqrt{K}$. Let $\sigma$ be the complete $k$-curve in the model space $M_{K}^{2}$.

(a) If $r$ is less than the diameter of $\sigma$, then $\gamma$ is either no longer than the minor arc of $\sigma$ with chordlength $r$ or no shorter than the major arc of $\sigma$ with chordlength $r$. The upper bound on length holds if the diameter of $\gamma$ is less than the diameter of $\sigma$.

(c) If $\gamma$ is closed and nonconstant (not necessarily closing smoothly), then $\gamma$ is no shorter than a complete $k$-curve (necessarily a circle) in $M_{K}^{2}$.
\end{theorem}

\vspace{0.4cm}

Let $s = s(r,k,K)$ be the arclength of a minor $k$-arc of chordlength $r$ in $M_{K}^{2}$. Then the function $s$ is an analytic function of $(r,k,K)$ with series in powers of $r$ in all cases:
\[
s = r + \frac{1}{24}k^2r^3 + \frac{9k^{4}+8k^{2}K}{1920}r^5 + \ldots
\]

In particular, for instance, if $K=0$ then
$s = \frac{2}{k}{\sin}^{-1}\left( \frac{kr}{2} \right)$.\\

By increasing the coefficient $\frac{1}{24}k^2$ of $r^3$ to a certain computable constant $C(k,K)$, it follows that 
$s \leq r + C(k,K)r^3$ for all $r\leq 1$.\\

\noindent
\bf Angle, width \rm\\

Changing notation, let
$\gamma: [0,1] \longrightarrow U$ be any curve,
$\sigma: [0,1] \longrightarrow U$ be the minimizing geodesic segment joining its endpoints.

Define the width of $\gamma$ as $W(\gamma):=\underset{s\in [0,1]}{\max} \underset{t \in [0,1]}{\min} d(\gamma(s),\sigma(t))$.

This width is the radius of the smallest tubular neighborhood about the chord $\sigma$ that contains $\gamma$.\\

\begin{theorem}[\protect{\cite[Theorem 6.1]{AB:96}}]
If $\gamma$ is a curve in $U$ with $\kappa \leq k$ and length no more than half a complete $k$-curve in $M_{K}^{2}$, then the width and base-angles of $\gamma$ are no greater than they are for a $k$-arc $\overline{\gamma}$ in $M_{K}^{2}$ of the same length.
\end{theorem}

In symbols,
\begin{align*}
W(\gamma) &\leq W(\overline{\gamma}) = kr^{2}/8 + O(r^3)\\
\text{ and }
{\angle}(\gamma,[\gamma(0)\gamma(1)]) &\leq {\angle}(\overline{\gamma},[\overline{\gamma(0)}\overline{\gamma(1)}]) = kr/2 + O(r^3),
\end{align*}
where $r=|\gamma(0)\gamma(1)|$ denotes chordlength.

\vspace{0.8cm}

\section{Gromov-Hausdorff Convergence}
\label{s:A3}

This appendix section gives background on Gromov-Hausdorff convergence. In particular, it details
two
functorial properties of maps: one for surjective, Lipschitz maps, and another for Lipschitz homotopy equivalences.
For supplementary references and additional background on Gromov-Hausdorff convergence, the reader may consult~\cite{F:90} or \cite{Pe:93}.\\

Let $Z$ be a metric space.
The \emph{Hausdorff distance $d_{H}^{Z}(X,Y)$} between two subsets $X,Y \subseteq Z$ is defined to be
$d_{H}^{Z}(X,Y) := \inf \{ \epsilon > 0 : B(X,\epsilon)\supseteq Y, \, B(Y,\epsilon)\supseteq X   \}$
where $B(X,\epsilon)=\{z \in Z : d(z,X)< \epsilon \}$ denotes the metric $Z$-ball about $X$, of radius $\epsilon$.

Suppose $X$, $Y$ are metric spaces.
The \emph{Gromov-Hausdorff distance} between $X$ and $Y$ is
\begin{align*}
d_{GH}(X,Y) = \inf \{ \,
 d_{H}^{Z}(i_{X}(X),i_{Y}(Y)) \}
\end{align*}
where the infimum is taken over all metric spaces $Z$
and all distance-preserving embeddings  $i_{X}: X \hookrightarrow Z$ and  $i_{Y} : Y \hookrightarrow Z$.

One says that metric spaces $X_i$ converge to $X$, and writes $X_{i} \overset{GH}{\longrightarrow} X$, if $d_{GH}(X,X_{i}) \longrightarrow 0$ as $i \longrightarrow \infty$.
In practice one usually uses the following formulation to verify that a convergence occurs.\\

\begin{definition}
An \emph{$\epsilon$-Hausdorff approximation} $f : X \longrightarrow Y$ is a (not necessarily continuous) map such that\\
(i) $f(X)$ is an $\epsilon$-net in $Y$, i.e.,
\hspace{0.5cm}     $B(f(X),\epsilon; Y) = Y$ and\\
(ii) $f$ is an $\epsilon$-almost isometry, i.e.,
\hspace{0.2cm}   $|d_{Y}(f(x_1),f(x_2)) - d_{X}(x_1,x_2)| \leq \epsilon$ $ \ \forall x_1,x_2 \in X$
\end{definition}

\vspace{0.5cm}

\noindent
\bf Fact: \rm $d_{GH}(X,Y) \leq 3\epsilon$ if there exists an $\epsilon$-Hausdorff approximation $f : X \longrightarrow Y$.\\

\noindent
\bf Fact: \rm For any metric spaces $X$ and $Y$, 
$\frac{1}{2}|d(X)-d(Y)| \leq d_{GH}(X,Y)$.

\vspace{0.5cm}

One has the notion of convergence of points.
\begin{definition}
If $X_i \overset{GH}{\longrightarrow} X$ via $\epsilon_i$-Hausdorff approximations $f_{i} : X_i \longrightarrow X$, then one says points $x_i \in X_i$ converge to a point $x \in X$ ($x_{i} \mapsto x$) if $d(f_{i}(x_i),x) \longrightarrow 0$.
\end{definition}

This permits one to define convergence of maps.

\begin{definition}
If $f_i : X_i \longrightarrow Y_i$ are maps, $X_i \overset{GH}{\longrightarrow} X$ and $Y_i \overset{GH}{\longrightarrow} Y$, then \emph{$f_i$ converge to a map $f : X \longrightarrow Y$} 
if $f_{i}(x_i) \mapsto f(x)$ whenever $X_i \ni x_i \mapsto x \in X$.
\end{definition}

\vspace{1cm}

For the proof of the following two propositions, see \cite{Wo:06}.

\begin{proposition}\label{prop:8.1}
Suppose $X_i$ and $Y_i$ are metric spaces with $X_i \overset{GH}{\longrightarrow} X$, $Y_i \overset{GH}{\longrightarrow} Y$, and with $X$, $Y$ compact. Suppose that for all $i$ there exists $L$-Lipschitz maps ${\psi}_{i} : Y_i \longrightarrow X_i$. Then there exists an $L$-Lipschitz map $\Psi : Y \longrightarrow X$. If the ${\psi}_i$ are in addition surjective, then the limit map $\Psi$ is also surjective.
\end{proposition}

\begin{proposition}\label{prop:8.2}

Suppose $X_j, Y_j$ are complete metric spaces
and 
$X_j \overset{GH}{\longrightarrow} X$, \,
$Y_j \overset{GH}{\longrightarrow} Y$
with $X$, $Y$ compact.

Suppose that for each $j$ there exist continuous
maps
$f_{j} : X_{j} \longrightarrow Y_{j}$, \,
$g_{j} : Y_{j} \longrightarrow X_{j}$, \,
$f : X \longrightarrow Y$ and
$g : Y \longrightarrow X$
with
\begin{align*}
f_{j} &\longrightarrow f \quad \text{ in the sense that }
f_{j}(x_j) \mapsto f(x)
\text{ whenever $x_{j} \mapsto x\in X$}\\
g_{j} &\longrightarrow g \quad \text{ in the sense that }
g_{j}(y_j) \mapsto g(y)
\text{ whenever $y_{j} \mapsto y\in Y$}
\end{align*}

Suppose that for each $j$, there exist maps $H_{j} : X_{j} \times I \longrightarrow X_{j}$ with
\begin{align*}
H_{j}(x,0) &= g_{j} \circ f_{j} (x)\\
H_{j}(x,1) &= id_{X_j}(x) = x\\
H_{j}(x,t) & \text{ globally Lipschitz in $x,t$, uniformly in $j$ }\\
& \text{(where $X_{j} \times I$ equipped with the direct product metric)}
\end{align*}

Suppose also for each $j$, there exist maps ${\overline{H}}_{j} : Y_{j} \times I \longrightarrow Y_{j}$ with
\begin{align*}
{\overline{H}}_{j}(x,0) &= f_{j} \circ g_{j} (x)\\
{\overline{H}}_{j}(x,1) &= id_{Y_j}(x) = x\\
{\overline{H}}_{j}(x,t) & \text{ globally Lipschitz in $x,t$, uniformly in $j$ }\\
& \text{(where $Y_{j} \times I$ equipped with the direct product metric)}
\end{align*}

Then
$X$ and $Y$ are homotopy equivalent (via a Lipschitz homotopy equivalence).
\end{proposition}
\begin{align*}
\xymatrix{\ar @{} [dr] |{}
X_j \ar[d]^{\, }_{GH} \ar@<1ex>[r]^{f_j}_{\simeq} & Y_j \ar[d]^{GH}_{\, } \ar@<1ex>[l]^{g_j}\\
X  \ar@{.>}@<1ex>[r]^{f}_{}  & Y \ar@{.>}@<1ex>[l]^{g}_{\bf \therefore\rm \simeq}
}
\end{align*}

\vspace{0.7cm}

\begin{lemma}\label{lem:4.4}
Let $X_i$ be compact length spaces such that $X_i \overset{GH}{\longrightarrow} X$, with $X$ compact.
If $A_i \subseteq X_i$ are convex, compact subsets, then they sub-converge (as subsets)
to a convex subset of $X$.
\end{lemma}

\begin{proof}
Let $f_{i} : X_{i} \longrightarrow X$ be $\epsilon_i$-Hausdorff approximations.
Note that $d_{H}^{X}(f_{i}(A_i),\overline{f_{i}(A_i)})=0$, where $\overline{f_{i}(A_i)}$ denotes the metric (Cauchy) completion of $f_{i}(A_i)$.
Since $X$ is compact, $\overline{ f_i(A_i) }$ is compact. 
Blaschke's theorem (see, e.g., \cite{BBI:01} p.253) implies that for some subsequence the Hausdorff limit $\underset{H}{\lim} \overline{f_i(A_i)}$ in $X$ exists and is compact.
So $A:=\overline{ \underset{H}{\lim} f_{i}(A_i) } = \underset{H}{\lim} \overline{f_i(A_i)}$ exists and is compact. 
In particular, $A$ is metrically complete.

Now convexity of a subspace $A\subseteq (X,d)$ means by definition that the restriction $d|_{A}$ is strictly intrinsic, or equivalently, that for all 
$x,y\in A$,
there exists a shortest $X$-path from $x$ to $y$ contained entirely in $A$. Under the assumption of completeness of $A$, this is implied by the existence of ($X$-)midpoints lying in the set.

Let $x,y \in A$. By definition of $A$, we can choose from $A_i$ approximating points $x_i \mapsto x$ and $y_i \mapsto y$, so $d(f_i(x_i),x)\longrightarrow 0$ and $d(f_i(y_i),y)\longrightarrow 0$ as $i\longrightarrow \infty$. The $X_i$ have midpoints because they are length spaces. Let $z_i:=$midpoint of a shortest segment $[x_{i}y_{i}]_{X_i}$. Then $z_i\in A_i$ by convexity of $A_i$ in $X_i$.
Since $f_i$ are $\epsilon_i$-Hausdorff approximations,
\begin{align*}
| |f_i(x_i)f_i(y_i)|_{X}-|x_{i}y_{i}|_{X_i} | &< \epsilon_i\\
| |f_i(x_i)f_i(z_i)|_{X}-|x_{i}z_{i}|_{X_i} | &< \epsilon_i\\
| |f_i(y_i)f_i(z_i)|_{X}-|y_{i}z_{i}|_{X_i} | &< \epsilon_i.
\end{align*}
By compactness of $X$, the sequence $\{ f_i(z_i) \}$ converges to some $z\in X$. And $z\in A$ by definition of $A$.
Combining the above yields that $z$ is an ($X$-)midpoint of $x$ and $y$.
Therefore the set $A:=\overline{ \underset{H}{\lim} f_{i}(A_i) }$ is convex, since it is complete and has ($X$-)midpoints%
\end{proof}

\vspace{0.3cm}

\rmk{
If the length spaces $X_i$ are $CAT(k)$ for some $k\in \mathbb{R}$, $CAT_{k}Rad(X_i)\geq i_0 > 0$ and $X$ is compact, then a sequence of locally convex, compact subsets $A_i \subseteq X_i$ will sub-converge (as subsets) to a locally convex subset of $X$. This holds because $X$ inherits the $CAT(k)$ property, and in $CAT(k)$ balls, locally convex paths, i.e. geodesics, are convex, i.e., minimizing segments.\\
}

Said somewhat differently, if $X$ is $CAT(k)$ with a lower bound $CAT_{k}Rad(X) \geq r > 0$, and $A \subseteq X$ is a subset for which $d_{A}(x,y) = d_{X}(x,y)$ for any two points $x,y\in A$ sufficiently close, then $A$ is $(0,2,r)$-convex in $X$ (see the definition on p.\pageref{def:cconvex}).\\

More generally, without reference to the existence or properties of geodesics, or any curvature bounds, one has

\begin{lemma}\label{lem:C.3}
Suppose $X_{i} \overset{GH}{\longrightarrow} X$ where $X$ and each $X_i$ are
metric spaces.
Suppose
$A_{i} \subseteq X_{i}$ are
$(C,2,r)$-convex subsets (with metrics $d_{A_i}$), where $C\geq 0$ and $r>0$ are any fixed constants.
Then, if
$A_{i} \longrightarrow A \subseteq X$ as subsets
and $A_{i} \overset{GH}{\longrightarrow} A$,
$A$ is $(C,2,r/2)$-convex.
\end{lemma}

\vspace{0.5cm}

Before beginning the proof, it should be remarked that in general, a subset convergence, such as a Hausdorff convergence, for example, neither implies nor is implied by a Gromov-Hausdorff convergence. This is why the lemma must assume both. On the other hand, there are many natural geometric examples where both types of convergence happen simultaneously. \\

\begin{proof}

Suppose $d_{GH}(X,X_{i}) \leq \epsilon_{i}$ and $d_{GH}(A,A_{i}) \leq \epsilon_{i}$, with $\epsilon_{i} \longrightarrow 0$.

Let $z\in A$, and $x,y\in A$ be any points such that $d_{X}(x,z)\leq r/2$, $d_{X}(y,z)\leq r/2$.
It is required to show that $d_{A}(x,y) \leq d_{X}(x,y) + C d_{X}^{3}(x,y)$.

Since $A_i \longrightarrow A$ as subsets,
there exist $z_{i} \in A_i$ with $z_{i} \mapsto z \in A$
and $x_i, y_i \in A_i$ with $x_i \mapsto x$, $y_i \mapsto y$.
For sufficiently large $i$,
\begin{align*}
d_{X_i}(x_i,z_i) &\leq d_{X}(x,z)+\epsilon_{i} \leq \frac{r}{2} + \epsilon_{i} \leq r
\intertext{and}
d_{X_i}(y_i,z_i) &\leq d_{X}(y,z)+\epsilon_{i} \leq \frac{r}{2} + \epsilon_{i} \leq r.
\end{align*}

By assumption, for all $x_i, y_i \in A_i$ with $X_{i}$-distance $\leq r$ to $z_i$,
\begin{align}\label{ineq:8.9}
d_{A_i}(x_i,y_i) \leq d_{X_i}(x_i,y_i) + C d_{X_i}^{3}(x_i,y_i).
\end{align}

Since $X_{i} \overset{GH}{\longrightarrow} X$, 
$|d_{X_i}(x_i,y_i) - d_{X}(x,y)| \leq \epsilon_{i}$.

Since $A_{i} \overset{GH}{\longrightarrow} A$, 
$|d_{A_i}(x_i,y_i) - d_{A}(x,y)| \leq \epsilon_{i}$.

So for all sufficiently large $i$,
\begin{align*}
d_{A}(x,y)
&\leq d_{A_i}(x_i,y_i) + \epsilon_{i}\\
&\leq  d_{X_i}(x_i,y_i) + C d_{X_i}^{3}(x_i,y_i)  + \epsilon_{i}
\quad \text{ by (\ref{ineq:8.9}), since $d_{X_i}(x_i,z_i), d_{X_i}(y_i,z_i) \leq r$ }\\
&\leq \left( d_{X}(x,y) + \epsilon_{i} \right) + C {\left(  d_{X}(x,y) + \epsilon_{i} \right)}^3 + \epsilon_{i}\\
&= d_{X}(x,y) + C \left ( d_{X}^{3}(x,y) + 3\epsilon_{i}d_{X}^{2}(x,y) + 3{\epsilon_{i}}^{2}d_{X}^{}(x,y) + {\epsilon_{i}}^{3} \right) + \epsilon_{i}
\end{align*}

Hence passing to the limit as $i \longrightarrow \infty$,
\begin{align*}
d_{A}(x,y) \leq d_{X}(x,y) + C d_{X}^{3}(x,y),
\end{align*}
which means by definition that $A$ is $(C,2,r/2)$-convex in $X$.
\end{proof}

\vspace{1cm}

To end the section, we give a commutation relation for limits of warped product metric spaces.
In the present work, it is used only for the proof of Proposition~\ref{prop:8.2} and Theorem~\ref{thm:Top}.
It may be applied to study limits of collars, as produced in Lemma~\ref{lem:2.1}.\\

\begin{proposition}\label{lem:4.7}
Let $X_i$ be geodesic metric spaces. Then
$\underset{GH}{\lim}(X_i \times_{\phi} Y) = (\underset{GH}{\lim}X_i)\times_{\phi} Y$
\noindent
if $\phi : Y \longrightarrow \mathbb{R}$ is continuous and $Y$ compact (whenever the limits exist).
\end{proposition}

\begin{proof}[Sketch of the proof]
Suppose $X = \underset{GH}{\lim}X_i$.

\noindent
The length of a curve $\gamma=(\alpha,\beta) : [0,1] \longrightarrow X {\times}_{\phi} Y$ is defined as
\begin{align*}
L(\gamma):=\underset{\underset{0\leq t_0 \leq \ldots \leq t_{N}=1}{N}}{\sup} {\sum}_{j=1}^{N-1}\sqrt{{\left[ \phi(\beta(t_{j}^{*}))d_{X}(\alpha(t_j),\alpha(t_{j+1}))\right]}^2 + {d_{Y}(\beta(t_j),\beta(t_{j+1})) }^2}
\end{align*}
where 
$t_{j}^{*} \in [t_{j},t_{j+1}]$ is an arbitrary evaluation point. The warped product $X {\times}_{\phi} Y$ is defined as the topological space $X {\times} Y$ equipped with the metric induced from the length structure above.

Any two points in $X {\times}_{\phi} Y$ can be joined by a minimizing segment.

Let $f_{i} : X \longrightarrow X_i$ be an $\epsilon_i$-Hausdorff approximation.

Then the map $F_{i} : X {\times}_{\phi} Y \longrightarrow X_i {\times}_{\phi} Y$ defined by $F_{i}(x,y):=(f_{i}(x),y)$, is
a $\left( \underset{y\in Y}{\sup} |\phi(y)| \cdot \epsilon_i \right)$-Hausdorff approximation.
For more details, see~\cite{Wo:06}.

\end{proof}

\ifnum \draft=0
  \backmatter
\fi

\bibliographystyle{siam}
\bibliography{paper}

\vspace{1cm}

\noindent
\sc{Department of Mathematics\\
University of Toronto\\
Toronto, ON, M5S 2E4}

\noindent
\it{E-mail address: }\tt{jawong1@math.toronto.edu}

\end{document}